\documentclass[12pt]{amsart}
\usepackage{amsmath,amssymb,array,longtable, yfonts}
\usepackage[mathscr]{euscript}
\usepackage[all]{xy} 
\usepackage[colorlinks=true,urlcolor=black,citecolor=black,linkcolor=black,%
pdftitle={Seshadri Constants on P1 times P1, and applications to the symplectic packing problem},%
pdfauthor={Chris Dionne and Mike Roth},%
pdfsubject={Algebraic Geometry},%
pdfkeywords={Algebraic geometry, Seshadri Constants, Symplectic Embeddings}]{hyperref}
\usepackage{nicematrix}  
\usepackage{graphicx}
\usepackage{pstricks}    
\usepackage{pst-plot}
\usepackage{multido}
\usepackage{graphicx}
\usepackage{xpatch}  

\usepackage{ifthen}

\newboolean{bourbaki}
\setboolean{bourbaki}{false}

\newcounter{figcount}[section]

\newcounter{eqcount}[theorem]
\renewcommand{\theeqcount}{\thetheorem.\arabic{eqcount}}
\newcounter{remarkcount}[theorem]
\newcommand{\RemNum}{\refstepcounter{remarkcount}(\Euler{\theremarkcount}) } 

\newcommand{\Osh}{\mathcal{O}}          
\newcommand{\nuhat}{\hat{\nu}}          
\newcommand{\eptilde}{\tilde{\ep}}      
\newcommand{\bpf}{\noindent{\em Proof:} \/}  
\newcommand{\epf}{$\Box$}               
\renewcommand{\emptyset}{\varnothing}   
\newcommand{\ep}{\epsilon}                 

\newcommand{\mult}{\operatorname{mult}}
\newcommand{\Vol}{\operatorname{Vol}}    
\newcommand{\NS}{\operatorname{NS}}      
\newcommand{\Id}{\operatorname{Id}}       
\newcommand{\Aut}{\operatorname{Aut}}     
\renewcommand{\geq}{\geqslant}        
\renewcommand{\leq}{\leqslant}        
\newcommand{\st}{\,\,|\,\,}
\renewcommand{\phi}{\varphi}  
\newcommand{\cY}{\mathcal{Y}} 
\newcommand{\cX}{\mathcal{X}} 
\newcommand{\cN}{\mathcal{N}} 
\newcommand{\cC}{\mathcal{C}} 
\newcommand{\slope}{\operatorname{slope}} 
\newcommand{\ord}{\operatorname{ord}}

\ifthenelse{\boolean{bourbaki}}
{
\newcommand{\CC}{\mathbf{C}} 
\newcommand{\FF}{\mathbf{F}} 
\newcommand{\NN}{\mathbf{N}} 
\newcommand{\PP}{\mathbf{P}} 
\newcommand{\QQ}{\mathbf{Q}} 
\newcommand{\RR}{\mathbf{R}} 
\newcommand{\ZZ}{\mathbf{Z}} 

\DeclareSymbolFont{frenchmath}{OT1}{cmr}{m}{n} 

\DeclareMathSymbol{A}{\mathalpha}{frenchmath}{`A}  
\DeclareMathSymbol{B}{\mathalpha}{frenchmath}{`B}  
\DeclareMathSymbol{C}{\mathalpha}{frenchmath}{`C}
\DeclareMathSymbol{D}{\mathalpha}{frenchmath}{`D}
\DeclareMathSymbol{E}{\mathalpha}{frenchmath}{`E}
\DeclareMathSymbol{F}{\mathalpha}{frenchmath}{`F}
\DeclareMathSymbol{G}{\mathalpha}{frenchmath}{`G}
\DeclareMathSymbol{H}{\mathalpha}{frenchmath}{`H}
\DeclareMathSymbol{I}{\mathalpha}{frenchmath}{`I}
\DeclareMathSymbol{J}{\mathalpha}{frenchmath}{`J}
\DeclareMathSymbol{K}{\mathalpha}{frenchmath}{`K}
\DeclareMathSymbol{L}{\mathalpha}{frenchmath}{`L}
\DeclareMathSymbol{M}{\mathalpha}{frenchmath}{`M}
\DeclareMathSymbol{N}{\mathalpha}{frenchmath}{`N}
\DeclareMathSymbol{O}{\mathalpha}{frenchmath}{`O}
\DeclareMathSymbol{P}{\mathalpha}{frenchmath}{`P}
\DeclareMathSymbol{Q}{\mathalpha}{frenchmath}{`Q}
\DeclareMathSymbol{R}{\mathalpha}{frenchmath}{`R}
\DeclareMathSymbol{S}{\mathalpha}{frenchmath}{`S}
\DeclareMathSymbol{T}{\mathalpha}{frenchmath}{`T}
\DeclareMathSymbol{U}{\mathalpha}{frenchmath}{`U}
\DeclareMathSymbol{V}{\mathalpha}{frenchmath}{`V}
\DeclareMathSymbol{W}{\mathalpha}{frenchmath}{`W}
\DeclareMathSymbol{X}{\mathalpha}{frenchmath}{`X}
\DeclareMathSymbol{Y}{\mathalpha}{frenchmath}{`Y}
\DeclareMathSymbol{Z}{\mathalpha}{frenchmath}{`Z}
} 
{ 
\newcommand{\CC}{\mathbb{C}} 
\newcommand{\FF}{\mathbb{F}} 
\newcommand{\NN}{\mathbb{N}} 
\newcommand{\PP}{\mathbb{P}} 
\newcommand{\QQ}{\mathbb{Q}} 
\newcommand{\RR}{\mathbb{R}} 
\newcommand{\ZZ}{\mathbb{Z}} 
}

\providecommand{\binom}[2]{{#1\choose#2}}

\newcommand{\point}{\vspace{3mm}\par\refstepcounter{theorem}\noindent{\bf \thetheorem.} }

\newcommand{\bpoint}[1]{\vspace{3mm}\par\refstepcounter{theorem}\noindent{\bf \thetheorem.} 
  {\bf #1.} }

\renewenvironment{equation}{\medskip\noindent\refstepcounter{eqcount}\makebox[0pt][l]{\rm ({\bf\theeqcount})}\begin{minipage}[b]{\textwidth}$$}{$$\end{minipage}\medskip\noindent}

\newcommand{\Fig}{\refstepcounter{eqcount}{\sc Figure \theeqcount.}} 
\DeclareMathAlphabet{\matheuler}{U}{zeur}{m}{n}
\DeclareMathAlphabet{\matholdstyle}{OT1}{pplj}{m}{n}
\newcommand{\Euler}[1]{\ensuremath{\matheuler{#1}}}  

\newcounter{enumcount}
\newcommand{\PauseEnumerate}{\end{enumerate}\setcounter{enumcount}{\arabic{enumi}}}
\newcommand{\ResumeEnumerate}{\begin{enumerate}\setcounter{enumi}{\theenumcount}}

\newcommand{\Old}[1]{\oldstylenums{#1}}

\newcommand{\AlphaList}{\renewcommand{\labelenumi}{{({\em\alph{enumi}})}}}

\newcommand{\RomanList}{\renewcommand{\labelenumi}{{({\em\roman{enumi}})}}}

\AlphaList  

\newcounter{dummycounter}

\makeatletter
\newcommand*{\da@rightarrow}{\mathchar"0\hexnumber@\symAMSa 4B }
\newcommand*{\da@leftarrow}{\mathchar"0\hexnumber@\symAMSa 4C }
\newcommand*{\xdashrightarrow}[2][]{%
  \mathrel{%
    \mathpalette{\da@xarrow{#1}{#2}{}\da@rightarrow{\,}{}}{}%
   }%
}
\newcommand{\xdashleftarrow}[2][]{%
  \mathrel{%
    \mathpalette{\da@xarrow{#1}{#2}\da@leftarrow{}{}{\,}}{}%
  }%
}
\newcommand*{\da@xarrow}[7]{%
  \sbox0{$\ifx#7\scriptstyle\scriptscriptstyle\else\scriptstyle\fi#5#1#6\m@th$}%
  \sbox2{$\ifx#7\scriptstyle\scriptscriptstyle\else\scriptstyle\fi#5#2#6\m@th$}%
  \sbox4{$#7\dabar@\m@th$}%
  \dimen@=\wd0 %
  \ifdim\wd2 >\dimen@
    \dimen@=\wd2 %
  \fi
  \count@=2 %
  \def\da@bars{\dabar@\dabar@}%
  \@whiledim\count@\wd4<\dimen@\do{%
    \advance\count@\@ne
    \expandafter\def\expandafter\da@bars\expandafter{%
      \da@bars
      \dabar@ 
    }%
  }%
  \mathrel{#3}%
  \mathrel{%
    \mathop{\da@bars}\limits
    \ifx\\#1\\%
    \else
      _{\copy0}%
    \fi
    \ifx\\#2\\%
    \else
      ^{\copy2}%
    \fi
  }%
  \mathrel{#4}%
}
\makeatother

\newcommand{\xmapsto}{\mapstochar\relbar\joinrel\xrightarrow}

\makeatletter
\newcommand*{\Relbarfill@}{\arrowfill@\Relbar\Relbar\Relbar}
\newcommand*{\xeq}[2][]{\ext@arrow 0055\Relbarfill@{#1}{#2}}
\makeatother

\makeatletter
\newcommand{\smallbullet}{} 
\DeclareRobustCommand\smallbullet{%
  \mathord{\mathpalette\smallbullet@{0.60}}
}
\newcommand{\smallbullet@}[2]{%
  \vcenter{\hbox{\scalebox{#2}{$\m@th#1\bullet$}}}%
}
\makeatother

\makeatletter
\xpatchcmd{\@tocline}
{\hfil\hbox to\@pnumwidth{\@tocpagenum{#7}}\par}
{\ifnum#1<0\hfill\else\dotfill\fi\hbox to\@pnumwidth{\@tocpagenum{#7}}\par}
{}{}
\makeatother

\makeatletter
 \def\l@subsection{\@tocline{2}{0pt}{4pc}{6pc}{}}
\def\l@subsubsection{\@tocline{3}{0pt}{8pc}{8pc}{}}
 \makeatother

\newgray{lgray}{0.9}
\newgray{vlgray}{0.95}
\newgray{midgray}{0.65}
\newgray{darkgray}{0.3}

\newcommand{\Trans}{%
exch 4 -1 roll add 3 1 roll add
}

\newcommand{\Scale}{%
dup 4 -1 roll mul 3 1 roll mul
}

\newcommand{\DrawArc}{%
\parametricplot{0}{1}{t dup mul t 1 t sub mul 2 sqrt mul t dup mul 1 t sub dup mul add 1 exch div \Scale}
}

\newcommand{\DrawPartArc}[2]{%
\parametricplot{#1}{#2}{t dup mul t 1 t sub mul 2 sqrt mul t dup mul 1 t sub dup mul add 1 exch div \Scale}
}

\newcommand{\DrawArcNeg}{%
\parametricplot{0}{1}{t dup mul t 1 t sub mul 2 sqrt mul t dup mul 1 t sub dup mul add 1 exch div \Scale -1 mul}
}

\newcommand{\DrawPartArcScale}[3]{%
\parametricplot{#1}{#2}{t dup mul t 1 t sub mul 2 sqrt mul t dup mul 1 t sub dup mul add 1 exch div \Scale 1 1 0 \SliceCoords -1 \Scale \Trans #3 \Scale 1 1 0 \SliceCoords \Trans }
}

\newcommand{\DrawArcNegScale}[1]{%
\parametricplot{0}{1}{t dup mul t 1 t sub mul 2 sqrt mul t dup mul 1 t sub dup mul add 1 exch div \Scale -1 mul 1 1 0 \SliceCoords -1 \Scale \Trans #1 \Scale 1 1 0 \SliceCoords \Trans }
}

\newcommand{\SliceCoords}{%
3 1 roll dup 3 -1 roll add 1 exch div \Scale exch 
}

\newcommand{\RSliceCoords}{%
sqrt mul \SliceCoords 
}

\newcommand{\ArcCoords}{%
2 copy add 3 1 roll dup 3 -1 roll mul 2 mul sqrt 3 -1 roll 1 exch div \Scale 
}

\newcommand{\TanToArc}[4]{%
\parametricplot{#3}{#4}{#2 2 #1 #2 mul mul sqrt #1 #2 add 1 exch div \Scale 2 copy exch -2 mul 1 add t \Scale \Trans }
}

\newcommand{\LabelTanToArc}[4]{%
\rput(!#2 2 #1 #2 mul mul sqrt #1 #2 add 1 exch div \Scale 2 copy exch -2 mul 1 add #3 \Scale \Trans){#4}
}

\newcommand{\SecToArc}[6]{%
\parametricplot{#5}{#6}{#1 #2 \ArcCoords 1 t sub \Scale #3 #4 \ArcCoords t \Scale \Trans}
}

\newcommand{\LabelSecToArc}[6]{%
\rput(!#1 #2 \ArcCoords 1 #5 sub \Scale #3 #4 \ArcCoords #5 \Scale \Trans){#6}
}

\newcommand{\SLabel}[2]{%
\psline[linecolor=gray](!#1 dup 1 add div 0.04)(!#1 dup 1 add div -0.04)
\rput(!#1 dup 1 add div -0.12){#2}
}

\newcommand{\SLabelInf}[1]{%
\psline[linecolor=gray](!0 1 0 \SliceCoords 0 0.04 \Trans)(!0 1 0 \SliceCoords 0 -0.04 \Trans)
\rput(!0 1 0 \SliceCoords 0 -0.12 \Trans){#1}
}

\newcommand{\AlphaVal}{%
4 sub dup dup mul 16 sub sqrt add 4 div
}

\newcommand{\EightNefPoint}{%
dup 1 sub dup dup mul 3 1 roll exch dup dup mul 3 1 roll mul 2 div 8 \RSliceCoords
}

\newcommand{\PutEightCircle}[1]{%
\pscircle(! #1 \EightNefPoint ){0.04}
}

\begin{document}
\title[Seshadri Constants on $\PP^1\times\PP^1$]{Seshadri constants on $\PP^1\times \PP^1$, and applications
to the symplectic packing problem}

\author{Chris Dionne}
\address{}
\email{anthonychristopherdionne@gmail.com}

\author{Mike Roth}
\address{Dept.\ of Mathematics and Statistics, Queens University, Kingston,
Ontario, Canada}
\email{mike.roth@queenu.ca}


\begin{abstract}
In this paper we compute the $r$-point Seshadri constant on $\PP^1\times\PP^1$ for those line
bundles where the answer might be expected to be governed by $(-1)$-curves.
As a consequence we obtain explicit formulas for the symplectic packing problem for $\PP^1\times \PP^1$.  
Some exact values for the Seshadri constant outside the region governed by Mori's cone theorem are also given. 
These latter results use a useful new ``reflection method''.

\smallskip
In the analysis there is a striking difference between the cases when $r$ is odd and when $r$ is even.  
When $r$ is even the problem admits an infinite order automorphism, and there are infinitely many 
$(-1)$-curves to consider.  In contrast, when $r$ is odd only a finite number (usually $4$) types of $(-1)$-curves 
are relevant to our answer.
\end{abstract}

\maketitle


\section{Introduction}
\label{sec:introduction}

\point\label{sec:r-point-def}
Let $Y$ be a smooth projective surface and $\pi\colon X\longrightarrow Y$ the blowup of $Y$ at $r$ very
general points $p_1$,\ldots, $p_r\in Y$.  We denote by $E_1$,\ldots, $E_r$ the exceptional divisors of $\pi$, with
$E_i$ lying over $p_i$, and use $E=\sum_{i=1}^{r} E_i$ for their sum.  Given an ample line bundle $L$ on $Y$,
the {\em $r$-point Seshadri constant} of $L$ is defined to be 

\vspace{-0.25cm}
\begin{equation}\label{eqn:ep-r-def-nef}
\ep_r(L):= \sup \left\{ \gamma \geq 0 \st \pi^{*}L - \gamma E\,\,\mbox{is nef}\,\rule{0cm}{0.40cm}\right\}.
\end{equation}
Equivalently

\vspace{-0.25cm}
\begin{equation}\label{eqn:ep-r-def-curve}
\ep_r(L) = \inf_C \left\{\frac{C\cdot \pi^{*}L}{C\cdot E}\right\},
\end{equation}
where the infimum is over the effective curves $C$ in $X$ which do not contain any $E_i$ as a component. 
We adopt the convention that $C$ does not have to intersect any of the $E_i$.  In such cases
$(C\cdot \pi^{*}L)/(C\cdot E)$, interpreted as $+\infty$, does not affect the infimum, and this convention allows 
us to avoid many repetitions of ``for those $C$ which also intersect at least one of the $E_i$''.

\point 
To our knowledge $r$-point Seshadri constants were first introduced by K\"{u}chle \cite{K}, for smooth
projective varieties $Y$ of arbitrary dimension.
In general very few exact values of $\ep_r(L)$ are known.  For instance, when $Y=\PP^2$, computing 
$\ep_r(\Osh_{\PP^2}(1))$ 
is equivalent to the Nagata conjecture, a problem which is open for all $r\geq 10$, $r$ not a square.

\point 
\label{sec:slope-def}
In this paper we restrict to the surface $Y=\PP^1\times\PP^1$. We use $L=\Osh_{Y}(e_1,e_2)$ 
for the line bundle on $\PP^1\times\PP^1$ of bidegree $(e_1,e_2)$.  Such a line bundle $L$ is nef if and only if
$e_1$, $e_2\geq 0$, and ample if and only if $e_1$, $e_2\geq 1$.   By the {\em slope} of $L$ we mean $e_2/e_1$, 
allowing $\infty$ if $e_1=0$ and $e_2\neq 0$. 

\point 
\label{sec:alpha-r-def}
For a positive integer $r$ we set 
$$\alpha_r:=\frac{(r-4)+\sqrt{r(r-8)}}{4}
\rule{0.25cm}{0cm}\mbox{and}\rule{0.25cm}{0cm}
\beta_r:=\frac{(r-4)-\sqrt{r(r-8)}}{4}.$$

The numbers $\alpha_r$ and $\beta_r$ are the roots of $t^2-\left(\frac{r-4}{2}\right)t+1=0$.  
When $r$ is even $\alpha_r$ and $\beta_r$ are mutually inverse units in the
ring of integers of $\QQ[\sqrt{r(r-8)}]$.  When $r\geq 10$ this ring is a real quadratic extension of $\QQ$, 
and $\alpha_r$ and $\beta_r$ are of infinite order. 

The numbers $\alpha_r$ and $\beta_r$ govern the problem of computing $\ep_r(L)$ on $\PP^1\times \PP^1$ in 
several ways.   Here is the first. 

\point
\label{sec:inner-outer-def}
We call a nef bundle $L=\Osh_{Y}(e_1,e_2)$ an {\em inner bundle} if $\frac{e_2}{e_1}\in [\beta_r,\alpha_r]$, 
and an {\em outer bundle} otherwise.
The motivation for this terminology comes from Figure \ref{fig:square-zero-cone}.
We note that whether any particular $L$ is an inner or outer bundle depends on the value of $r$. 

Let $p_1$,\ldots, $p_r\in Y$ be very general points, and 
$\pi\colon X\longrightarrow Y$ the blowup of $Y$
at $p_1$,\ldots, $p_r$ as in \S\ref{sec:r-point-def}.  For a line bundle $L$ on $\PP^1\times \PP^1$, and 
any $\gamma\geq 0$ we set 

\vspace{-0.25cm}
\begin{equation}\label{eqn:L-gamma-def}
L_{\gamma} := \pi^{*}L-\gamma\sum_{i=1}^{r} E_i = \pi^{*}L-\gamma E.
\end{equation}

By definition of the $r$-point Seshadri constant, if $\gamma>\ep_r(L)$ then the class $L_{\gamma}$ is not nef, 
and therefore there is an irreducible curve $C\subset X$ such that $L_{\gamma}\cdot C<0$.   
In \S\ref{sec:heuristic-argument} we give a heuristic argument that if $L$ is an outer bundle then one might expect
that $C$ is $K_X$-negative, and thus 
by the fact that in such a case one must also have $C^2<0$, that $C$ is a $(-1)$-curve. 

One consequence of our analysis in the paper is that this guess is correct, and we are able to explicitly 
compute $\ep_r(L)$ for all outer bundles and all $r$.  The answer appears in \S\ref{sec:intro-seshadri-outer} after
discussing another appearance of $\alpha_r$ and $\beta_r$.

\point 
\label{sec:V_r-def}
Set $F_1$ and $F_2$ to be the pullback to $X$ of the fibre classes 
$\Osh_{Y}(1,0)$ and $\Osh_{Y}(0,1)$ respectively, and let $V_r\subset H^2(X,\RR)$
be the subspace of the real N\'{e}ron-Severi group spanned by $F_1$, $F_2$, and $E$.  Thus $V_r$ is a 
three-dimensional real vector space, and for vectors 
$v=d_1F_1+d_2F_2-mE$ and $w=e_1F_1+e_2F_2-nE$ in $V_r$ the intersection pairing between $v$ and $w$ is given by 

\vspace{-0.25cm}
\begin{equation}\label{eqn:int-form-def}
v\cdot w= d_1e_2+d_2e_2-mnr.
\end{equation}

By \eqref{eqn:L-gamma-def} $L_{\gamma}$ is a class in $V_r$.    If $L_{\gamma}$ is not nef the
following argument shows that there is always
an effective curve $C$ with class in $V_r$ such that $L_{\gamma}\cdot C<0$. 

Let $C'$ be an irreducible curve such that $L_{\gamma}\cdot C'<0$, and 
let $d_1F_1+d_2F_2-\sum_{i=1}^{r} m_iE_i$ be the class of $C'$.  Since the points $p_1$,\ldots, $p_r$ are general
it follows that for any permutation $\sigma$ on $\{1,\ldots, r\}$ there is an irreducible curve $C'_{\sigma}$
with class $d_1F_2+d_2F_2-\sum_{i=1}^{r} m_{\sigma(i)}E_i$.
Moreover, $L_{\gamma}\cdot C'=L_{\gamma}\cdot C'_{\sigma}$.

Let $\sigma$ be an $r$-cycle, and set $C$ to be the sum $C := \sum_{i=1}^{r} C'_{\sigma^i}$.
Then $C$ is an effective curve, of class $rd_1F_1+rd_1F_2-\left(\sum_{i=1}^{r} m_i\right)E\in V_{r}$, and
$L_{\gamma}\cdot C = r(L_{\gamma}\cdot C')<0$.

More generally this symmetrization argument shows that the restriction of the nef and effective cones 
to $V_r$ are cones which are still dual in $V_r$.  
Thus, to understand $\ep_{r}(L)$ we may restrict our attention to $V_r$. 

\point It is easy to verify that the linear transformation $T_{r}\colon V_{r}\longrightarrow V_{r}$ given,
in the basis $F_1$, $F_2$, and $E$, by the matrix 

\vspace{-0.25cm}
\begin{equation}\label{eqn:Tr-def}
T_r:= \left[
\begin{array}{crr}
0 & 1 & 0 \\
1 & \frac{r}{2} & r \\
0 & -1 & -1 \rule{0cm}{0.45cm}\\
\end{array}
\right]
\end{equation}

preserves the intersection form. 
When $r$ is even we show in Theorem \ref{thm:Tr-transformation} that $T_r$ is an automorphism of the problem, in
the sense that if $\xi\in V_r$ is any class, then $\xi$ is nef, or effective, or represents a curve with $s$
irreducible components, if and only if $T_r(\xi)$ is respectively nef, effective, or represents a curve with
$s$ irreducible components.  These statements are not true when $r$ is odd.

The transformation $T_r$ has eigenvalues $\alpha_r$, $\beta_r$, and $1$, with respective eigenvectors (in
coordinates given by $F_1$, $F_2$, and $E$)

\begin{equation}\label{eqn:va-vb-v1-def}
\rule{0.25cm}{0cm}
v_{\alpha_r} = \left(\tfrac{1}{\alpha_r+1},\, \tfrac{\alpha_r}{\alpha_r+1},\, -\tfrac{2}{r}\right),
\rule{0.25cm}{0cm}
v_{\beta_r} = \left(\tfrac{\alpha_r}{\alpha_r+1},\, \tfrac{1}{\alpha_r+1},\, -\tfrac{2}{r}\right),
\rule{0.25cm}{0cm}\mbox{and} \rule{0.15cm}{0cm}
v_{1} = \left(-2,-2,1\right).
\end{equation}

We note that $v_1$ is the class of $K_{X}$, that $v_{\alpha_r}$ and $v_{\beta_r}$ are exchanged by 
the automorphism exchanging $F_1$ and $F_2$. Additionally, since $\alpha_r\beta_r=1$, $v_{\beta_r}$
may also be written as $v_{\beta_r}=\left(\tfrac{1}{\beta_r+1},\,\tfrac{\beta_r}{\beta_r+1},\,-\tfrac{2}{r}\right).$
Thus the Galois automorphism of $\QQ[\sqrt{r(r-8)}]$ exchanging $\alpha_r$ and $\beta_r$ exchanges
$v_{\alpha_r}$ and $v_{\beta_r}$, although we will not use this fact.

If $r\geq 10$ then $\alpha_r$ and $\beta_r$ are units of infinite order, and thus $T_r$ also has
infinite order.  Starting with a nef or effective class and iterating $T_r$ then
allows us to produce infinitely many other nef or effective classes.   This is the key to our computation 
of $\ep_{r}(L)$ for even $r$ and outer bundles $L$.

When $r\geq 10$ we have $0<\beta_r < 1 < \alpha_r$.  Thus, in forward iterations of $T_r$ vectors generally converge
(modulo scaling) to $v_{\alpha_r}$, and under backwards iterations to $v_{\beta_r}$.  
As a consequence if $r$ is even then both $v_{\alpha_r}$ and $v_{\beta_r}$ 
are limits of nef classes, and are therefore also nef.
They are also {\em square-zero classes}, $v_{\alpha_r}^2=0=v_{\beta_r}^2$, and so on the boundary of the nef cone.

\point 
\label{sec:digression}
Before proceeding to the results for $\ep_r(L)$ we make two more digressions.  For a nef line bundle 
$L=\Osh_{Y}(e_1,e_2)$ we define the {\em numerical bound}, $\eta_{r}(L)$ by 

\vspace{-0.25cm}
\begin{equation}\label{eqn:eta-r-def}
\eta_r(L) := \sqrt{\tfrac{L^2}{r}} = \sqrt{\tfrac{2e_1e_2}{r}}.
\end{equation}

The value $\eta_{r}(L)$ is precisely the value of $\gamma$ so that $L_{\gamma}^2=0$. In other words, 
$\eta_r(L)$ is the value of $\gamma$ which puts $L_{\gamma}$ on the cone of square-zero classes. The number
$\eta_r(L)$ is therefore also 
an upper bound for the Seshadri constant : $\ep_r(L)\leq \eta_{r}(L)$. 

For a vector $v\in V_r$ with $v^2=0$, and not a multiple of $v_{\alpha_r}$ or $v_{\beta_r}$, we put

\vspace{-0.25cm}
\begin{equation}\label{eqn:phi-r-def}
\phi_{r}(v): = \frac{\log\left(\frac{v\cdot v_{\beta_r}}{v\cdot v_{\alpha_r}}\right)}{2\log(\alpha_r)}.
\end{equation}

This formula is justified by the following properties 
(see \S\ref{sec:properties-of-phi-r}--\S\ref{sec:properties-of-phi-r-end}).
For such a vector $v$, $\phi_r(\lambda v)=\phi_r(v)$ for any
$\lambda\in \RR$, $\lambda\neq 0$; $\phi_r(T_r^{n}(v)) = \phi_r(v)+n$ for all $n\in \ZZ$;
and if $\tilde{v}$ is the vector obtained from $v$ by the automorphism switching $F_1$ and $F_2$, 
then $\phi_r(\tilde{v}) = - \phi_r(v)$.    Thus, $\phi_r$ is a map from the square-zero cone (up to scaling, and 
minus the lines spanned by $v_{\alpha_r}$ and $v_{\beta_r}$) to $\RR$ which takes symmetries of the 
problem to similar symmetries on $\RR$.

\bpoint{Seshadri constants for outer bundles}
\label{sec:intro-seshadri-outer}
Here we concentrate on the cases $r\geq 9$.  When $r\leq 7$ the blowup
of $\PP^1\times\PP^1$ at $r$ general points is Fano, and the answers in those cases have a different
character than the general case.  In addition, one minor aspect of our description below is not valid for $r=8$. 
These cases are discussed in \S\ref{sec:Small-r}. 

In order to describe the answers in the even case, here and in the symplectic
packing problem, it will be convenient to define several sequences $\{ s_{n,r}\}_{n\in \ZZ}$ by giving the terms
$s_{-1,r}$, $s_{0,r}$, and $s_{1,r}$, and defining all other terms by the recursion equation coming from the 
characteristic polynomial of $T_r$~:

\vspace{-0.1cm}
\begin{equation} \label{eqn:sequence-recursion}
s_{n,r} = \tfrac{r-2}{2}\left(s_{n-1,r}-s_{n-2,r}\rule{0cm}{0.40cm}\right) + s_{n-3,r}.
\end{equation}

We define the sequence $\{p_{n,r}\}_{n\in \ZZ}$ by $p_{-1,r}=0$, $p_{0,r}=0$, $p_{1,r}=r$, and determine all other
$p_{n,r}$ by the recursion \eqref{eqn:sequence-recursion}.  Similarly we define $\{m_{n,r}\}_{n\in \ZZ}$
by $m_{-1,r}=1$, $m_{0,r}=-1$, $m_{1,r}=1$, and \eqref{eqn:sequence-recursion}.   We note that 
$m_{n,r}=m_{-n,r}$ and $p_{n,r}=p_{-1-n,r}$ for all $n\in \ZZ$.

\bpoint{Theorem (Seshadri Constants for Outer Bundles)} 
\label{thm:Seshadri-for-outer-bundles} 
Suppose that $L=\Osh_{Y}(e_1,e_2)$ with $e_1,e_2\geq 1$, and
that $L$ is an outer bundle, i.e., that $\frac{e_2}{e_1}\not\in [\beta_r,\alpha_r]$.  

{\bf If $r$ is odd, $r\geq 9$}.  Then

\vspace{-0.50cm}
\begin{equation}\label{eqn:Seshadri-odd-n}
\ep_r(L) = \left\{ 
\begin{array}{cl}
e_2 & \mbox{if $\frac{e_2}{e_1}\leq \frac{2}{r+1}$}, \\
\frac{2e_1+(r-1)e_2}{2r} & \mbox{if $\frac{e_2}{e_1} \in [\frac{2}{r+1},\beta_r]$}, \rule{0cm}{0.6cm}\\
\frac{(r-1)e_1+2e_2}{2r} & \mbox{if $\frac{e_2}{e_1} \in [\alpha_r,\frac{r+1}{2}]$}, \rule{0cm}{0.6cm}\\
e_1 & \mbox{if $\frac{r+1}{2}\leq \frac{e_2}{e_1}$}. \rule{0cm}{0.6cm} \\
\end{array}
\right.
\end{equation}

{\bf If $r$ is even, $r\geq 10$}. 
Set $v_{L}=(e_1,e_2,-\eta_{r}(L))\in V_r$, i.e., set $v_{L}$ to be the class of
$L_{\gamma}$ with $\gamma=\eta_r(L)$, and put $n=\lfloor \phi_r(v_{L})+\frac{1}{2}\rfloor$, where
$\lfloor x \rfloor$ denotes the largest integer $\leq x$.  Then

\begin{equation}\label{eqn:Seshadri-even-n}
\ep_r(L) = \frac{e_1p_{n,r}+e_2p_{n-1,r}}{rm_{n,r}}.
\end{equation}

\vspace{-0.15cm}

\point
\label{sec:curve-explanation}
The explanation for \eqref{eqn:Seshadri-even-n} is as follows.  
Consider the classes $C_{n,r}=(p_{n-1,r},p_{n,r},-m_{n,r})$ defined
by the sequences $\{p_{n,r}\}_{n\in \ZZ}$ and $\{m_{n,r}\}_{n\in \ZZ}$.  Then $C_{0,r}=(0,0,1)$ is the class of $E$,
and for all $n\in \ZZ$, $C_{n,r}=T_{r}^{n}(C_{0,r})$.   In Theorem \ref{thm:Outer-nef} we show that when $n$ 
is even the $C_{n,r}$ generate the effective cone of curves in $V_{r}$ whose slopes lie outside of 
$[\beta_r,\alpha_r]$.   
It follows that for an outer bundle $L$ one of these curves determines $\ep_r(L)$. 
The formula with $\phi_{r}$ above is one possible method of locating the correct $n$, and intersecting 
with $C_{n,r}$ then gives \eqref{eqn:Seshadri-even-n}.

For the proofs of the results in the even case see \S\ref{sec:Outer-Even},
and for the odd case, \S\ref{sec:Outer-Odd}. 

\bpoint{Applications to the symplectic packing problem} 
\label{sec:Symplectic-Packing-applications}
We recommend \cite[\S1]{MP} and \cite{B2} for a discussion of the history of this problem and the reasons for its 
interest.  Here we give a brief outline oriented towards our application of the previous results to the 
symplectic packing problem for $\PP^1\times\PP^1$. 

Let $(M,\omega_{M})$ be a closed symplectic manifold of real dimension $4$,  and let $B_{\lambda}\subset \RR^4$
denote the ball of radius $\lambda$ centred at $0$, equipped with the restriction of the standard symplectic form
on $\RR^4$~: $\omega_{\RR^4} = dx_1\wedge dy_1 + dx_2\wedge dy_2$. 

For a given $r$, consider all possible symplectic embeddings of $r$ disjoint copies of $B_{\lambda}$ into $M$,
and denote by $\nuhat_r(M)$ the supremum of the volumes which can be filled by such embeddings (i.e., the supremum
of $r\pi^2\lambda^4$, over those $\lambda$ for which there exists such a symplectic embedding of $r$
disjoint $B_{\lambda}$ into $M$).  
Finally, set $\nu_{r}(M) = \nuhat_{r}(M)/\Vol(M)$, where the volume of $M$ (like the volume of $B_{\lambda}$)
is computed using the volume form $\omega_{M}\wedge\omega_{M}$.

Two basic questions are : (\Euler{1}) What is $\nu_r(M)$ for different values of $r$?; and 
(\Euler{2}) for which $r$ does $\nu_r(M)=1$?  
If $\nu_r(M)<1$ one says that there is a {\em packing obstruction}, while if 
$\nu_{r}(M)=1$ one says that there is a {\em full packing}.

Let $Y$ be a smooth projective surface, and $L$ a real ample class on $Y$. The first Chern class $c_1(L)$ can be
represented by a K\"{a}hler form $\omega_{L}$, which, when written out in terms of the underlying real coordinates,
is a real symplectic form.
We consider the packing problem for the real manifold $M$ underlying $Y$, with symplectic form $\omega_{L}$.
To align our notation with the notation in the rest of the paper, we will use $\nu_{r}(L)$ for the value of
$\nu_{r}(M)$ in this situation.

\point
A remarkable discovery of \cite{MP} is that $(-1)$-curves in $X$ provide obstructions to full packings.  
Even more striking is that this obstruction looks much like the Seshadri constant, except with the test curves $C$ 
in \eqref{eqn:ep-r-def-curve} limited to $(-1)$-curves.
To set this up we first extend definition \eqref{eqn:eta-r-def} to any such pair $(Y,L)$ by setting 
$\eta_{r}(L) = \sqrt{L^2/r} = \sqrt{\Vol(M)/r}$. 
Then we set

\vspace{-0.2cm}
\begin{equation}\label{eqn-ep-tilde-r-def} 
\eptilde_r(L) = \min\left(\inf_C \left\{\frac{C\cdot \pi^{*}L}{C\cdot E}\right\}, \eta_{r}(L)\right)
\end{equation}

\vspace{-0.2cm}
where this time the infimum is over irreducible $(-1)$-curves $C$ in $X$ distinct from the exceptional divisors $E_i$.
Following our convention in \S\ref{sec:r-point-def}, if there is no such $(-1)$-curve $C$ which intersects any of the
$E_i$ the infimum is interpreted as $\infty$, and then $\eptilde_{r}(L) = \eta_{r}(L)$. 

The obstruction result of \cite{MP} is that one always has 
$\nu_{r}(L) \leq \left(\frac{\eptilde_{r}(L)}{\eta_{r}(L)}\right)^2$.
The paper is concerned with the case $Y=\PP^2$, but the obstruction argument does not depend on this. 
Even more remarkably, a result of Biran, \cite[Theorem 6.A]{B1}, asserts that there is a class of surfaces, which 
includes $\PP^2$ and ruled surfaces, where one has $\nu_{r}(L) = \left(\frac{\eptilde_{r}(L)}{\eta_{r}(L)}\right)^2$
for all $L$.

The theorem of Biran completely determines the packing constant $\nu_{r}(L)$ for real ample line bundles $L$. But,
for a given $L$ it is not immediate what the value of $\nu_{r}(L)$ is, nor the values of $r$ such that $\nu_{r}(L)=1$.  
Using Theorem \ref{thm:Seshadri-for-outer-bundles} we give formulas for $\nu_r(L)$ for all real ample 
line bundles on $Y=\PP^1\times\PP^1$;  equivalently, by \cite[Theorem 1.1]{LM}, for all symplectic forms on the 
underlying real manifold.
As before we list the results for $r\geq 9$; the results for $r\leq 8$ appear in \S\ref{sec:packing-constants-small-r}.

\bpoint{Theorem (Formulas for the symplectic packing constant)}
\label{thm:symplectic-packing-results}
Let $L$ be a real ample line bundle of type $(e_1,e_2)$ (i.e., $e_1$ and $e_2$ are positive real numbers). 

{\bf If $r$ is odd, $r\geq 9$}.  Then 

\begin{equation}\label{eqn:Seshadri-tilde-odd-n}
\nu_r(L) = \left\{ 
\begin{array}{cl}
\frac{re_2}{2e_1} & \mbox{if $\frac{e_2}{e_1}\leq \frac{2}{r+1}$}, \\
\frac{(2e_1+(r-1)e_2)^2}{8re_1e_2} & 
\mbox{if $\frac{e_2}{e_1} \in [\frac{2}{r+1},\frac{2}{(\sqrt{r}-1)^2}]$}, \rule{0cm}{0.6cm}\\
1 & \mbox{if $\frac{e_2}{e_1} \in [\frac{2}{(\sqrt{r}-1)^2},\frac{(\sqrt{r}-1)^2}{2}]$}, \rule{0cm}{0.6cm} \\
\frac{((r-1)e_1+2e_2)^2}{8e_1e_2} & 
\mbox{if $\frac{e_2}{e_1} \in [\frac{(\sqrt{r}-1)^2}{2},\frac{r+1}{2}]$}, \rule{0cm}{0.6cm}\\
\frac{re_1}{2e_2} & \mbox{if $\frac{r+1}{2}\leq\frac{e_2}{e_1}$}. \rule{0cm}{0.6cm} \\
\end{array}
\right.
\end{equation}

{\bf If $r$ is even, $r\geq 10$}.  Then 

\vspace{-0.25cm}
\begin{equation}\label{eqn:Seshadri-tilde-even-n}
\rule{2.0cm}{0cm}
\nu_r(L) = \left\{ 
\begin{array}{cl}
1 & \mbox{if $\frac{e_2}{e_1} \in [\beta_r,\alpha_r]$} \\
\frac{r\left(\ep_{r}(L)\right)^2}{2e_1e_2}  & \mbox{if $\frac{e_2}{e_1}\not\in [\beta_r,\alpha_r]$
(with $\ep_r(L)$ computed by the rule in \eqref{eqn:Seshadri-even-n})}. \rule{0cm}{0.6cm} \\
\end{array}
\right.
\end{equation}

\bpoint{Conditions for full packings} 
\label{sec:full-packings}
Define sequences $\{q_{n,r}\}_{n\in \ZZ}$ by 
$q_{-1,r}=1$, $q_{0,r}=0$, $q_{1,r}=1$, and the recursion \eqref{eqn:sequence-recursion}.
For this sequence one has $q_{n,r}=q_{-n,r}$ for all $n$.  
Taking into account the cases $r\leq 8$ (see \S\ref{sec:Small-r}), and reversing the formulae 
in Theorem \ref{thm:symplectic-packing-results}, we get the following answer to question (\Euler{2}).  

\bpoint{Theorem (Conditions for full packings)}
\label{thm:full-packings}

{\bf If $r$ is odd}.  Then $\nu_{r}(L)=1$ if and only if
$r\geq \max\left(\left(\sqrt{\frac{2e_2}{e_1}}+1\right)^2, 
\left(\sqrt{\frac{2e_1}{e_2}}+1\right)^2,9\right)$.

{\bf If $r$ is even}.  Then $\nu_{r}(L)=1$ if and only if 

\RomanList
\begin{enumerate}
\item $r\geq \dfrac{2(e_1+e_2)^2}{e_1e_2}$, or 
\item $r$ is a value for which $\frac{e_2}{e_1}$ is equal to $\frac{q_{n+1,r}}{q_{n,r}}$ for some $n$.
\rule{0cm}{0.5cm}
\end{enumerate}

For a given $(e_1,e_2)$, there is at most one value of $r$ for which ({\em ii}) occurs, see
Theorem \ref{thm:uniqueness-of-r}.

\vspace{-0.2cm}
\bpoint{Examples}
\label{sec:packing-examples}
The lower bounds on $r$ differ in the even and odd cases.
Consequently by picking a line bundle with an extreme slope, we can find examples of line bundles $L$ 
with ranges where full packings exist only for even $r$. 
Here are two examples with similar slopes.

For the first, we give an example where case ({\em ii}) above does not occur.
If $L=\Osh_{Y}(2,401)$ then by the formulas above there is 

\begin{itemize}
\item no full packing for any $r\leq 405$;
\item a full packing for every $r\geq 443$; and
\item for $r\in [406,442]$ a full packing only for even $r$.
\end{itemize}

Similarly, If $L=\Osh_{Y}(1,200)$ then there is

\begin{itemize}
\item no full packing for any $r\leq 399$;
\item a full packing for every $r\geq 441$; and
\item for $r\in [400,440]$ a full packing only for even $r\geq 406$ and $r=400$.
\end{itemize}

In this second example $r=400$ is an ``unusual'' $r$, i.e., appears because of case ({\em ii}). 

This phenomena seems very surprising to the authors, even knowing the proofs of the formulas. For instance, 
returning to the first example, there is a full packing when $r=410$.  If we look for a packing with $r=409$, then
we could start with a packing for $r=410$ and use $409$ of the balls.  Admittedly, that is not yet a full packing,
but surely it would be possible to increase the radius and move the centres just a little bit to make up for it, 
and not have the balls intersect \ldots ?  Of course, the results above say that it is not possible.     
As a consistency check, Theorem \ref{thm:Seshadri-for-outer-bundles} gives  
$\nu_{409}(L) = \frac{654481}{656036}\approx 0.9976297\ldots$, larger than 
the ratio $\frac{409}{410} \approx 0.99756097\ldots$ achievable using $409$ out of the $410$ balls but without 
increasing the radius. 

The proofs for the above results on symplectic packing, using the previous results about Seshadri constants, appear 
in \S\ref{sec:Symplectic-Packing}.

\bpoint{Results for inner bundles}
\label{sec:inner-bundle-summary}
One implication of the SGSH conjecture (see for example \cite[\S1.4]{CHMR}) and our analysis of which $(-1)$-curves 
affect Seshadri constants, is that there should be a portion of the nef cone which is round.
Specifically, for $r\geq 9$, we should have $\ep_r(L)=\nu_{r}(L)$ for all $L$ whose slopes are in, respectively,
$[\beta_r,\alpha_r]$ if $r$ is even, and $[\frac{2}{(\sqrt{r}-1)^2},\frac{(\sqrt{r}-1)^2}{2}]$ if $r$ is odd.

If this description of the nef cone is correct, then the boundary of the nef cone, for slopes in the 
ranges indicated above, consists of classes $\xi$ which are nef, square-zero (i.e., $\xi^2=0$) and 
$K_X$-positive : $K_X\cdot \xi>0$.    In this paper we call such classes inner square-zero nef classes.

Finding such classes is quite useful.  By definition, if $\xi$ is nef, then there are no effective classes
on the half plane $\xi^{<0}$. If $\xi$ is $K_X$-positive, this half plane will contain a large proportion of
$K_X$-positive classes, and it is exactly these classes whose existence we usually have the greatest difficulty
in ruling out.  In addition, if $\xi^2=0$ it means that $\xi$ is on the boundary of the nef cone, and so provides
the strongest condition on restricting effective classes.

One of the contributions of this paper is to construct such inner square-zero nef classes for all $r\geq 9$.  
To our knowledge, this is the first construction of such classes on the blowup of $\PP^1\times\PP^1$ at $r$
general points.   These classes are constructed in \S\ref{sec:Inner-square-zero-classes-reflections}, using a 
new ``reflection method''. 
If $r$ is even we obtain, using $T_r$, infinitely many such classes, but if $r$ is odd we only construct finitely many.
In \S\ref{sec:Inner-square-zero-classes-pullbacks} we use pullback maps to construct other infinite families of
such classes when $r$ is even.

Given such a class $\xi$, say of the form $\xi=\left(e_1,e_2,-\sqrt{\frac{2e_1e_2}{r}}\right)$ (the last coordinate
is determined by the condition that $\xi^2=0$), for the bundle $L=\Osh_{Y}(e_1,e_2)$ we then have
$\ep_r(L)=\eta_r(L)$.    We are thus able to exhibit classes which achieve the predicted value of the Seshadri 
constant.

\bpoint{Relation with other work, I}
Seshadri constants on $\PP^1\times\PP^1$ and the related symplectic packing problem
were studied in the 2005 Ph.D.\ thesis of W.\ Syzdek, the published version of which appears as \cite{Sy}. 

In \cite{Sy} Syzdek finds the same curves $C_{n,r}$ (from \S\ref{sec:curve-explanation}) which we use to compute 
the Seshadri constant in the even case.  More precisely, our $C_{n,r}$ are each the disjoint union of $r$ $(-1)$-curves 
(for instance $C_{n,r}=T^{n}_{r}(E)$, and $E$ is the disjoint union of the $r$ exceptional 
divisors), and Syzdek finds instead the classes of these $(-1)$-curves.  Specifically, our curve $C_{n,r}$ 
with $n\geq 1$ is the symmetric orbit of the curve called $M_{D_{n+2}}$ in \cite[Proposition 3.9]{Sy}, 
with $l=\frac{r-6}{2}$ and with the bidegrees switched.   These curves ($C_{n,r}$, or its components $D_{n+2}$) 
impose the same conditions on the $r$-point Seshadri constants.

The approach in \cite{Sy} cannot rule out that there may be other curves which affect the Seshadri constant
when $\frac{e_2}{e_1}\not\in [\beta_r,\alpha_r]$, and so for the corresponding line bundles can only give
an upper bound on $\ep_{r}(L)$, and an upper bound on the symplectic packing number.  
In our argument we can conclude that these upper bounds are the actual values.
The extra piece of information in our method is that we know that the iterates $T^{n}_r(F_2)$ are nef, 
and duality of the nef and effective cones then eliminates the possibility of other such curves. 

In the case that $r$ is odd, the $(-1)$-curves we find also already appear in \cite[Table 4]{Sy}.   For instance,
our curve of bidegree $(\frac{r-1}{2},1)$ is the curve of bidegree $(k+n+3,1)$ in table 4 when $r=2k+2n+7$.
Theorem 3.17 of \cite{Sy} seems to claim, in the case $r$ is odd, that the curves in \cite[Table 4]{Sy} compute the
Seshadri constant for all $(e_1,e_2)$ with $\frac{e_2}{e_1}\not\in [\frac{2}{(\sqrt{r}-1)^2},\frac{(\sqrt{r}-1)^2}{2}]$.
We do not know how to justify this claim since we cannot rule out the possibility that there may be curves $C$
with $C^2<0$ and $C\cdot K_{X}>0$ which could impose a stronger condition on the Seshadri constant.  
In fact, we have had to do some work in \S\ref{sec:Outline-of-outer-odd}--\S\ref{sec:proof-of-outer-odd-theorem}
to show that if such a curve exists, it at least could not affect 
line bundles with $\frac{e_2}{e_1}\not\in [\beta_r,\alpha_r]$.  

As part of question (\Euler{2}) in \S\ref{sec:Symplectic-Packing-applications},
one may ask for an $r_0$ so that for $r\geq r_0$ one has $\nu_r(L)=1$.  In general, as the examples
in \S\ref{sec:packing-examples} suggest, it is the odd $r$ which determine $r_0$.  Our bound in 
Theorem \ref{thm:full-packings} is sharp.
Setting $s=\frac{e_2}{e_1}$, and assuming $s\geq 2$ to simplify the discussion,
we get that there is a full packing for all $r$ greater than $(\sqrt{2s}+1)^2 = 
2s+2\sqrt{2s}+1$.  In contrast, \cite{Sy} (Definition 3.2, formula for $R_0$ with $a=1$, $b=s$)
gives a slightly worse estimate of 
$$\frac{3+2s+3s^2}{2s} + \frac{(1+s)\sqrt{2(1+s^2)}}{s} = \left(\tfrac{3}{2}+\sqrt{2}\right)s+ 
(1+\tfrac{3}{\sqrt{2}})+O\left(\tfrac{1}{s}\right)\rule{0.25cm}{0cm}\mbox{as $s\to\infty$}.$$
Finally, we should note that the dichotomy of behaviour between even and odd $r$, one aspect of the problem which
we find surprising, already appears in \cite{Sy}.
For instance, in the estimates $R_0$ and $r_0$ in \cite[Definition 3.2]{Sy}.

In summary, the improvements in this paper over the results of \cite{Sy} are~: 
(\Euler{1}) When $r$ is even to give exact values of $\ep_r(L)$ for outer bundles $L$ and exact values of
$\nu_r(L)$ for all ample $L$; 
(\Euler{2}) When $r$ is odd to justify the calculation of $\ep_r(L)$ for outer bundles $L$ and to give the 
exact region where $\nu_r(L)=1$.
We are thus able to calculate explicit answers to the symplectic packing problem on $\PP^1\times\PP^1$.
(\Euler{3}) To produce inner square-zero nef classes for all $r\geq 9$. Thus, to produce
inner bundles where the Seshadri constant can be computed exactly.

The authors also think that the introduction of the numbers $\alpha_r$ and $\beta_r$, the realization that the problem 
has an infinite order automorphism when $r$ is even, $r\geq 8$, and and the graphical reasoning from 
\S\ref{sec:airplane-hanger}--\S\ref{sec:graphical-arguments} greatly simplify the analysis of the problem.

\bpoint{Relation with other work, II}
\label{sec:other-work-II}
The paper \cite{DTG} gives lower bounds on $\ep_r(L)$ for those $L$ whose Seshadri constant is not affected
by $(-1)$-curves. 

When $r$ is odd, this means line bundles $L$ with 
$\frac{e_2}{e_1}\in [\frac{2}{(\sqrt{r}-1)^2},\frac{(\sqrt{r}-1)^2}{2}]$. Then \cite[Theorem 5]{DTG} gives
the lower bound $\ep_r(L)\geq \eta_{r}(L)\cdot (1-\frac{1}{5r})^{\frac{1}{2}}$.    

When $r$ is even, the results of \cite{DTG} apply to all inner bundles, and \cite[Theorem 4]{DTG} gives the lower bound
$\ep_r(L) \geq \eta_{r}(L)\cdot(1-\frac{2}{9r})^{\frac{1}{2}}$. 

As discussed in \S\ref{sec:inner-bundle-summary} in this case we are able to find inner bundles $L$
where $\ep_r(L)=\eta_r(L)$.  Using these bundles and convexity of the Seshadri constant then gives lower bounds
on $\ep_r$ which are better than the lower bound above on various regions of 
the intervals above.
See the discussion in \S\ref{sec:comparison-of-lower-bounds}.

\bpoint{Limitations of this paper} In the study of Seshadri problems on blowups of rational surfaces, and in particular the Nagata conjecture, the sticking point is our inability to either rule out all
$K_{X}$-positive curves $C$ with $C^2<0$, or exhibit one which exists.  Unfortunately this paper is no exception.  

However, the construction of the inner square-zero nef classes in 
\S\ref{sec:Inner-square-zero-classes-reflections}--\S\ref{sec:Inner-square-zero-classes-pullbacks} does eliminate
a large range of such classes, and seems to the authors to be a useful step forward.

Second, the most precise results about Seshadri constants in this paper are for outer bundles, those
bundles $L=\Osh_{Y}(e_1,e_2)$ with $\frac{e_2}{e_1}\not\in [\beta_r,\alpha_r]$.  As $r\to\infty$ we have
$\beta_r\to 0$ and $\alpha_r\to\infty$.  Thus, as $r$ increases, the region of our ignorance also increases, and
the region of complete understanding shrinks to zero.

\bpoint{Organization of the paper}
In \S\ref{sec:setup} we describe a graphical way of representing and arguing about the problem 
and give a heuristic argument that irreducible curves affecting the Seshadri constant of outer bundles should 
be $(-1)$-curves.  This picture also explains the appearance of $\alpha_r$ and $\beta_r$ in the problem.

In \S\ref{sec:Outer-Even} we show that $T_r$ is an automorphism of the problem when $r$ is even, 
calculate the nef and effective cones for classes whose slope is outside of $[\beta_r,\alpha_r]$, and compute
the Seshadri constants for outer bundles for even $r\geq 10$.  
In \S\ref{sec:Outer-Odd} we compute the Seshadri constants for outer bundles for odd $r\geq 9$.
In \S\ref{sec:Small-r} we give the results for all $r$, even and odd, with $r\leq 8$.

In \S\ref{sec:brief-study-of-slopes} we study the slopes $\frac{q_{n+1,r}}{q_{n,r}}$ which show up in the
exceptional case ({\em ii}) in the symplectic packing problem (\S\ref{thm:full-packings}), as well
as establish the properties of the map $\phi_r$ defined in \eqref{eqn:phi-r-def}.

In \S\ref{sec:Symplectic-Packing} we use the results of the previous sections to establish the results
on symplectic packings, Theorems \ref{thm:symplectic-packing-results} and \ref{thm:full-packings}.
In \S\ref{sec:Reflections} we give the reflection theorem, a method of producing nef classes using certain
types of specializations.
In \S\ref{sec:Inner-square-zero-classes-reflections} we use the reflection theorem to construct inner square-zero
nef classes in both the odd and even cases. 
Finally in \S\ref{sec:Inner-square-zero-classes-pullbacks} we use pullback maps to produce other families
of inner square-zero nef classes when $r$ is even, and an interesting family of bounds when $8\mid r$. 

\bpoint{Acknowledgements}
Most of the results of this paper form part of the Ph.D.\ thesis of the first author, supervised by the second author. 
We thank Piotr Pokora for helpful comments on previous versions of these results. 
The second author also thanks the Mathematics Department of Humboldt-Universit\"{a}t zu Berlin, 
where some of this work was carried out, for their excellent working conditions.

\tableofcontents

\newpage
\section{The square-zero cone and graphical arguments}
\label{sec:setup}

\bpoint{The square-zero cone} 
\label{sec:airplane-hanger} 
Let $X$ be the blowup of $\PP^{1}\times\PP^{1}$ at $r$ general points. 
As in \S\ref{sec:V_r-def}, let $V_r \subset H^2(X,\RR)$ be the real subspace spanned by the pullbacks $F_1$, $F_2$,
of the fibre classes from $\PP^{1}\times\PP^{1}$, and the sum $E$ of the exceptional divisors. 
We are interested in studying the restriction of the nef and effective cones to $V_r$. 

Let $C= d_1F_1+d_2F_2-mE$ be a class in $V_{r}$.  If $d_1<0$ or $d_2<0$ then $C$ is neither effective nor nef.
If $d_1$, $d_2\geq 0$, but $m<0$, then $C$ is effective but not nef.   The real interest is therefore when 
$d_1$, $d_2$, and $m\geq 0$, and we restrict to that octant from now on.

\hfill
\begin{tabular}{c}
\includegraphics[scale=0.1]{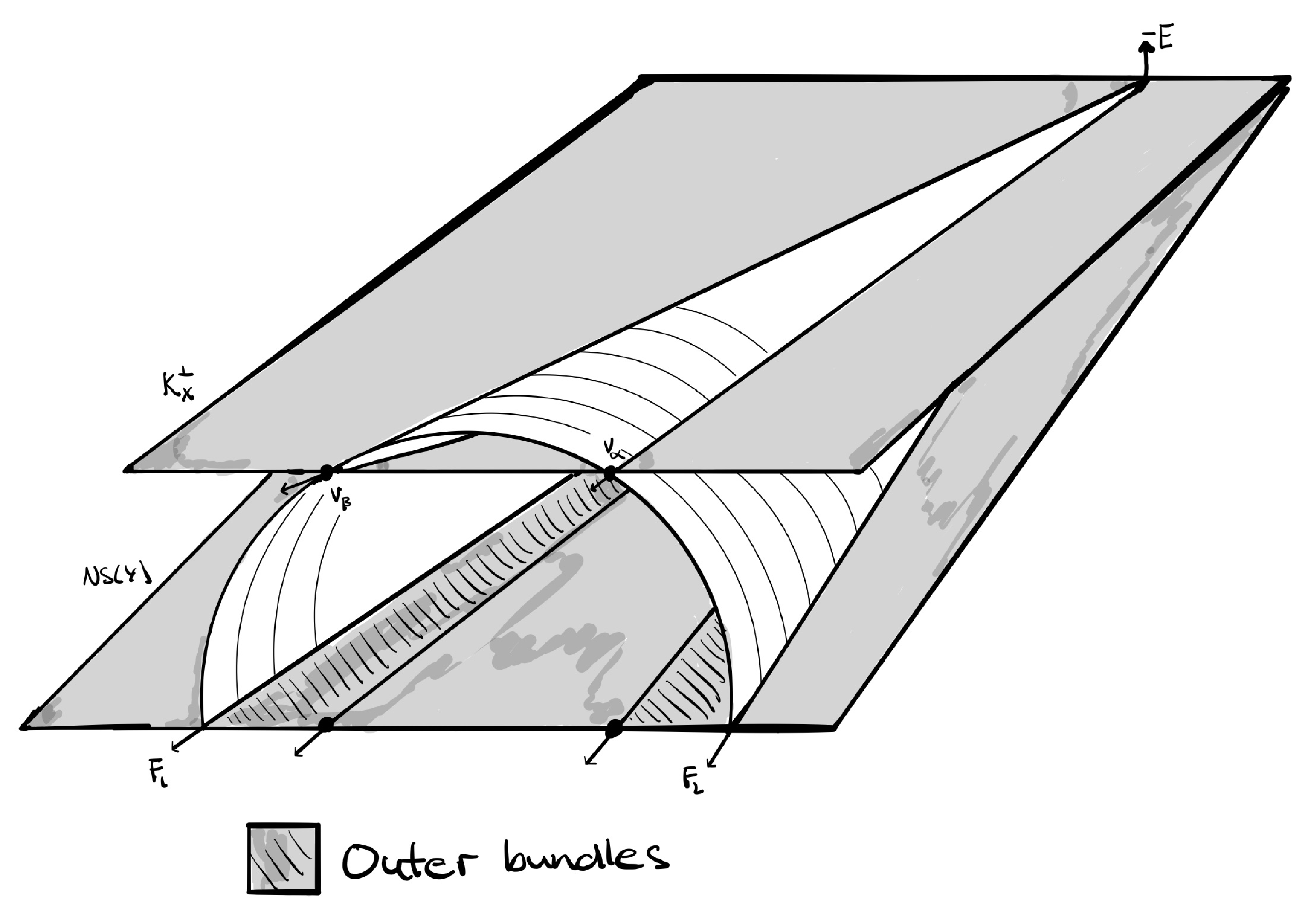}  \\
\small \Fig\label{fig:square-zero-cone} \\
{\small The square-zero cone, $r\geq 9$} \\
\end{tabular}

\vspace{-6.7cm}
\parshape 2 0cm 11.5cm 0cm 11.0cm
In that octant a key object of interest for us is the {\em square-zero cone}, those classes $\xi\in V_r$ such
that $\xi^2=0$. A picture of this cone in the octant where $d_1$, $d_2$, $m\geq 0$ is shown at right.

\parshape 1 0cm 8.5cm
In the picture, the plane on the base is the subspace spanned by $F_1$ and $F_2$, i.e., the image of the real 
N\'{e}ron-Severi group of $\PP^1\times \PP^1$ under the pullback map to $X$.   The curved shape is the square-zero
cone, and it meets the base plane in the rays spanned by $F_1$ and $F_2$.

\parshape 5 0cm 9.0cm 0cm 9.0cm 0cm 9.0cm 0cm 9.0cm 0cm \textwidth
The plane at the top of the picture is the subspace of classes orthogonal to $K_{X}$, i.e., those classes $\xi$
so that $\xi\cdot K_X=0$.  
The $K_{X}$-negative classes lie below the plane, and the $K_{X}$-positive classes lie above. 

When $r\leq 7$ this plane lies strictly above the square-zero cone,  when $r=8$ this plane
is tangent to the cone, and when $r\geq 9$ this plane intersects the cone in two rays.   These rays are
the rays spanned by the vectors $v_{\alpha_r}$ and $v_{\beta_r}$ defined in \eqref{eqn:va-vb-v1-def}. 

\vspace{-0.3cm}
\hfill
\begin{tabular}{c}
$
\begin{array}{|c|ccc|}
\hline\hline
\cdot & v_{\alpha} & v_{\beta} & K_{X} \\
\hline
v_{\alpha} & 0 & \frac{r-8}{r} & 0 \rule{0cm}{0.4cm}\\
v_{\beta}  & \frac{r-8}{r} & 0 & 0 \rule{0cm}{0.4cm}\\
K_{X} & 0 & 0 & 8-r \\
\hline\hline
\end{array}
$
\\
\small \Fig \label{fig:intersection-table} \rule{0cm}{0.4cm}\\
\small Intersection table \\
\end{tabular}\hspace{-0.3cm}

\vspace{-3.1cm}

\parshape 1 0cm 12cm
The intersection matrix for $v_{\alpha_r}$, $v_{\beta_r}$ and $K_X$ (= $v_1$ in the notation in 
\eqref{eqn:va-vb-v1-def}) is shown at right, and is easily verified from the formulas for those classes,
and the formula \eqref{eqn:int-form-def} for the intersection form. 

\parshape 4 0cm 12cm 0cm 12cm 0cm 12cm 0cm \textwidth 
This is perhaps the quickest way to check that $v_{\alpha_r}$ and $v_{\beta_r}$ span the rays above.  The
table shows that they are both square-zero classes, and orthogonal to $K_X$.  Note that when $r=8$ we have
$v_{\alpha_8}=v_{\beta_8}=-\frac{1}{4} K_{X}$; this is the case where the plane $K_{X}^{\perp}$ is tangent to the 
square-zero cone.

The projection of the rays spanned by $v_{\alpha}$ and $v_{\beta}$ onto the base plane are rays of slopes $\alpha_r$
and $\beta_r$ respectively.   Those rays in the base plane whose slopes are outside of $[\beta_r,\alpha_r]$ are the 
outer bundles, and those with slopes in the interval $[\beta_r,\alpha_r]$ are the inner bundles 
(\S\ref{sec:inner-outer-def}). 

\bpoint{Graphical arguments}
\label{sec:graphical-arguments}
There are several places in the paper where an argument can be simply expressed by a picture which would otherwise
require a chain of uninformative inequalities.  This graphical way of thinking has also guided our approach to the 
problem.  In this subsection we explain our graphical notation, and several elementary facts which can be seen
from this point of view. 

We restrict ourselves to the octant $d_1$, $d_2$, $m\geq 0$ of \S\ref{sec:airplane-hanger}.  The nef and effective
cones are stable under scaling by positive real numbers, and so it is sufficient to consider 
Figure \ref{fig:square-zero-cone} up to scaling, which we represent as a diagram of the type 
in Figure \ref{fig:basic-diagram}. 

\hfill
\begin{tabular}{c}
\begin{tabular}{c}
\psset{yunit=3cm, xunit=5cm}
\begin{pspicture}(0,-0.2)(1,0.8)
\psset{linecolor=gray}
\psset{linecolor=black}
\psset{linecolor=gray}
\psline(!1 0 0 \SliceCoords)(!0 1 0 \SliceCoords) 
\psline(! 12 \AlphaVal 1 \ArcCoords exch pop 0 exch)(! 12 \AlphaVal 1 \ArcCoords exch pop 1 exch)
\SLabel{0}{\tiny $0$}
\SLabel{1 12 \AlphaVal div}{\tiny $\beta_r$} 
\SLabel{12 \AlphaVal }{\tiny $\alpha_r$} 
\SLabelInf{\tiny $\infty$} 
\psset{linecolor=black,fillstyle=none}
\DrawArc
\psset{fillstyle=solid,fillcolor=white,linecolor=black,linestyle=solid,linecolor=black}
\pscircle(!12 \AlphaVal 1 \ArcCoords){0.04}
\rput(!12 \AlphaVal 1 \ArcCoords 0.1 add){$v_{\beta}$}
\pscircle(!1 12 \AlphaVal \ArcCoords){0.04}
\rput(!1 12 \AlphaVal \ArcCoords 0.1 add){$v_{\alpha}$}
\rput(!1 12 \AlphaVal \ArcCoords 0.3 0 \Trans){\tiny $K_X^{\perp}$}
\rput(!1 1 \ArcCoords 3 div){\color{gray}\small $+$}  
\rput(!1 1 \ArcCoords 0.25 add){\color{gray}\small $-$}
\rput(-0.05,0.3){\color{gray}\small $-$}
\rput(1.05,0.3){\color{gray}\small $-$}
\end{pspicture} \\
\end{tabular} \\
\small \Fig \label{fig:basic-diagram} \\
\small (Figure \ref{fig:square-zero-cone} up to scaling) \\
\end{tabular}

\vspace{-3.7cm}
\parshape 1 0cm 10.5cm
In this picture the curve represents the square-zero cone, the line on the bottom the portion of the nef cone
spanned by $F_1$ and $F_2$, and the upper line the plane $K_{X}^{\perp}$.  We label a class $(e_1,e_2,0)$ along the 
bottom by its slope $\frac{e_2}{e_1}$, 
so that $F_1$ corresponds to slope $0$ and $F_2$ to slope $\infty$. 

\parshape 4 0cm 10.5cm 0cm 10.5cm 0cm 10.5cm 0cm \textwidth
The signs $+$ and $-$ in this diagram are a reminder that classes inside the square-zero cone have positive 
self-inter\-section, and classes outside have negative self-intersection, and we will omit them from further
diagrams.
We will also sometimes omit the line for $K_{X}^{\perp}$.

Here are three arguments we will use frequently.   We first give the associated pictures, and then
explain what the statements are. 

\vspace{0.2cm}
\begin{centering}
\begin{tabular}{|c|c|c|}
\hline
\begin{tabular}{c}
\psset{yunit=3cm, xunit=5cm}
\begin{pspicture}(-0.2,0)(0.7,1.0)
\psset{linestyle=dashed,dash = 3pt 2pt, linecolor=gray}
\TanToArc{1}{0.175}{-0.55}{0.55}
\LabelTanToArc{1}{0.175}{-0.65}{\small $\xi^{\perp}$}
\rput(0.4,0.77){\tiny\color{gray}$+$}
\rput(0.35,0.87){\tiny\color{gray}$-$}
\psset{linecolor=black,linestyle=solid}
\DrawPartArc{0.05}{0.47}
\psset{fillstyle=solid,fillcolor=white}
\pscircle(!1 0.175 \ArcCoords){0.04}
\rput(!1 0.175 \ArcCoords 0.12 add){\small $\xi$}
\end{pspicture}
\end{tabular}
&
\begin{tabular}{c}
\psset{yunit=3cm, xunit=5cm}
\begin{pspicture}(0,0)(1,1)
\psset{linecolor=gray}
\SecToArc{9}{2}{1}{2}{-0.4}{1.4}
\LabelSecToArc{9}{2}{1}{2}{1.6}{\small $C^{\perp}$}
\rput(0.03,0.45){\color{gray}\tiny$+$}
\rput(0.02,0.57){\color{gray}\tiny$-$}
\psset{linestyle=dashed,dash = 3pt 2pt, linecolor=gray}
\TanToArc{9}{2}{-0.25}{0.55}
\TanToArc{1}{2}{-0.55}{0.25}
\psset{linecolor=black,linestyle=solid}
\DrawPartArc{0.00}{0.85}
\psset{fillstyle=solid,fillcolor=white}
\pscircle(!9 2 \ArcCoords){0.04}
\rput(!9 2 \ArcCoords 0.1 add){\small $\xi_1$}
\pscircle(!1 2 \ArcCoords){0.04}
\rput(!1 2 \ArcCoords 0.1 add){\small $\xi_2$}
\psset{fillstyle=solid,fillcolor=black}
\pscircle(!3 2 4 \SliceCoords){0.04}
\rput(!3 2 4 \SliceCoords 0.1 add){\small $C$}
\end{pspicture}
\end{tabular}
&
\begin{tabular}{c}
\psset{yunit=3cm, xunit=5cm}
\begin{pspicture}(0,0)(1,1)
\psset{linecolor=gray,linestyle=solid}
\psline(0,0)(1,0)
\psset{linecolor=black,linestyle=dotted}
\psline(!2 1 \ArcCoords)(!2 1 \ArcCoords pop 0)
\psset{linecolor=gray,linestyle=solid}
\psset{linestyle=dashed,dash = 3pt 2pt, linecolor=gray}
\psset{linecolor=black,linestyle=solid}
\DrawArc
\psset{fillstyle=solid,fillcolor=white}
\pscircle(!2 1 \ArcCoords){0.04}   
\rput(!2 1 \ArcCoords 0.1 add){\small $\xi$}
\rput(!2 1 \ArcCoords -0.1 0.1 \Trans){\small ({\em i})}
\psset{fillstyle=solid,fillcolor=black}
\pscircle(!6 12 13 \SliceCoords){0.04}
\rput(!6 12 13 \SliceCoords 0.05 0.065 \Trans){\small $C$}
\rput(!6 12 13 \SliceCoords -0.08 0.065 \Trans){\small ({\em ii})}
\rput(0.2,-0.15){\small $\pi^{*}L$}
\rput(0.05,-0.15){\small ({\em iii})}
\psset{linecolor=gray,linestyle=solid,fillstyle=none}
\psline(0.2,0.04)(0.2,-0.04)
\end{pspicture}
\end{tabular}
\\
({\em a}) & ({\em b}) & ({\em c}) \rule[-0.3cm]{0cm}{0.3cm}\\
\hline
\multicolumn{3}{c}{\small \Fig \label{fig:three-facts}} \rule{0cm}{0.8cm}\\
\multicolumn{3}{c}{\small Three graphical arguments} \\
\end{tabular} \\
\end{centering}

\AlphaList
\begin{enumerate}
\item If $\xi$ is a class on the square-zero cone, then the hyperplane $\xi^{\perp}$ is the tangent line
to the cone at $\xi$.
\PauseEnumerate

This is the well-known fact that for a variety $Q$, given as the zeros of a quadratic form 
$\langle\cdot,\cdot\rangle$ on some vector space, then for any point $x\in Q$, the tangent plane to $Q$ at $x$ 
consists of those vectors $v$ such that $\langle x,v\rangle=0$ (since then $\langle x+\epsilon v, x+\epsilon v\rangle$ 
vanishes to first order in $\epsilon$). 

The classes which intersect $\xi$ positively are below this line, and the ones which intersect $\xi$ negatively are
above. 

\ResumeEnumerate
\item If $C$ is a class with $C^2<0$, then the hyperplane $C^{\perp}$ is spanned by $\xi_1$ and $\xi_2$, where
$\xi_1$, $\xi_2$ are the two points on the square-zero cone whose tangent lines contain $C$.
\PauseEnumerate

By ({\em a}) both $\xi_1$ and $\xi_2$ intersect $C$ in zero, therefore all classes on the plane
spanned by $\xi_1$ and $\xi_2$ intersect $C$ in zero. By reason of dimension, this plane is all of $C^{\perp}$. 
The classes above the line, including $C$ itself, intersect $C$ negatively, and the classes below intersect $C$ 
positively. 

\ResumeEnumerate
\item If $\xi$ is a class on the square-zero cone which is nef ({\em i}), then no effective class $C$ which is to
the right of $\xi$ ({\em ii}) can effect the Seshadri constant of a line bundle $L$ which is to the left 
of $\xi$ ({\em iii}), and similarly with right and left reversed in ({\em ii}) and ({\em iii}). 
\PauseEnumerate

The visual interpretation of \eqref{eqn:ep-r-def-nef} is that one starts at $\pi^{*}L$, and moves upwards in
the direction of $-E$ until $L_{\gamma}$ ($ := \pi^{*}L-\gamma E$) either runs into a plane of the type
$C^{\perp}$ with $C^2<0$, or hits the square-zero cone (e.g.\ see ({\em c}\Old{3}) in Figure \ref{fig:proof-of-c} 
below). In the first case, $\ep_{r}(L)$ is computed by $C$ 
(if $C^{\perp}$ is the first such plane encountered) and in the second case $\ep_{r}(L)$ is the maximum
possible value, $\eta_{r}(L)$. 

With reference to Figure \ref{fig:proof-of-c} below, the argument for ({\em c}) is then that,
since $\xi$ is nef, the class $C$ must be below $\xi^{\perp}$ ({\em c}\Old{1}).  But, this means that
$\xi_1$ and $\xi_2$, the points on the square-zero cone whose tangent lines contain $C$, must both be to the right
of $\xi$ ({\em c}\Old{2}).  Therefore the line spanned by $\xi_1$ and $\xi_2$ exits the square-zero cone to the
right of $\xi$ (at worst at $\xi$ if $C$ is on $\xi^{\perp}$) and so $C^{\perp}$ (for this $C$) cannot
affect the Seshadri constant of $L$ ({\em c}\Old{3}). 

\vspace{0.2cm}
\begin{centering}
\begin{tabular}{|c|c|c|}
\hline
\begin{tabular}{c}
\psset{yunit=3cm, xunit=5cm}
\begin{pspicture}(0,0)(1,1)
\psset{linecolor=gray,linestyle=solid}
\psline(0,0)(1,0)
\psset{linecolor=black,linestyle=dotted}
\psset{linecolor=gray,linestyle=solid}
\psset{linestyle=dashed,dash = 3pt 2pt, linecolor=gray}
\TanToArc{2}{1}{-0.3}{0.6}  
\LabelTanToArc{2}{1}{-0.4}{\small $\xi^{\perp}$} 
\psset{linecolor=black,linestyle=solid}
\DrawArc
\psset{fillstyle=solid,fillcolor=white}
\pscircle(!2 1 \ArcCoords){0.04}   
\rput(!2 1 \ArcCoords 0.1 add){\small $\xi$}
\psset{fillstyle=solid,fillcolor=black}
\pscircle(!6 12 13 \SliceCoords){0.04}
\rput(!6 12 13 \SliceCoords 0.05 0.065 \Trans){\small $C$}
\rput(0.2,-0.15){\small $\pi^{*}L$}
\psset{linecolor=gray,linestyle=solid,fillstyle=none}
\psline(0.2,0.04)(0.2,-0.04)
\end{pspicture}
\end{tabular}
&
\begin{tabular}{c}
\psset{yunit=3cm, xunit=5cm}
\begin{pspicture}(0,0)(1,1)
\psset{linecolor=gray,linestyle=solid}
\psline(0,0)(1,0)
\psset{linecolor=black,linestyle=dotted}
\psset{linecolor=gray,linestyle=solid}
\psset{linestyle=dashed,dash = 3pt 2pt, linecolor=gray}
\TanToArc{9}{8}{-0.3}{0.5}  
\TanToArc{2}{9}{-0.5}{0.3}  
\psset{linecolor=black,linestyle=solid}
\DrawArc
\psset{fillstyle=solid,fillcolor=white}
\pscircle(!2 1 \ArcCoords){0.04}   
\rput(!2 1 \ArcCoords 0.1 add){\small $\xi$}
\pscircle(!9 8 \ArcCoords){0.04} 
\rput(!9 8 \ArcCoords 0.1 add){\small$\xi_1$} 
\pscircle(!2 9 \ArcCoords){0.04} 
\rput(!2 9 \ArcCoords 0.1 add){\small $\xi_2$}
\psset{fillstyle=solid,fillcolor=black}
\pscircle(!6 12 13 \SliceCoords){0.04}
\rput(!6 12 13 \SliceCoords 0.05 0.065 \Trans){\small $C$}
\rput(0.2,-0.15){\small $\pi^{*}L$}
\psset{linecolor=gray,linestyle=solid,fillstyle=none}
\psline(0.2,0.04)(0.2,-0.04)
\end{pspicture}
\end{tabular}
&
\begin{tabular}{c}
\psset{yunit=3cm, xunit=5cm}
\begin{pspicture}(0,0)(1,1)
\psset{linecolor=gray,linestyle=solid}
\psline(0,0)(1,0)
\psset{linecolor=black,linestyle=dotted}
\psset{linecolor=gray,linestyle=solid}
\SecToArc{9}{8}{2}{9}{-0.3}{1.3}
\LabelSecToArc{9}{8}{2}{9}{1.5}{\tiny $C^{\perp}$}
\psset{linestyle=dashed,dash = 3pt 2pt, linecolor=gray}
\psset{linecolor=black,linestyle=solid}
\DrawArc
\psset{fillstyle=solid,fillcolor=white}
\pscircle(!2 1 \ArcCoords){0.04}   
\rput(!2 1 \ArcCoords 0.1 add){\small $\xi$}
\pscircle(!9 8 \ArcCoords){0.04} 
\pscircle(!2 9 \ArcCoords){0.04} 
\psset{fillstyle=solid,fillcolor=black}
\pscircle(!6 12 13 \SliceCoords){0.04}
\rput(!6 12 13 \SliceCoords 0.05 0.08 \Trans){\small $C$}
\rput(0.2,-0.15){\small $\pi^{*}L$}
\psset{linecolor=lightgray,linestyle=solid,fillstyle=none,arrows=->}
\psline(0.2,0)(0.2,0.4)
\rput(0.27,0.4){\small\color{gray} $L_{\gamma}$}
\psset{linecolor=gray,linestyle=solid,fillstyle=none,arrows=-}
\psline(0.2,0.04)(0.2,-0.04)
\end{pspicture}
\end{tabular}
\\
({\em c\Old{1}}) & ({\em c}\Old{2}) & ({\em c}\Old{3}) \rule[-0.3cm]{0cm}{0.3cm}\\
\hline
\multicolumn{3}{c}{\small \Fig \label{fig:proof-of-c}} \rule{0cm}{0.8cm}\\
\multicolumn{3}{c}{\small Argument for Figure \ref{fig:three-facts}({\em c})}  \\
\end{tabular} \\
\end{centering}

\bpoint{Reasons for interest in the square-zero cone} 
The square-zero cone is a natural upper (respectively lower) bound for the nef cone (respectively the effective 
cone).   The nef cone can extend at most up to the square-zero cone, although it is not clear how close it can get,
and the effective cone extends past the square-zero cone, although it is not clear how far. 

Second, if $\xi$ is a class on the square-zero cone which is nef, then not only is $\xi$ an example of an extreme
nef class (one which reaches the maximum possible boundary), but, by 
\S\ref{sec:graphical-arguments}({\em c}) above, $\xi$ also splits the problem
of understanding the nef and effective cones into two pieces.   Essentially, there is no information transfer
across the dotted line in Figure \ref{fig:three-facts}({\em c}); knowledge about nef or effective classes on
one side does not allow one to conclude anything about nef or effective classes on the other side. 
An exception to this principle is when one can use $T_r$ to transport information from one part of the cone to another.

Finally, we note that a square-zero class $\xi$ which is nef is not only an example of an extreme nef class, it
is also an example of an extreme effective class.  If $\xi$ is nef, the effective cone must lie below the
tangent line at $\xi$ as in Figure \ref{fig:three-facts}({\em a}), and so at $\xi$ the effective cone is pinched 
down to $\xi$.  I.e., at such a point the boundaries of the nef and effective cones coincide. 

\bigskip

\bpoint{A heuristic argument}
\label{sec:heuristic-argument} 
In this section we provide an argument, with several gaps, which suggests the following principle~:

\smallskip
\begin{centering} 
\begin{minipage}{0.95\textwidth}
If $L$ is an outer bundle on $\PP^1\times\PP^1$, and if $\ep_r(L)\neq \eta_r(L)$, then the Seshadri constant of $L$ 
is computed by a $(-1)$ curve (equivalently, by the symmetrization of a $(-1)$-curve). 
\end{minipage}\\
\end{centering}

\newpage
\hfill
\begin{tabular}{c}
\begin{tabular}{c}
\psset{yunit=3cm, xunit=5cm}
\begin{pspicture}(0,-0.2)(1,0.8)
\psset{linecolor=gray}
\psset{linecolor=black}
\psset{linecolor=gray}
\psline(!1 0 0 \SliceCoords)(!0 1 0 \SliceCoords) 
\psline(!10 \AlphaVal 1 \ArcCoords exch pop 0 exch)(!10 \AlphaVal 1 \ArcCoords exch pop 1 exch)
\SLabel{0}{\tiny $0$}
\SLabel{1 9 \AlphaVal div}{\tiny $\beta_r$} 
\SLabel{9 \AlphaVal }{\tiny $\alpha_r$} 
\SLabelInf{\tiny $\infty$} 
\psset{linecolor=black,fillstyle=none}
\DrawArc
\psset{fillstyle=solid,fillcolor=white,linecolor=black,linestyle=solid,linecolor=black}
\pscircle(!10 \AlphaVal 1 \ArcCoords){0.04}
\rput(!10 \AlphaVal 1 \ArcCoords 0.1 add){$v_{\beta}$}
\pscircle(!1 10 \AlphaVal \ArcCoords){0.04}
\rput(!1 10 \AlphaVal \ArcCoords 0.1 add){$v_{\alpha}$}
\rput(!1 10 \AlphaVal \ArcCoords 0.4 0 \Trans){\tiny $K_X^{\perp}$}
\rput(0.8,-0.15){\small $\pi^{*}L$}
\psset{linecolor=lightgray,linestyle=solid,fillstyle=none,arrows=->}
\psline(0.8,0)(0.8,0.4)
\end{pspicture} \\
\end{tabular} \\
\small \Fig \label{fig:heurstic-argument} \\
\small Heuristic argument  \\
\end{tabular}

\vspace{-3.7cm}
\parshape 1 0cm 10.5cm
The argument is the following.  Consider $L_{\gamma}$ ($:=\pi^{*}L-\gamma E$) for increasing $\gamma$.
Since $L$ is an outer bundle, $L_{\gamma}$ exits the square-zero cone before it crosses the line $K_{X}^{\perp}$. 
If $\ep_r(L)\neq\eta_r(L)$ then $\ep_r(L)$ is computed by some irreducible curve $C'$, with symmetrization $C$
(as in \S\ref{sec:V_r-def}).  By \cite[Theorem 2.6.2({\em f}\hspace{0.5mm})]{Dio} 
$C$ must be quite close to the square-zero cone. 

\parshape 2 0cm 10.5cm 0cm \textwidth
Thus the line $C^{\perp}$ (as in Figure \ref{fig:three-facts}({\em b})) must be quite small, and so $C$ quite
close to the point where $L_{\gamma}$ exits the square-zero cone.   This suggests that $C$ will also be below
the line $K_{X}^{\perp}$, and so $K_{X}$-negative.  If so, then $C'$ is also $K_X$-negative (all curves in the
symmetrization have the same intersection with $K_X$). 

Therefore we have $(C')^2\leq -1$ and $C'\cdot K_{X}\leq -1$.  On a smooth irreducible surface $X$, and with $C'$
an irreducible curve, one always has $(C'+K_{X})\cdot C'\geq -2$.
Thus both inequalities above are equalities, and $C'$ is a $(-1)$-curve. \epf 

\point Despite the gaps in the argument above, the conclusion is correct.   When $r$ is even we show in 
Corollary \ref{cor:Outer-nef} that both $v_{\alpha_r}$ and $v_{\beta_r}$ are nef classes.  As a result, 
first, by \S\ref{sec:graphical-arguments}({\em c}) only curves $C$ whose slope is outside 
$[\beta_r,\alpha_r]$ can affect the Seshadri constant of an outer bundle.  Second, since $v_{\alpha_r}$ and
$v_{\beta_r}$ are nef, any symmetric effective curve $C$ must be below the tangent lines to $v_{\alpha_r}$ 
and $v_{\beta_r}$ (\S\ref{sec:graphical-arguments}({\em a})), and therefore a curve $C$ with slope 
outside $[\beta_r,\alpha_r]$ is $K_X$-negative, as suggested in \S\ref{sec:heuristic-argument}.

In the case that $r$ is odd (and $r\neq 9$)  the classes $v_{\alpha_r}$ and $v_{\beta_r}$ are not nef, and 
we require a different argument to show that $K_X$-positive curves (or $K_X$-null curves) cannot influence the 
Seshadri constant of any outer bundle.   This argument appears in 
\S\ref{sec:Outline-of-outer-odd}--\S\ref{sec:proof-of-outer-odd-theorem}.

\section{Even $r$, $r\geq 10$}
\label{sec:Outer-Even} 

The following result is the key to our analysis of the Seshadri constants of outer bundles when $r$ is even. 

\bpoint{Theorem} 
\label{thm:Tr-transformation} 
Let $\pi\colon X\longrightarrow \PP^1\times \PP^1$ be the blowup of $\PP^1\times\PP^1$ at $r$ general points, $p_1$,
\ldots, $p_r$, with $r$ even.
As in \S\ref{sec:V_r-def} let $V_r\subset H^2(X,\RR)$ be the subspace generated by
the fibre classes $F_1$ and $F_2$, along with the sum of the exceptional divisors $E$.  
Then

\begin{enumerate}
\item The linear transformation $T_r\colon V_r\longrightarrow V_r$ given, in the basis $F_1$, $F_2$, $E$ by 
$$
T_r:= \left[
\begin{array}{crr}
0 & 1 & 0 \\
1 & \frac{r}{2} & r \\
0 & -1 & -1 \rule{0cm}{0.45cm}\\
\end{array}
\right]
$$
is an automorphism of $V_r$, preserving the intersection form. 

\smallskip
\item The eigenvalues of $T_r$ are $\alpha_r$, $\beta_r$, and $1$, with respective eigenvectors 
$v_{\alpha_r}$, $v_{\beta_r}$, and $K_X$, where $v_{\alpha_r}$ and $v_{\beta_r}$ are the vectors given in 
\eqref{eqn:va-vb-v1-def}.

\smallskip
\item If $\xi\in V_r$ is any class, then $\xi$ is nef, or effective, or represents
a curve with $s$ irreducible components if and only if $T_r(\xi)$ is respectively nef, effective, or represents
a curve with $s$ irreducible components.
\end{enumerate}

Thus, by ({\em c}), when $r$ is even $T_r$ induces an automorphism of the nef and effective
cones restricted to $V_r$. 
This automorphism is of infinite order whenever $r\geq 8$. 

\bpf
({\em a}) 
The identity 
$$
\left[
\begin{array}{crr}
0 & 1 & 0 \\
1 & \frac{r}{2} & r \\
0 & -1 & -1 \rule{0cm}{0.45cm}\\
\end{array}
\right]^{t}
\left[
\begin{array}{ccr}
0 & 1 & 0 \\
1 & 0 & 0 \\
0 & 0 & -r \\
\end{array}
\right]
\left[
\begin{array}{crr}
0 & 1 & 0 \\
1 & \frac{r}{2} & r \\
0 & -1 & -1 \rule{0cm}{0.45cm}\\
\end{array}
\right]
=
\left[
\begin{array}{ccr}
0 & 1 & 0 \\
1 & 0 & 0 \\
0 & 0 & -r \\
\end{array}
\right],
$$
where $t$ denotes matrix transpose, shows that $T_r$ preserves the intersection form on $V_r$. 

({\em b}) The characteristic polynomial of $T_r$ is $\left(t^2-\left(\frac{r-4}{2}\right)t+1\right)(t-1)$,
where now $t$ denotes a variable, and therefore the eigenvalues of $T_r$ are $\alpha_r$, $\beta_r$, and $1$. 
It is straightforward to verify that $K_X$ is an eigenvector of eigenvalue $1$. 
Using the identity $\alpha_r^2=\left(\frac{r-4}{2}\right)\alpha_r-1$ as is, and in the form
$2(\alpha_r+1)-r\alpha_r=-2(\alpha_r+1)\alpha_r$, we compute that
$$
\left[
\begin{array}{crr}
0 & 1 & 0 \\
1 & \frac{r}{2} & r \\
0 & -1 & -1 \rule{0cm}{0.45cm}\\
\end{array}
\right]
\left[
\begin{array}{c}
\frac{1}{\alpha_r+1} \\
\frac{\alpha_r}{\alpha_r+1} \rule{0cm}{0.5cm}\\
-\frac{2}{r}\rule{0cm}{0.5cm} \\ 
\end{array}
\right]
=
\left[
\begin{array}{c}
\frac{\alpha_r}{\alpha_r+1} \\
\frac{\left(\frac{r-4}{2}\right)\alpha_r-1}{\alpha_r+1} \rule{0cm}{0.6cm}\\
\frac{2}{r} -\frac{\alpha_r}{\alpha_r+1} \rule{0cm}{0.6cm}\\
\end{array}
\right]
=
\left[
\begin{array}{c}
\frac{\alpha_r}{\alpha_r+1} \\
\frac{\alpha_r^2}{\alpha_r+1} \rule{0cm}{0.6cm}\\
\frac{-2(\alpha_r+1)\alpha_r}{r(\alpha_r+1)} \rule{0cm}{0.6cm}\\
\end{array}
\right]
=
\alpha_r
\left[
\begin{array}{c}
\frac{1}{\alpha_r+1} \\
\frac{\alpha_r}{\alpha_r+1} \rule{0cm}{0.5cm}\\
-\frac{2}{r}\rule{0cm}{0.5cm} \\ 
\end{array}
\right],
$$
and so $T_{r}(v_{\alpha_r})=\alpha_r\,v_{\alpha_r}$.  Similarly $T_r(v_{\beta_r}) = \beta_r\,v_{\beta_r}$.

The real value of the theorem is in part ({\em c}).  
The idea of the argument is that the proper transform of each fibre of type $F_1$ passing through a point $p_i$ is a
$(-1)$-curve.  Blowing down these $r$-different $(-1)$-curves gives another way to realize $X$ as a blowup of 
$\PP^1\times\PP^1$. Comparing the two descriptions as blowups and switching the factors of $\PP^1\times\PP^1$
gives $T_r$. 

To carry this out, 
consider the linear series $|\Osh_{Y}(\frac{r}{2},1)|$ on $Y=\PP^1\times\PP^1$, of dimension $r+1$.
The curves in the series have self intersection $r$, and intersection number $1$ with curves in $|\Osh_{Y}(1,0)|$.
When the $r$ points are general, the series $|\pi^{*}\Osh_{Y}(\frac{r}{2},1)-E|$ (i.e., the proper transforms
of the curves in the series passing through the points) is therefore a basepoint free pencil of curves on $X$. 
The curves have self intersection $0$, and intersection number $1$ with $F_1$. 

The pencils $|F_1|$ and $|\pi^{*}(\Osh_{Y}(\frac{r}{2},1)-E|$  give a birational morphism
$\mu\colon X\longrightarrow\PP^1\times \PP^1$.   This map blows down the curves of class $F_1-E_i$, 
$i=1$,\ldots, $r$, since these classes have intersection number $0$ with the curves in each pencil. 
Since $\mu$ is birational, and since the Picard ranks of $X$ and $\PP^1\times\PP^1$ are $r+2$ and $2$ respectively, 
these are the only curves blown down.

Thus $\mu$ also expresses $X$ as the blowup of $\PP^1\times\PP^1$ at $r$ points, say at $q_1$,\ldots, $q_r$.  
(The map $\mu$ is only really well-defined when we have fixed bases for these pencils; we will do this below.)
Let $F_1'$, $F_2'$ and $E'$ be the pullback of the fibre classes via $\mu$, and the sum of the exceptional divisors
of $\mu$ respectively.  The matrix expressing the change of coordinates on $V_r$ from the second
basis to the first is 
$$T'_r:=
\begin{bNiceMatrix}[first-row,first-col]
& \mbox{\color{gray}\scriptsize $F'_1$} & \mbox{\color{gray}\scriptsize $F'_2$} & \mbox{\color{gray}\scriptsize $E'$} \\
\mbox{\color{gray}\scriptsize $F_1$} & 1 & \phantom{-}\frac{r}{2} & \phantom{-} r \\
\mbox{\color{gray}\scriptsize $F_2$} & 0 & \phantom{-}1 & \phantom{-}0\rule{0cm}{0.4cm}\\
\mbox{\color{gray}\scriptsize $E$} & 0 & -1 & -1 \rule{0cm}{0.4cm}\\
\end{bNiceMatrix}.
$$
Since this matrix represents the identity transformation on $V_r$, albeit between two different bases, a vector
$v$ (in the basis $F_1'$, $F_2'$, and $E'$) is nef, or effective, or represents a curve with $s$ irreducible
components if and only if the vector $T'_r(v)$ (in the basis $F_1$, $F_2$, and $E$) respectively is 
nef, or effective, or represents a curve with $s$ irreducible components. 

If $p_1$, \ldots, $p_r$ are in very general position, then $q_1$,\ldots, $q_r$ are also in
very general position.  We will check this below, but first show how this is enough
to finish the argument.

If $p_1$, \ldots, $p_r$, and $q_1$, \ldots, $q_r$ are in very general position, then 
a class $d_1F_1+d_2F_2-mE$ is nef, effective, or represents a curve with $s$ irreducible components if and only if
the class $d_1F_1'+d_2F_2'-mE'$ has the respective property.

Thus the matrix
$$T''_r:=
\begin{bNiceMatrix}[first-row,first-col]
& \mbox{\color{gray}\scriptsize $F_1$} & \mbox{\color{gray}\scriptsize $F_2$} & \mbox{\color{gray}\scriptsize $E$} \\
\mbox{\color{gray}\scriptsize $F_1$} & 1 & \phantom{-}\frac{r}{2} & \phantom{-} r \\
\mbox{\color{gray}\scriptsize $F_2$} & 0 & \phantom{-}1 & \phantom{-}0\rule{0cm}{0.4cm}\\
\mbox{\color{gray}\scriptsize $E$} & 0 & -1 & -1 \rule{0cm}{0.4cm}\\
\end{bNiceMatrix}
=
\begin{bNiceMatrix}[first-row,first-col]
& \mbox{\color{gray}\scriptsize $F'_1$} & \mbox{\color{gray}\scriptsize $F'_2$} & \mbox{\color{gray}\scriptsize $E'$} \\
\mbox{\color{gray}\scriptsize $F_1$} & 1 & \phantom{-}\frac{r}{2} & \phantom{-} r \\
\mbox{\color{gray}\scriptsize $F_2$} & 0 & \phantom{-}1 & \phantom{-}0\rule{0cm}{0.4cm}\\
\mbox{\color{gray}\scriptsize $E$} & 0 & -1 & -1 \rule{0cm}{0.4cm}\\
\end{bNiceMatrix}
\cdot
\begin{bNiceMatrix}[first-row,first-col]
& \mbox{\color{gray}\scriptsize $F_1$} & \mbox{\color{gray}\scriptsize $F_2$} & \mbox{\color{gray}\scriptsize $E$} \\
\mbox{\color{gray}\scriptsize $F'_1$} & 1 & 0 & 0 \\
\mbox{\color{gray}\scriptsize $F'_2$} & 0 & 1 & 0 \\
\mbox{\color{gray}\scriptsize $E'$}   & 0 & 0 & 1 \\
\end{bNiceMatrix}
$$
gives a linear transformation $V_r\longrightarrow V_r$ in the basis $(F_1,F_2,E)$ preserving
each of those properties.

The transformation $T''_r$ is a reflection. 
The transformation $S_r$ which fixes $E$ and switches $F_1$ and $F_2$ (i.e., the transformation
induced by the automorphism of $\PP^1\times\PP^1$ switching the factors) is also a reflection, and also preserves
all the properties we are interested in.  The product $S_r\cdot T''_r$ is $T_r$, and therefore $T_r$ preserves classes
which are nef, or effective, or which represent curves with $s$ irreducible components, as claimed. 

Thus, to complete the proof of ({\em c}) it is sufficient to verify that $q_1$, \ldots, $q_r$ is in very 
general position if $p_1$, \ldots, $p_r$ are.  This is clear when $r=2$ (any two points not on the same fibres
are in very general position), and so from now on we assume that $r\geq 4$. 

By acting by $\Aut(\PP^1)\times\Aut(\PP^1)$ we may assume that $p_1=([1:0],[1:0])$,
$p_2=([0:1],[0:1])$, and $p_3=([1:1],[1:1])$.  Similarly we may choose a basis for the pencils 
$|F_1|$ and $|\pi^{*}(\Osh_{Y}(\frac{r}{2},1)-E|$ so that under the map $\mu$, $q_1=([1:0],[1:0])$, $q_2=([0:1],[0:1])$,
and $q_3=([1:1],[1:1])$, where $q_i$ is the image of $F_1-E_i$.   This choice is enough to fix the map $\mu$ uniquely. 

Let $U\subset (\PP^1\times\PP^1)^{r-3}$ be the Zariski open subset set of the configuration space of $r-3$ points 
$(p_4,\ldots, p_r)$ so that $p_1$,\ldots, $p_r$ (with $p_1$, $p_2$, and $p_3$ as above) are distinct,
and such that the curves in  $|\Osh_{\PP^1\times\PP^1}(\frac{r}{2},1)|$ passing through
$p_1$,\ldots, $p_r$ form a pencil whose generic member is irreducible.   

By $T_r'$ above, $\frac{r}{2}F_1'+F_2'-E' = F_2$, and thus the points 
$q_1$,\ldots, $q_r$ (with $q_1$, $q_2$, $q_3$ also as above) are sufficiently independent 
so that the linear series $|\frac{r}{2}F_1'+F_2'-E'|$ is a basepoint-free pencil of curves, generically irreducible.

Thus the process $(p_1,\ldots, p_r) \mapsto (q_1,\ldots, q_r)$ induces a map 
$$I\colon (\PP^1\times\PP^1)^{r-3}\longrightarrow (\PP^1\times\PP^1)^{r-3}$$
such that $I(U)\subseteq U$.  Moreover, the identity 
$\frac{r}{2}F_1'+F_2'-E' = F_2$ shows that $I(q_4,\ldots, q_r) = (p_4,\ldots, p_r)$, i.e,. that $I$ is an involution. 

Therefore given a family $V_n\subset U$, $n\in \NN$, of proper closed subsets of $U$, the set 
$$W:=\bigcup_{n\in \NN} \left(V_n\cup I(V_n)\right)$$ 
is a countable union of proper closed subsets of $U$, stable under $I$.  Thus if $(p_4,\ldots, p_r)\not\in W$, 
then $(q_4,\ldots, q_r)=I(p_4,\ldots, p_r)\not\in W$.  Therefore, if
$p_1$,\ldots, $p_r$ are very general, then so are $q_1$,\ldots, $q_r$. 
This finishes the proof of ({\em c}).

When $r=8$, $T_{8}$ is a unipotent matrix consisting of a single Jordan block.  When $r\geq 10$,
$\alpha_r$ is a real number with $\alpha_r>1$.  Thus $T_{r}$ is of infinite order for all $r\geq 8$.  
\epf

\bpoint{Remarks} 
\label{sec:remarks-on-Tr}

\RemNum The process of blowing down, and switching factors in the proof of 
Theorem \ref{thm:Tr-transformation}({\em c}) gives an automorphism $\tilde{T}_r$ of all of $H^2(X,\RR)$ 
(or $H^2(X,\ZZ)$), and not just $V_r$.
In $H^2(X,\RR)$ the orthogonal complement to $V_r$ is the subspace $0\cdot F_1 + 0\cdot F_2 -\sum_{i=1}^{r} m_i E_i$ 
with $\sum_{i} m_i=0$.  On this subspace one can check that $\tilde{T}_r$ acts as multiplication by $-1$.

\RemNum When $r\geq 9$, $\PP^2$ blown up at $r$ general points (and $\PP^1\times\PP^1$ blown up at an odd number of 
points) 
has many such infinite order Cremona ``relabelling'' automorphisms acting on $H^2$ of the blowup.   
The great advantage in the case of $\PP^1\times\PP^1$ blown up at an even number of points is that $\tilde{T}_r$ 
preserves the subspace of equal multiplicity curves.   No such automorphisms exist in the cases of $\PP^2$ blown up 
at $r$ general points, nor for $\PP^1\times\PP^1$ blown up at an odd number of points.   

In the case of $\PP^2$, the corresponding 
space $V_r$ is two dimensional, spanned by the hyperplane class $H$ and $E$. Since $H$ is on the boundary of the nef
cone, any such automorphism has to take $H$ to $H$.  But to preserve the intersection form, the automorphism must
now take $E$ to $\pm E$, and in order to preserve the nef cone, it must take $E$ to $+E$.  Thus, on $\PP^2$ blown
up at any number of points, any such automorphism which preserves the equal multiplicity subspace acts as the 
identity on $V_r$.
We will see in \S\ref{sec:no-automorphisms-when-r-is-odd} that, other than switching the factors, there is no such 
automorphism for $\PP^1\times\PP^1$ blown up at an odd number of points.

\RemNum Parts ({\em a}) and ({\em b}) of Theorem \ref{thm:Tr-transformation} still hold when $r$ is odd.
However, since $T_r$ does not preserve nef or effective classes when $r$ is odd, this transformation is meaningless
for our problem.

\bpoint{Theorem (nef cone for outer bundles, $r$ even)} 
\label{thm:Outer-nef}
Let $X$ be the blowup of $\PP^1\times\PP^1$ at $r$ general points, with $r\geq 10$ even. Then
the nef cone in $V_{r}$, restricted to the half plane $K_{X}^{\leq 0}$ is spanned by $v_{\alpha_r}$, $v_{\beta_r}$,
and the classes $T^{n}_r(F_2)$, with $n\in \ZZ$. 

\bpf
For $n\in \ZZ$ set $\xi_{n}:=T_r^{n}(F_2)$.  Since $F_2$ is a nef, square-zero class, intersecting $K_X$ negatively,
by Theorem \ref{thm:Tr-transformation} each of the $\xi_n$ are also nef square-zero classes intersecting $K_X$ 
negatively. 

Since $F_2\cdot v_{\beta_r} = \frac{\alpha_r}{\alpha_r+1}$, we conclude from the table in 
Figure \ref{fig:intersection-table} that when writing $F_2\in V_r$ in the basis  $v_{\alpha_r}$, $v_{\beta_r}$, and
$K_X$, the coefficient of $v_{\alpha_r}$ is $\frac{r}{r-8}\cdot \frac{\alpha_r}{\alpha_r+1}$.  Thus, since
$\alpha_r$ is the dominant eigenvalue of $T_r$, we have 
$$v_{\alpha_r} = \frac{(r-8)(\alpha_r+1)}{r\cdot \alpha_r} \lim_{n\to\infty} \frac{1}{\alpha_r^n}\,\xi_n,$$
and so $v_{\alpha_r}$ is nef. Similarly, $v_{\beta_r}$ is nef.  We have already checked that both classes are 
square-zero.

Therefore the intersection of the nef cone and $K_{X}^{\leq 0}$ contains the convex hull of 
$v_{\alpha_r}$, $v_{\beta_r}$, and the classes $\xi_{n}$, $n\in \ZZ$, as in 
Figure \ref{fig:K_X-neg-nef-cone}({\em a}) below.

\vspace{0.2cm}
\begin{centering}
\begin{tabular}{|cc|cc|}
\hline
\rule{1.0cm}{0cm}
\begin{tabular}{c}
\psset{yunit=3cm, xunit=5cm}
\begin{pspicture}(0,-0.2)(1,0.8)
\psset{linecolor=gray}
\psset{linecolor=black}
\psset{linecolor=gray}
\psline(!1 0 0 \SliceCoords)(!0 1 0 \SliceCoords) 
\psline(!3 1 \ArcCoords exch pop 0 exch)(!1 3 \ArcCoords exch pop 1 exch)
\rput(!3 1 \ArcCoords -0.4 0 \Trans){\tiny $K_X^{\perp}$}
\psset{fillcolor=vlgray,linecolor=gray,fillstyle=solid}
\pspolygon(!0 1 \ArcCoords)(!1 12 \ArcCoords)(!7 48 \ArcCoords)(!37 192 \ArcCoords)(!175 768 \ArcCoords)(!781 3072 \ArcCoords)(!1 3 \ArcCoords)(!3 1 \ArcCoords)(!3072 781 \ArcCoords)(!768 175 \ArcCoords)(!192 37 \ArcCoords)(!48 7 \ArcCoords)(!12 1 \ArcCoords)(!1 0 \ArcCoords)
\psset{linecolor=black,fillstyle=none}
\DrawArc
\psset{fillstyle=solid,fillcolor=white,linecolor=black,linestyle=solid,linecolor=black}
\pscircle(!3 1 \ArcCoords){0.04} 
\pscircle(!1 3 \ArcCoords){0.04} 
\pscircle(!175 768 \ArcCoords){0.04} 
\pscircle(!37 192 \ArcCoords){0.04} 
\pscircle(!7 48 \ArcCoords){0.04} 
\pscircle(!1 12 \ArcCoords){0.04} 
\pscircle(!0 1 \ArcCoords){0.04} 
\pscircle(!1 0 \ArcCoords){0.04} 
\pscircle(!12 1 \ArcCoords){0.04} 
\pscircle(!48 7 \ArcCoords){0.04} 
\pscircle(!192 37 \ArcCoords){0.04} 
\pscircle(!768 175 \ArcCoords){0.04} 
\rput(!0 1 \ArcCoords 0.06 0 \Trans){\tiny $\xi_0$} 
\rput(!1 12 \ArcCoords 0.06 0 \Trans){\tiny $\xi_1$} 
\rput(!7 48 \ArcCoords 0.06 0 \Trans){\tiny $\xi_2$} 
\rput{-37}(!37 192 \ArcCoords 0.02 0.03 \Trans){\tiny $\cdots$} 
\rput{37}(!192 37 \ArcCoords -0.02 0.03 \Trans){\tiny $\cdots$} 
\rput(!1 3 \ArcCoords 0 0.07 \Trans){\tiny $v_{\alpha}$} 
\rput(!3 1 \ArcCoords 0 0.07 \Trans){\tiny $v_{\beta}$} 
\rput(!1 0 \ArcCoords -0.07 0 \Trans){\tiny $\xi_{-1}$} 
\rput(!12 1 \ArcCoords -0.07 0 \Trans){\tiny $\xi_{-2}$} 
\rput(!48 7 \ArcCoords -0.07 0 \Trans){\tiny $\xi_{-3}$} 
\end{pspicture} 
\end{tabular}
& & &
\begin{tabular}{c}
\psset{yunit=3cm, xunit=5cm}
\begin{pspicture}(0,-0.2)(1,0.8)
\psset{linecolor=gray}
\psset{linecolor=black}
\psset{linecolor=gray}
\psline(!1 0 0 \SliceCoords)(!0 1 0 \SliceCoords) 
\psset{fillcolor=vlgray,linecolor=gray,fillstyle=solid}
\pspolygon(!0 1 \ArcCoords)(!1 12 \ArcCoords)(!7 48 \ArcCoords)(!37 192 \ArcCoords)(!175 768 \ArcCoords)(!781 3072 \ArcCoords)(!1 3 \ArcCoords)(!3 1 \ArcCoords)(!3072 781 \ArcCoords)(!768 175 \ArcCoords)(!192 37 \ArcCoords)(!48 7 \ArcCoords)(!12 1 \ArcCoords)(!1 0 \ArcCoords)
\psset{linecolor=black}
\SecToArc{48}{7}{12}{1}{0}{1}
\SecToArc{12}{1}{1}{0}{0}{1}
\SecToArc{1}{0}{0}{1}{0}{1}
\SecToArc{0}{1}{1}{12}{0}{1}
\SecToArc{1}{12}{7}{48}{0}{1}
\psset{linestyle=dashed,dash = 3pt 2pt, linecolor=gray}
\TanToArc{0}{1}{-0.3}{0}
\TanToArc{1}{12}{-0.25}{0.25}
\TanToArc{7}{48}{-0.25}{0.25}
\TanToArc{1}{0}{0}{0.3}
\TanToArc{12}{1}{-0.25}{0.25}
\TanToArc{48}{7}{-0.25}{0.25}
\psset{linecolor=black,fillstyle=none,linestyle=solid}
\DrawArc
\psset{fillstyle=solid,fillcolor=white,linecolor=black,linestyle=solid,linecolor=black}
\pscircle(!3 1 \ArcCoords){0.04} 
\pscircle(!1 3 \ArcCoords){0.04} 
\pscircle(!175 768 \ArcCoords){0.04} 
\pscircle(!37 192 \ArcCoords){0.04} 
\pscircle(!7 48 \ArcCoords){0.04} 
\pscircle(!1 12 \ArcCoords){0.04} 
\pscircle(!0 1 \ArcCoords){0.04} 
\pscircle(!1 0 \ArcCoords){0.04} 
\pscircle(!12 1 \ArcCoords){0.04} 
\pscircle(!48 7 \ArcCoords){0.04} 
\pscircle(!192 37 \ArcCoords){0.04} 
\pscircle(!768 175 \ArcCoords){0.04} 
\psset{fillstyle=solid,fillcolor=black}
\pscircle(!0 6 sqrt 2 mul 1 \SliceCoords){0.04}
\rput(!0 6 sqrt 2 mul 1 \SliceCoords 0.06 0 \Trans){\tiny $C_1$}
\pscircle(!0.232497 2.10901 1 \SliceCoords){0.04}
\rput(!0.232497 2.10901 1 \SliceCoords 0.06 0 \Trans){\tiny $C_2$}
\pscircle(!6 sqrt 2 mul 0 1 \SliceCoords){0.04}
\rput(!6 sqrt 2 mul 0 1 \SliceCoords -0.08 0 \Trans){\tiny $C_{-1}$}
\pscircle(!2.10901 0.232497 1 \SliceCoords){0.04}
\rput(!2.10901 0.232497 1 \SliceCoords -0.08 0 \Trans){\tiny $C_{-2}$}
\end{pspicture} 
\end{tabular}
\rule{0.5cm}{0cm}
\\
\rule{0.6cm}{0cm}({\em a}) & \rule{1cm}{0cm} & \rule{1cm}{0cm} & ({\em b}) \rule[-0.3cm]{0cm}{0.3cm}\rule{0.3cm}{0cm}\\
\hline
\multicolumn{4}{c}{\small \Fig \label{fig:K_X-neg-nef-cone}} \rule{0cm}{0.8cm}\\
\multicolumn{4}{c}{\small Argument for Theorem \ref{thm:Outer-nef}} \\
\end{tabular} \\
\end{centering}

\vspace{0.3cm}
Next, for $n\in \ZZ$ set $C_n:=T_r^{n}(E)$.  Since $F_1\cdot E =0$ and $F_2\cdot E=0$, i.e., since $\xi_{-1}\cdot C_0=0$
and $\xi_{0}\cdot C_0=0$, and since $T_r$ preserves intersections, we have $\xi_{n-1}\cdot C_n=0$ and
$\xi_{n}\cdot C_n=0$ for all $n\in \ZZ$.

By \ref{sec:graphical-arguments}({\em b}) this means that for each $n$, the nef cone cannot pass the line spanned by
$\xi_{n-1}$ and $\xi_n$.  Thus (as illustrated in Figure \ref{fig:K_X-neg-nef-cone}({\em b}))
the intersection of the nef cone with $K_{X}^{\leq 0}$ can be no larger than the
cone spanned by $v_{\alpha_r}$, $v_{\beta_r}$, and the $\xi_{n}$, $n\in \ZZ$.  \epf

\bpoint{Corollary (Case of even $r$ in Theorem \ref{thm:Seshadri-for-outer-bundles})} 
\label{cor:Outer-nef}
Let $X$ be the blowup of $\PP^1\times\PP^1$ at $r$ general points, with $r\geq 10$ even.

\AlphaList
\begin{enumerate}
\item The classes $v_{\alpha_r}$ and $v_{\beta_r}$ are nef.
\item For any ample bundle $L=\Osh_{Y}(e_1,e_2)$ with $\frac{e_2}{e_1}\not\in [\beta_r,\alpha_r]$, 
the Seshadri constant of $L$ is computed by $C_n$, where $n$ is an
integer such that $L_{\eta_r(L)}$ is between $\xi_{n}$ and $\xi_{n-1}$ on the square zero cone. 
\end{enumerate}

\hfill
\begin{tabular}{c}
\psset{yunit=9cm, xunit=15cm}
\begin{pspicture}(0.7,0.2)(1,0.6)
\psset{linecolor=gray}
\psset{linecolor=black}
\psset{linecolor=gray}
\psset{linecolor=black}
\SecToArc{1}{18}{1}{6}{0}{1}
\psset{linestyle=dashed,dash = 3pt 2pt, linecolor=gray}
\TanToArc{1}{18}{-0.15}{0.10}
\TanToArc{1}{6}{-0.07}{0.13}
\psclip[linecolor=vlgray,fillstyle=none,linestyle=none]{\pspolygon(0.7,0.2)(1,0.2)(1,0.6)(0.7,0.6)}
\psset{fillcolor=vlgray,linecolor=gray,fillstyle=solid}
\pspolygon(!0 1 \ArcCoords)(!1 18 \ArcCoords)(!1 6 \ArcCoords)(!1 4 \ArcCoords)(!1 2 \ArcCoords)(!2 1 \ArcCoords)(!1 0 \ArcCoords)
\psset{linecolor=black,fillstyle=none,linestyle=solid}
\DrawPartArc{0.62}{0.9}
\psset{linecolor=darkgray,fillstyle=none,linestyle=solid}
\psline[arrows=->](!1 8 0 \SliceCoords)(!1 8 3 \SliceCoords)
\rput(!1 8 3 \SliceCoords -0.02 0 \Trans){\tiny $L_{\gamma}$} 
\endpsclip
\psset{linecolor=black,fillstyle=none,linestyle=solid}
\psset{fillstyle=solid,fillcolor=white,linecolor=black,linestyle=solid,linecolor=black}
\pscircle(!1 18 \ArcCoords){0.04} 
\rput(!1 18 \ArcCoords 0.03 0 \Trans){\tiny $\xi_{n-1}$} 
\pscircle(!1 6 \ArcCoords){0.04} 
\rput(!1 6 \ArcCoords 0.02 0 \Trans){\tiny $\xi_{n}$}
\pscircle(!1 4 \ArcCoords){0.04} 
\rput(!1 4 \ArcCoords 0.03 0 \Trans){\tiny $\xi_{n+1}$} 
\psset{fillstyle=solid,fillcolor=black}
\pscircle(!0.211324 2.19615 1 \SliceCoords){0.04}
\rput(!0.211324 2.19615 1 \SliceCoords 0.02 0 \Trans){\tiny $C_n$}
\psset{fillstyle=solid,fillcolor=lightgray,linecolor=black,linestyle=solid,linecolor=black}
\pscircle(!1 8 \ArcCoords){0.04}
\rput(!1 8 \ArcCoords 0.01 0.02 \Trans){\tiny $L_{\eta}$}
\end{pspicture} \\
\small \Fig \label{fig:Closeup-of-Cn} \\
\end{tabular}

\vspace{-4.5cm}
\parshape 1 0cm 11cm
\bpf
Part ({\em a}) was established in the proof of Theorem \ref{thm:Outer-nef}, and is included here for reference.
Part ({\em b}) is immediate from the visual interpretation of \eqref{eqn:ep-r-def-nef} 
(as in \S\ref{sec:graphical-arguments}({\em c})), and the description of the nef cone in Theorem \ref{thm:Outer-nef}.

\parshape 6 0cm 11cm 0cm 11cm 0cm 11cm 0cm 11cm 0cm 11cm 0cm \textwidth
In more detail, the line segment $L_{\gamma}$, $\gamma\geq 0$, meets the square-zero cone when 
$\gamma=\eta_{r}(L)$.  By Theorem \ref{thm:Outer-nef} if $L$ is an outer bundle the line $L_{\gamma}$ exits the 
nef cone through a secant line 
spanned by $\xi_{n}$ and $\xi_{n-1}$ for some $n$ (unique 
unless $L_{\eta_r(L)}$ is one of the $\xi_{m}$), as illustrated in Figure \ref{fig:Closeup-of-Cn}.
Thus $L_{\eta_r(L)}$ is between $\xi_{n}$ and $\xi_{n-1}$, and can be used to identify $n$.

The proof of Theorem \ref{thm:Outer-nef} shows that this secant line is also the line $C_{n}^{\perp}$, and 
therefore $C_{n}$ computes the Seshadri constant of $L$. \epf

\bpoint{Computing the Seshadri constant of an outer bundle, $r$ even}  
\label{sec:computing-Seshadri-r-even}
Given a line bundle $L=\Osh_{Y}(e_1,e_2)$ in order 
to use Corollary \ref{cor:Outer-nef}({\em b}) to compute $\ep_{r}(L)$ one needs to find 
the correct value of $n$, and then find $C_n$.

Since $C_n=T_r^{n}(F_2)$, the coordinates of the $C_n$ satisfy the recursion relation coming from the
characteristic polynomial of $T_r$, i.e., satisfy \eqref{eqn:sequence-recursion}.  
It follows that, as in \S\ref{sec:intro-seshadri-outer}, if one defines sequences 
$\{p_n\}_{n\in \ZZ}$ and $\{m_n\}_{n\in \ZZ}$ by 
$p_{-1}=0$, $p_{0}=0$, $p_{1}=r$, and $m_{-1}=1$, $m_{0}=-1$, $m_{1}=1$, and the recursion 
\eqref{eqn:sequence-recursion}, that for all $n$ the class of $C_n$ is $(p_{n-1},p_{n},-m_{n})$.
If $L=\Osh_{Y}(e_1,e_2)$, $\ep_r(L)$ is then computed by \eqref{eqn:Seshadri-even-n}, for the right value of $n$.

\point 
\label{sec:methods-to-find-n} 
One can use several methods to find $n$.  If one considers the square zero cone in the region that 
$d_1$, $d_2\geq 0$ (including $m<0$, i.e., not just in the octant $m\geq 0$ from \S\ref{sec:airplane-hanger}) and 
removes the rays spanned by $v_{\alpha_r}$ and $v_{\beta_r}$, then up to scaling by $\RR_{>0}$ what remains
is the union of two disjoint open intervals.  The map $\phi_r$ defined by \eqref{eqn:eta-r-def} gives a 
diffeomorphism of each of these intervals with $\RR$, converting the action of $T_r$ into addition by $1$, and
converting the operation of swapping the fibre classes into multiplication by $-1$
(see \S\ref{sec:applications-to-square-zero-cone}).

It is straightforward to check that $\phi_r(F_2)=\frac{1}{2}$, from which it follows that 
$\phi_r(\xi_n) = n+\frac{1}{2}$ for all $n\in \ZZ$.   Thus, if $L_{\eta_r(L)}$ is between $\xi_n$ and
$\xi_{n-1}$ on the square zero cone, we have $n+\frac{1}{2} \geq \phi_r(L_{\eta_r(L)}) \geq n-\frac{1}{2}$
and so (as long as $L_{\eta_r(L)}\neq \xi_n$), $n=\lfloor \phi_r(L_{\eta_r(L)})+\frac{1}{2}\rfloor$.
If $L_{\eta_r(L)}=\xi_n$, then this formula produces $n+1$ instead, but in this case 
$\ep_r(L)=\eta_r(L)$, and both $C_{n}$ and $C_{n+1}$ compute this answer.
This is the method given in Theorem \ref{thm:Seshadri-for-outer-bundles}.

\point The previous method has the advantage of giving a formula for $n$, however it is not very useful computationally. 
Evaluating $\phi_r$ correctly requires a high degree of precision in real number calculations, too large to
be of much use in general.  Instead, it is computationally more efficient to find $n$ so that the slope of $L$ is 
between the slopes of $\xi_n$ and $\xi_{n-1}$.

As in \S\ref{sec:full-packings}, define a sequence $\{q_n\}_{n\in \ZZ}$ by $q_{-1}=1$, $q_{0}=0$, $q_{1}=1$,
and the recursion \eqref{eqn:sequence-recursion}.   
Then $\xi_{n}=T_r^{n}(F_2) = (q_{n},q_{n+1},-\sqrt{\frac{2q_{n+1}q_{n}}{r}})$ for all $n\in \ZZ$, and thus
the slope of $\xi_n$ is $\frac{q_{n+1}}{q_{n}}$ (this is the reason for case ({\em ii}) in 
Theorem \ref{thm:full-packings}).

If $L=\Osh_{Y}(e_1,e_2)$, with $e_1\leq e_2$ and $e_1\neq 0$, then the relevant value of $n$ is the smallest
$n\geq 1$ so that $\frac{q_{n+1}}{q_{n}}\leq \frac{e_2}{e_1}$, i.e., the smallest $n\geq 1$ so that 
$e_2q_{n}-e_1q_{n+1}\geq 0$.  Computing the $q_{n}$'s and checking the previous condition only involves
integer arithmetic. 

If instead $e_1>e_2$, one can either use the fact that, by symmetry, $\ep_r(L)=\ep_r(\Osh_{Y}(e_2,e_1))$ and the
method above, or use a similar argument for negative $n$. 

\bpoint{Automorphisms preserving $V_r$} In light of the utility of $T_r$ it is interesting
to ask if there are other integral linear automorphisms of $V_r$ which preserve all aspects of the problem (i.e., which
preserve the intersection form, the nef and effective cones, and the canonical class).  
Let $G$ denote the group of such automorphisms. 

In the proof of Theorem \ref{thm:Tr-transformation}({\em c}) we have seen that $S_r$, the automorphism switching $F_1$
and $F_2$ and fixing $E$, is in $G$.   The following argument shows that $T_r$ and $S_r$ generate $G$. 
Since $S_rT_rS_r^{-1} = T_r^{-1}$ (and $S_r=S_r^{-1}$), effectively this means that, up to switching $F_1$ and $F_2$, 
there are really no linear automorphisms of the problem other than the $T_r^{n}$. 

{\em Proof of claim.} 
Suppose that $g\in G$, and consider $g(\xi_0)$.    Since $g$ preserves the intersection form, the canonical
class, and the nef cone, $g(\xi_0)$ is a square-zero nef class which intersects $K_{X}$ negatively.
By Theorem \ref{thm:Outer-nef} this means that $g(\xi_0)$ must be multiple of one of the $\xi_n$.  Since 
$g$ is an integral linear transformation (i.e., coming from an action on $V_{r,\ZZ}$, the underlying integral
lattice), and since each $\xi_{n}$ is the integral generator on the ray it spans, we conclude that 
$g(\xi_0)=\xi_n$ for some unique $n\in \ZZ$.  Multiplying on the left by $T_{r}^{-n}$ we may assume
that $g(\xi_0)=\xi_0$.

Now consider $g(\xi_{-1})$ and $g(\xi_{1})$.  By the previous reasoning, we must have 
$g(\xi_{-1})=\xi_{i}$ and $g(\xi_{1})=\xi_{j}$ for unique $i$, $j\in \ZZ$.  Since $g$ preserves the intersection
form, we have $\xi_0\cdot \xi_{i} = g(\xi_0)\cdot g(\xi_{-1})=\xi_{0}\cdot \xi_{-1} = 1$
and similarly $\xi_0\cdot\xi_{j}=1$. 

But, it is easy to verify that the only $m$ for which $\xi_{0}\cdot \xi_{m}=1$ are $m=-1$, $1$.  Therefore either
$g(\xi_{-1})=\xi_{-1}$ and $g(\xi_{1})=\xi_{1}$, or $g(\xi_{-1})=\xi_{1}$ and $g(\xi_{1})=\xi_{-1}$.

In the first case, $g$ now fixes $\xi_{-1}$, $\xi_{0}$, and $\xi_{1}$, and so must act as the identity on $V_r$,
since these three classes span $V_r$.  In the second case, $T_rS_r$ is also a transformation which fixes $\xi_0$
and swaps $\xi_{-1}$ and $\xi_{1}$.  Multiplying $g$ on the left by $T_rS_r$ then reduces us to the first case.

Thus, in both cases $g$ is in the group generated by $T_r$ and $S_r$.  \epf

\section{Odd $r$, $r\geq 9$}
\label{sec:Outer-Odd} 

\bpoint{Portrait of the outer nef cone, $r$ odd} 
\label{sec:Outer-odd-r-intro}
Figure \ref{fig:General-outer-picture-odd-r}({\em a}) below shows, for odd $r$, $r\geq 11$,  the $(-1)$-curves, 
or symmetrizations of $(-1)$-curves, which can affect the Seshadri constants of line bundles.   In contrast
to the case when $r$ is even, there are only four such curves; these are labelled $C_1$,\ldots, $C_4$ below.
The picture is somewhat cluttered, so in ({\em b}) and ({\em c}) we show separately the two curves $C_4$ and $C_3$
on the right hand side of ({\em a}).

\hspace{-1.0cm}
\begin{tabular}{|c|c|c|}
\hline
\begin{tabular}{c}
\psset{yunit=3cm, xunit=5cm}
\begin{pspicture}(-0.5,0)(1.6,1)
\psset{linecolor=gray}
\psline(0,0)(1,0)
\psset{linecolor=black}
\psset{linestyle=dashed,dash = 3pt 2pt, linecolor=lightgray}
\TanToArc{1}{0}{-0.1}{0.4}
\TanToArc{0}{1}{-0.4}{0.1}
\psset{linestyle=solid,linecolor=black}
\SecToArc{1}{0}{11}{2}{-0.2}{1.2}
\SecToArc{0}{1}{2}{11}{-0.2}{1.2}
\SecToArc{2}{11 sqrt 1 sub dup mul}{2}{11 sqrt 1 add dup mul}{-0.2}{1.7}
\SecToArc{11 sqrt 1 add dup mul}{2}{11 sqrt 1 sub dup mul}{2}{-0.7}{1.2}
\psset{linecolor=gray}
\SecToArc{11 \AlphaVal}{1}{1}{11 \AlphaVal}{-0.2}{1.2}
\rput(!1 11 \AlphaVal \ArcCoords 0.15 0.1 \Trans){\tiny\color{gray} $K_{X}^{\perp}$}
\psset{linecolor=black,fillstyle=none}
\DrawArc
\psset{fillstyle=solid,fillcolor=white,linecolor=black,linestyle=solid,linecolor=black}
\psset{fillstyle=solid,fillcolor=black}
\pscircle(!11 0 1 11 \RSliceCoords){0.04}
\rput(!9 0 1 9 \RSliceCoords -0.27 0 \Trans){\tiny $C_1:=(r,0,-1)$}
\pscircle(!0 11 1 11 \RSliceCoords){0.04}
\rput(!0 11 1 11 \RSliceCoords 0.27 0 \Trans){\tiny $(0,r,-1)=:C_4$}
\pscircle(!5 1 1 11 \RSliceCoords){0.04}
\rput(!5 1 1 11 \RSliceCoords -0.32 0 \Trans){\tiny $C_2:=(\frac{r-1}{2},1,-1)$}
\pscircle(!1 5 1 11 \RSliceCoords){0.04}
\rput(!1 5 1 11 \RSliceCoords 0.30 0 \Trans){\tiny $(1,\frac{r-1}{2},-1)=:C_3$}
\psset{fillstyle=solid,fillcolor=white}
\pscircle(!1 0 \ArcCoords){0.04} 
\pscircle(!0 1 \ArcCoords){0.04} 
\psset{fillcolor=lightgray}
\pscircle(!11 \AlphaVal 1 \ArcCoords){0.04} 
\rput(!11 \AlphaVal 1 \ArcCoords 0 0.1 \Trans){\tiny $v_{\beta}$} 
\pscircle(!1 11 \AlphaVal \ArcCoords){0.04} 
\rput(!1 11 \AlphaVal \ArcCoords 0 0.1 \Trans){\tiny $v_{\alpha}$} 
\pscircle(!6 1 1 11 \RSliceCoords){0.04}
\pscircle(!1 6 1 11 \RSliceCoords){0.04}
\SLabel{1 6 div}{\tiny $\frac{2}{r+1}$}
\SLabel{6}{\tiny $\frac{r+1}{2}$}
\end{pspicture} 
\end{tabular}
&
\begin{tabular}{c}
\psset{yunit=3cm, xunit=5cm}
\begin{pspicture}(0.6,-0.2)(1.15,1)
\psclip[linestyle=none]{\pspolygon(0.65,-0.2)(1.15,-0.2)(1.15,0.9)(0.65,0.9)}
\psset{linecolor=gray,linestyle=solid}
\psline(0,0)(1,0)
\psset{linecolor=black}
\psset{linestyle=dashed,dash = 3pt 2pt, linecolor=lightgray}
\TanToArc{0}{1}{-0.4}{0.1}
\psset{linestyle=solid,linecolor=black}
\SecToArc{0}{1}{2}{11}{-0.2}{1.2}
\LabelSecToArc{0}{1}{2}{11}{1.4}{\tiny $C_4^{\perp}$}
\psset{linecolor=gray}
\psset{linecolor=black,fillstyle=none}
\DrawArc
\psset{fillstyle=solid,fillcolor=white,linecolor=black,linestyle=solid,linecolor=black}
\psset{fillstyle=solid,fillcolor=black}
\pscircle(!0 11 1 11 \RSliceCoords){0.04}
\rput(!0 11 1 11 \RSliceCoords 0.07 0 \Trans){\tiny $C_4$}
\psset{fillstyle=solid,fillcolor=white}
\pscircle(!0 1 \ArcCoords){0.04} 
\psset{fillcolor=lightgray}
\endpsclip
\end{pspicture} 
\end{tabular}
&
\begin{tabular}{c}
\psset{yunit=3cm, xunit=5cm}
\begin{pspicture}(0.40,-0.2)(1.10,1)
\psclip[linestyle=none]{\pspolygon(0.5,-0.2)(1.15,-0.2)(1.15,0.9)(0.5,0.9)}
\psset{linecolor=gray,linestyle=solid}
\psline(0,0)(1,0)
\psset{linecolor=black}
\psset{linestyle=dashed,dash = 3pt 2pt, linecolor=lightgray}
\psset{linestyle=solid,linecolor=black}
\SecToArc{2}{11 sqrt 1 sub dup mul}{2}{11 sqrt 1 add dup mul}{-0.2}{1.7}
\LabelSecToArc{2}{11 sqrt 1 sub dup mul}{2}{11 sqrt 1 add dup mul}{2.0}{\tiny $C_3^{\perp}$}
\psset{linecolor=gray}
\psset{linecolor=black,fillstyle=none}
\DrawArc
\psset{fillstyle=solid,fillcolor=white,linecolor=black,linestyle=solid,linecolor=black}
\psset{fillstyle=solid,fillcolor=black}
\pscircle(!1 5 1 11 \RSliceCoords){0.04}
\rput(!1 5 1 11 \RSliceCoords 0.07 0 \Trans){\tiny $C_3$}
\psset{fillstyle=solid,fillcolor=white}
\psset{fillcolor=lightgray}
\pscircle(!1 11 \AlphaVal \ArcCoords){0.04} 
\rput(!1 11 \AlphaVal \ArcCoords 0 0.1 \Trans){\tiny $v_{\alpha}$} 
\endpsclip
\end{pspicture} 
\end{tabular} \\
({\em a}) & ({\em b}) & ({\em c}) \rule{0cm}{0.8cm}\\
\hline
\multicolumn{3}{c}{\Fig \label{fig:General-outer-picture-odd-r} \rule{0cm}{0.8cm}} \\
\end{tabular}

The curve $C_4$ of class $(0,r,-1)$ is the union of the proper transforms of the fibres of type $F_2$ through each of 
the $p_i$, i.e, is the union of the $r$ $(-1)$-curves $F_2-E_i$, $i=1$, \ldots, $r$.   The plane $C_4^{\perp}$ 
intersects the square-zero cone in rays spanned by the classes $F_2$ and $(2,r,-2)$.  A picture of $C_4$, 
and $C_4^{\perp}$ appears in Figure \ref{fig:General-outer-picture-odd-r}({\em b}).

The linear series $|\Osh_{Y}(1,\frac{r-1}{2})|$ has dimension $r$, curves in the series have self intersection $r-1$,
and smooth curves in the series are rational. 
When the $r$ points are general, $|\pi^{*}\Osh_{Y}(1,\frac{r-1}{2})-E|$ therefore consists of a single smooth
rational curve of self intersection $-1$, i.e., the class $C_3=(1,\frac{r-1}{2},-1)$ is represented by an irreducible
$(-1)$-curve. 

The plane $C_3^{\perp}$ intersects the square-zero cone in rays spanned by the classes 

\begin{equation}
\label{eqn:C_3-boundaries}
\left(1,\frac{(\sqrt{r}-1)^2}{2},-(1-\frac{1}{\sqrt{r}})\right)
\rule{0.25cm}{0cm}\mbox{and}\rule{0.25cm}{0cm} 
\left(1,\frac{(\sqrt{r}+1)^2}{2},-(1+\frac{1}{\sqrt{r}})\right).
\end{equation}
When $r$ is odd, $r\geq 11$, the class $v_{\alpha_r}$ is strictly above the plane $C_{3}^{\perp}$, and thus (in contrast
to the case when $r$ is even) in those cases neither $v_{\alpha_r}$ nor $v_{\beta_r}$ are nef.

The planes $C_{3}^{\perp}$ and $C_4^{\perp}$ intersect in the ray spanned by the class $(1,\frac{r+1}{2},-1)$, of 
slope $\frac{r+1}{2}$. 

The existence of the effective classes $C_3$ and $C_4$ gives the following inequalities on Seshadri constants for
a line bundle $L= \Osh_{Y}(e_1,e_2)$.

\begin{itemize}
\item If $\frac{e_2}{e_1}\in [\frac{(\sqrt{r}-1)^2}{2},\frac{r+1}{2}]$ then $\ep_{r}(L) \leq 
\frac{(r-1)e_1+2e_2}{2r}$  (inequality imposed by $C_3$);

\smallskip
\item If $\frac{r+1}{2}\leq \frac{e_2}{e_1}$ then $\ep_r(L) \leq e_1$ (inequality imposed by $C_4$).
\end{itemize}

Our goal in this section is to show that the above inequalities, and the symmetric inequalities involving $C_1$ and $C_2$,
are equalities whenever $\frac{e_2}{e_1}\not\in (\beta_r,\alpha_r)$.   That is, our goal is to prove the following
result.

\bpoint{Theorem (Case of odd $r$ in Theorem \ref{thm:Seshadri-for-outer-bundles})} 
\label{thm:Outer-nef-odd}
Suppose that $r$ is odd, $r\geq 9$, that $L=\Osh_{\PP^1\times\PP^1}(e_1,e_2)$ with $e_1,e_2\geq 1$, and
that $L$ is an outer bundle, i.e., that $\frac{e_2}{e_1}\not\in (\beta_r,\alpha_r)$.  Then

\vspace{-0.25cm}
\begin{equation}
\label{eqn:Outer-nef-odd}
\ep_r(L) = \left\{ 
\begin{array}{cl}
e_2 & \mbox{if $\frac{e_2}{e_1}\leq \frac{2}{r+1}$}, \\
\frac{2e_1+(r-1)e_2}{2r} & \mbox{if $\frac{e_2}{e_1} \in [\frac{2}{r+1},\beta_r]$}, \rule{0cm}{0.6cm}\\
\frac{(r-1)e_1+2e_2}{2r} & \mbox{if $\frac{e_2}{e_1} \in [\alpha_r,\frac{r+1}{2}]$}, \rule{0cm}{0.6cm}\\
e_1 & \mbox{if $\frac{r+1}{2}\leq \frac{e_2}{e_1}$}. \rule{0cm}{0.6cm} \\
\end{array}
\right.
\end{equation}

\bpoint{Outline of the argument} 
\label{sec:Outline-of-outer-odd}
The essential claim of the theorem is that there is no effective curve $C$ which, for an outer bundle $L$,
imposes a stronger condition on the Seshadri constant of $L$ than those imposed by $C_1$, \ldots, $C_4$.    

When $r=9$ the picture is slightly different from that of the general case shown in 
Figure \ref{fig:General-outer-picture-odd-r}({\em a}).  The difference is that $v_{\alpha_9}$ (respectively $v_{\beta_9}$)
is on the line $C_3^{\perp}$ (respectively $C_2^{\perp}$).  Specifically, 
$$v_{\alpha_9} = \frac{1}{3}\left( 1, \frac{(\sqrt{9}-1)^2}{2}, -(1-\frac{1}{\sqrt{9}})\right),$$
(i.e., one-third the first class in \eqref{eqn:C_3-boundaries})
and similarly for $v_{\beta_9}$.
In the case $r=9$ we show in Corollary \ref{cor:v-alpha-and-v-beta-nef-for-9}
that $v_{\alpha_9}$ and $v_{\beta_9}$ are nef.   Thus, by the principle of 
\ref{sec:graphical-arguments}({\em c}), no curve $C$ with slope in $(\beta_{9},\alpha_{9})$ can 
affect the Seshadri constant of an outer bundle.  This establishes the theorem for $r=9$. 

For $r\geq 11$ we may make several reductions.  First, by symmetry it is enough to restrict to the case 
that $e_1\leq e_2$.
Second, the effective classes $C$ which could effect a Seshadri constant must satisfy $C^2<0$ and be $K_{X}$-positive 
(any irreducible class, or symmetrization of an irreducible class as in \S\ref{sec:V_r-def}, which is not 
$C_1$,\ldots, $C_4$ must be below the planes $C_1^{\perp}$,\ldots, $C_4^{\perp}$, and to also satisfy $C^2<0$ must 
therefore be above the plane $K_X^{\perp}$). 

We are not able to show that such curves $C$ don't exist.  However, if one does exist, we are
able to show that for any nef bundle $L=\Osh_{Y}(e_1,e_2)$ with  $\frac{e_2}{e_1}\geq \alpha_r$, that 
$$\frac{C\cdot\pi^{*}L}{C\cdot E} \geq \frac{(r-1)e_1+2e_2}{2r} = \frac{C_3\cdot \pi^{*}L}{C_3\cdot E}.$$
That is, we are able to show that there are no such curves $C$ which impose a stronger condition on a bundle of slope 
$\geq\alpha_r$ than that imposed by $C_3$, and this is enough to establish the theorem.

The non-existence of such a $C$ is the result of an estimate on how strong a condition such a $C$ could put on
a bundle of slope $\alpha_r$, combined with an estimate on the size of $\alpha_r$.    
This material appears in 
\S\ref{thm:lower-bound-on-Seshadri-quotient}--\S\ref{prop:t-lower-estimate}.
The concluding arguments of the proof, using these steps, 
appears in \S\ref{sec:proof-of-outer-odd-theorem}.

\bpoint{Proposition (Weak lower bound on multiplicity for $K_X$-positive curves)} 
\label{prop:weak-lower-bound-K_X-pos}
Let $X$ be the blowup of $\PP^1\times\PP^1$ at $r$ general points ($r$ may be even or odd, and
there is no restriction on the size of $r$).  Then there does not exist an effective curve $C$ such that
$C^2<0$, $K_X\cdot C\geq 0$ and all multiplicities of $C$ are in $\{0,1\}$ (i.e., such that the class of $C$
is $d_1F_1+d_2F_2-\sum_{i=1}^{r} m_i E_i$ with each $m_i\in \{0,1\}$). 

\bpf
The vector space of curves of bidegree $(d_1,d_2)$ on $\PP^1\times \PP^1$ with multiplicity $m_1$,\ldots, $m_r$
at $p_1$, \ldots, $p_r$ has expected dimension 
$\max\left(0, (d_1+1)(d_2+1)-\tfrac{1}{2}\sum_{i=1}^{r}m_i^2-m_i\right)$.
Letting $C$ be of class $d_1F_1+d_2F_2 -\sum_{i=1}^{r} m_iE_i$, this is the same as 

\begin{equation}
\label{eqn:Expected-dim}
\max\left(0, \tfrac{1}{2}(C^2-K_X\cdot C)+1\right).
\end{equation}

General points of multiplicity $1$ impose independent conditions.  Therefore, if all the multiplicities are $0$ or $1$,
and if the $p_i$ are general, the dimension of the vector space of such curves is the expected one, as given by
\eqref{eqn:Expected-dim}.  If $C^2<0$ and $K_X\cdot C\geq 0$ then \eqref{eqn:Expected-dim} gives zero 
(if $C^2=-1$ then $K_X\cdot C>0$, since $C^2-K_X\cdot C$ must be even), 
and therefore no such effective curves exist.  \epf

\bpoint{Lemma} 
\label{lem:eta-is-two-means-nef}
Let $\pi\colon X\longrightarrow Y$ be the blowup of a smooth surface $Y$ at $r$ general points,
$L$ an integral nef line bundle on $Y$ such that $\eta_r(L)\in \QQ\cap [0,2]$, and set
$\tilde{L}=L_{\eta_r(L)} = \pi^{*}L-\eta_r(L)E$.
If there is no irreducible curve $C$ on $X$ with multiplicities in $\{0,1\}$ such that $\tilde{L}\cdot C<0$, then
$\tilde{L}$ is a square-zero nef class on $X$. 

\bpf
With $[0,2]$ replaced by $[0,1]$, this result was previously known, and in that version one does not 
even need to check for possible $C$'s such that $\tilde{L}\cdot C<0$.  The version above, increasing $[0,1]$ to $[0,2]$, 
but requiring one to eliminate certain possible classes of $C$, appears as \cite[Corollary 2.7.2]{Dio}. \epf

\bpoint{Corollary ($v_{\alpha_9}$ and $v_{\beta_9}$ are nef)} 
\label{cor:v-alpha-and-v-beta-nef-for-9}
Let $X$ be the blowup of $\PP^1\times\PP^1$ at $9$ general points.  Then the classes
$v_{\alpha_9}$ and $v_{\beta_9}$ are nef on $X$.

\bpf
Setting $L=\Osh_{\PP^1\times\PP^1}(3,6)$, we have $\eta_9(L)=2$ and $\tilde{L}=\pi^{*}L-2E=3F_1+6F_2-2E=9 v_{\alpha_9}$. 
By Lemma \ref{lem:eta-is-two-means-nef} the only way that $\tilde{L}$ could not be nef is if there is an irreducible
curve $C$ with multiplicities in $\{0,1\}$ such that $\tilde{L}\cdot C<0$.  Thus, it suffices to show that there is
no such curve where $v_{\alpha_9}\cdot C<0$; such a curve must satisfy $C^2<0$. 

By Proposition \ref{prop:weak-lower-bound-K_X-pos} there is no such effective class with $K_X\cdot C\geq 0$, and
therefore we must have $K_X\cdot C<0$.   But, such curves (with all multiplicities in $\{0,1\}$ and $C^2<0$) are,
in the notation of \S\ref{sec:Outer-odd-r-intro},
either $C_2$, $C_3$, or the components which make up $C_1$ and $C_4$.

Since $C_3\cdot v_{\alpha_9}=0$, and $C_i\cdot v_{\alpha_9}>0$ for $i=1$, $2$, $4$, we conclude that $v_{\alpha_9}$ is nef.
Similarly, $v_{\beta_9}$ is nef.  \epf

\bpoint{Theorem} 
\label{thm:lower-bound-on-Seshadri-quotient} 
Let $\pi\colon X\longrightarrow Y$ be the blowup of a smooth surface $Y$ at $r$ general points.  
Let $C$ be the class of an effective curve on $X$ with all multiplicities equal to $m$, 
$m\geq 2$, and such that $C$ is not the symmetrization (as in \S\ref{sec:V_r-def}) of a curve $C'$ with
all multiplicities in $\{0,1\}$, and put $t=\frac{m-1}{m^2}$.  Then for all nef classes $L$ on $Y$, 

\begin{equation}
\label{eqn:lower-bound-Seshadri-quotient}
\frac{C\cdot \pi^{*}L}{C\cdot E} \geq \eta_r(L)\cdot \sqrt{\frac{r-t}{r}}.
\end{equation}

\bpf
This lower bound, in contrapositive form, is \cite[Theorem 4.1.1]{Dio}.  \epf

For use below we note that this result applies to real nef classes $L$, and not just integral ones.
(For instance, accepting that \eqref{eqn:lower-bound-Seshadri-quotient} holds for integral nef bundles on $Y$,
since both sides scale the same way, the inequality must also hold for rational nef classes, and then, by continuity,
for real nef classes.)

\bpoint{Lemma} 
\label{lem:alpha-estimate}
For $r\geq 9$ the estimates 
$$\tfrac{r-5}{2} \leq \alpha_r < \tfrac{r-4}{2}$$
hold.  Furthermore, the lower bound is strict whenever $r>9$.

\bpf
From the formula 
$$\alpha_r = \frac{(r-4)+\sqrt{(r-4)^2-16}}{4} $$
it is clear that $\alpha_r < \frac{r-4}{2}$. 
Set $u:=r-4$ to simplify the notation.  The lower bound to be proved is
$\frac{u}{2}-\frac{1}{2}\leq \frac{u+\sqrt{u^2-16}}{4}.$
Multiplying by $4$ and subtracting $u$ this becomes 
$(u-2)\leq \sqrt{u^2-16}.$
After squaring both sides the inequality becomes $4(u-5)\geq 0$, which clearly holds when $u\geq 5$, 
and is a strict inequality when $u>5$. \epf

\bpoint{Proposition} 
\label{prop:t-lower-estimate}
Suppose that $r\geq 11$ and that $t$, with $t\leq r$, is such that 
$$ \frac{2\alpha_r+(r-1)}{4(\alpha_r+1)}\geq \sqrt{\frac{r-t}{r}}.$$
Then $t>\frac{1}{4}$. 

\bpf
Squaring both sides, and using the identity $\alpha_r^2 = \left(\frac{r-4}{2}\right)\alpha_r-1$, in that form,
and in the form $(\alpha_r+1)^2=\frac{r}{2}\alpha_r$, the inequality above becomes 
$$\frac{6(r-2)\alpha_r+(r-1)^2-4}{8r\alpha_r}\geq \frac{r-t}{r}$$
or

\begin{equation}
\label{eqn:t-lower-estimate}
t \geq r - \frac{6(r-2)\alpha_r+(r-1)^2-4}{8\alpha_r}
= \frac{2(r+6)\alpha_r-(r-1)^2+4}{8\alpha_r}.
\end{equation}

Using Lemma \ref{lem:alpha-estimate} gives 

\begin{equation}
\label{eqn:t-lower-r-estimate}
t-\frac{1}{4} \geq 
\frac{2(r+6)\left(\frac{r-5}{2}\right)-(r-1)^2+4}{8\left(\frac{r-4}{2}\right)}-\frac{1}{4} 
=\frac{2r-23}{4r-16}.
\end{equation}

When $r\geq 12$ the right hand side of \eqref{eqn:t-lower-r-estimate} is clearly $>0$, proving the proposition
in those cases.  
For $r=11$ the right hand side of \eqref{eqn:t-lower-estimate} is $\approx 0.4836\ldots$, and so again $t>\frac{1}{4}$. \epf

\point 
\label{sec:proof-of-outer-odd-theorem}
{\em Proof of Theorem \ref{thm:Outer-nef-odd}.}  
When $r=9$ we have shown in Corollary  \ref{cor:v-alpha-and-v-beta-nef-for-9} that $v_{\alpha_9}$ and $v_{\beta_9}$
are nef, thus we may assume that $r\geq 11$. 

Suppose that there is a line bundle $L=\Osh_{Y}(e_1,e_2)$ with $\alpha_r\leq \frac{e_2}{e_1}$ such
that $\ep_r(L)$ is not equal to the value given in \eqref{eqn:Outer-nef-odd}.  Let $C$ be an equal multiplicity curve,
either irreducible, or the symmetrization of an irreducible curve as in \S\ref{sec:V_r-def}, 
which computes $\ep_r(L)$ (or even one which just imposes a stronger condition than imposed by $C_3$ or $C_4$). 

We have $C^2<0$ and, as explained in \S\ref{sec:Outline-of-outer-odd}, since $C$ is not equal to any of the 
$C_i$, $i=1$, \ldots, $4$, we also have $C\cdot K_X > 0$. 

Let $L'$ be the real nef class $\left(\frac{1}{\alpha_r+1},\frac{\alpha_r}{\alpha_r+1}\right)$ (i.e., the
first two coordinates of $v_{\alpha_r}$) on $\PP^1\times\PP^1$.   
Since $\alpha_r\leq \frac{e_2}{e_1}$, it follows that $C$ must also impose a stronger condition on the Seshadri constant 
of $L'$ than that imposed by $C_3$, i.e., that 

\begin{equation}
\label{eqn:C_3-lower-estimate}
\frac{C_3\cdot\pi^{*}L'}{C_3\cdot E} \geq \frac{C\cdot \pi^{*}L'}{C\cdot E}.
\end{equation}

\vspace{-0.5cm}
\hfill
\begin{tabular}{c}
\psset{yunit=3cm, xunit=5cm}
\begin{pspicture}(0.3,0.2)(1.0,0.9)
\psset{linecolor=gray}
\psclip{\pspolygon[linestyle=none](0.28,0.2)(1.2,0.2)(1.2,0.9)(0.28,0.9)}
\psline(0,0)(1,0)
\psline[linecolor=lightgray,linestyle=dashed,dash=1.5pt 1pt](!1 3 \ArcCoords 0.1 add)(!1 3 \ArcCoords pop 0.45)
\psline[arrows=->](!1 3 \ArcCoords pop 0)(!1 3 \ArcCoords pop 0.45)
\rput(!1 3 \ArcCoords pop 0.45 -0.07 0 \Trans){\tiny \color{gray} $L_{\gamma}$} 
\psset{linecolor=black}
\psset{linestyle=solid,dash = 3pt 2pt, linecolor=lightgray,linewidth=0.2mm}
\psline(!1 3.846954977 \ArcCoords pop 0)(!1 3.846954977 \ArcCoords pop 0.6)
\psset{linestyle=solid,linecolor=black}
\SecToArc{2}{3}{1}{10}{-0.2}{1.2} 
\LabelSecToArc{2}{3}{1}{10}{1.35}{\tiny $C_3^{\perp}$}
\SecToArc{6}{5}{1}{7}{-0.2}{1.2} 
\LabelSecToArc{6}{5}{1}{7}{-0.3}{\tiny $C^{\perp}$} 
\psset{linecolor=gray}
\psset{linecolor=black,fillstyle=none}
\DrawArc
\psset{fillstyle=solid,fillcolor=white,linecolor=black,linestyle=solid,linecolor=black}
\psset{fillstyle=solid,fillcolor=black}
\psset{fillstyle=solid,fillcolor=white}
\psset{fillcolor=lightgray}
\pscircle(!1 2.5 \ArcCoords){0.04} 
\rput(!1 2.5 \ArcCoords 0 0.1 \Trans){\tiny $v_{\alpha}$} 
\pscircle(!1 3.846954977 \ArcCoords pop 0.5134403289){0.03} 
\endpsclip
\rput(!1 3.846954977 \ArcCoords pop 0.15){\tiny \color{gray} $s$}
\end{pspicture} \\
\Fig \label{fig:troublesome-curve} \rule{0cm}{0.6cm} \\
\end{tabular} 

\vspace{-2.5cm}
\parshape 6 0cm 12cm 0cm 12cm 0cm 12cm 0cm 12cm 0cm 12cm 0cm \textwidth
The reason why is shown in Figure \ref{fig:troublesome-curve} at right.  Since $C$ must be below $C_3^{\perp}$ and
$C_4^{\perp}$ it follows
that $C^{\perp}$ must exit the square-zero cone farther to the left than the point where $C_3^{\perp}$ 
does, as shown in the picture.  
Thus $C^{\perp}$ and $C_3^{\perp}$ must intersect along a ray, of some 
slope $s$.  Since $C$ imposes a stronger condition on the Seshadri constant of $L$ than $C_3$, this means that 
$\frac{e_2}{e_1}\leq s$.  But then $\alpha_r\leq s$ too, and so similarly $C$ imposes a stronger condition on the 
Seshadri constant of $L'$ than that imposed by $C_3$, giving \eqref{eqn:C_3-lower-estimate}.

By Proposition \ref{prop:weak-lower-bound-K_X-pos}, $C$ has multiplicities $m\geq 2$, and is not the symmetrization
of a curve $C'$ with all multiplicities in $\{0,1\}$.  Thus we may apply Theorem 
\ref{thm:lower-bound-on-Seshadri-quotient} to $L'$, to get 
$$\frac{C\cdot\pi^{*}L'}{C\cdot E} \geq \eta_{r}(L')\cdot \sqrt{\frac{r-t}{r}},$$
with $t=\frac{m-1}{m^2}$.  Since $\eta_{r}(L') = \frac{2}{r}$, combining the previous inequality with 
\eqref{eqn:C_3-lower-estimate} gives
$$
\frac{\frac{r-1}{2}\left(\frac{1}{\alpha_r+1}\right)+1\left(\frac{\alpha_r}{\alpha_r+1}\right)}{r}
= \frac{C_3\cdot\pi^{*}L'}{C_3\cdot E} \geq 
\eta_r(L')\cdot \sqrt{\frac{r-t}{r}}
=
\frac{2}{r}\cdot \sqrt{\frac{r-t}{r}},$$
or 
$$ \frac{2\alpha_r+(r-1)}{4(\alpha_r+1)}\geq \sqrt{\frac{r-t}{r}}.$$

Applying Proposition \ref{prop:t-lower-estimate} we conclude that $t>\frac{1}{4}$.  
But $t=\frac{m-1}{m^2}$, and for $m\geq 1$ the maximum value of $\frac{m-1}{m^2}$ is $\frac{1}{4}$, occurring when $m=2$.
This contradiction shows that there can be no such curve $C$, concluding the proof of Theorem \ref{thm:Outer-nef-odd}.
\epf

\bpoint{Remarks} 
\RemNum Here is an outline of the argument above. The final step is that the inequality 

\begin{equation}
\label{eqn:contradiction-inequality}
\frac{C_3\cdot\pi^{*}L'}{C_3\cdot E} \geq \eta_{r}(L')\cdot \sqrt{\frac{r-t}{r}}
\end{equation}
leads to a contradiction whenever $t=\frac{m-1}{m^2}$ with $m\geq 1$.

Note that we cannot get \eqref{eqn:contradiction-inequality} by applying 
Theorem \ref{thm:lower-bound-on-Seshadri-quotient} to $C_3$, since $C_3$ does not satisfy the hypothesis of the 
theorem --- all the multiplicities of $C_3$ are equal to $1$. 
(And it is good that we cannot --- the curve $C_3$ exists, and the inequality leads to a contradiction!) 

But, if we assume the existence of a curve $C$ which imposes a stronger condition than $C_3$
on a line bundle of slope $\geq \alpha_r$, we get the inequality \eqref{eqn:C_3-lower-estimate}.  
Applying Theorem \ref{thm:lower-bound-on-Seshadri-quotient} to $C$, and combining the inequality which results with
\eqref{eqn:C_3-lower-estimate} we arrive at \eqref{eqn:contradiction-inequality}, and thus a contradiction.  
Therefore no such curve $C$ can exist.

\RemNum 
As is clear from \eqref{eqn:t-lower-r-estimate}, as $r$ gets large, the lower bound estimate on $t$ goes 
to $\frac{3}{4}$, larger than the $\frac{1}{4}$ needed to give a contradiction.  This suggests that one can 
improve the region on which the formulas in \eqref{eqn:Outer-nef-odd} hold.

For $s>0$, let $L(s)$ denote the real nef class of type $(1,s)$ on $\PP^1\times\PP^1$. We note that
this is different than the $L_{\gamma}$ used throughout the rest of the paper.  As the idea for the proof of Theorem
\ref{thm:Outer-nef-odd} shows, whenever (for a fixed $r$, $r$ odd, $r\geq 9$) $s$ is such that 

\begin{equation}
\label{eqn:s-inequality} 
\frac{C_3\cdot \pi^{*}L(s)}{r\cdot \eta_r(L(s))} < \sqrt{\frac{r-\frac{1}{4}}{r}},
\end{equation}
one can conclude that the formulas in \eqref{eqn:Outer-nef-odd} hold for all line bundles $L$ with slope outside
$(\frac{1}{s},s)$.   Solving the inequality \eqref{eqn:s-inequality}, one finds that the smallest
$s$ which works is $s(r)=\frac{1}{4}(2r-\sqrt{12r+1}+1)$.   One has
$$\frac{(\sqrt{r}-1)^2}{2} < s(r) < \alpha_r.$$
(The leftmost term is the slope where $C_3^{\perp}$ exits the square-zero cone, see \eqref{eqn:C_3-boundaries}.)
In particular, $(\frac{1}{s(r)},s(r))\subseteq (\beta_r,\alpha_r)$.

Thus the previous argument can be used to prove the following result, which, since
$(\frac{1}{s(r)},s(r))\subseteq (\beta_r,\alpha_r)$, is stronger than Theorem \ref{thm:Outer-nef-odd}, and gives
an exact value of the Seshadri constant for a narrow range of inner bundles.

\bpoint{Theorem (extension of Theorem \ref{thm:Outer-nef-odd})} 
\label{thm:Outer-nef-odd-extended}
Suppose that $r$ is odd, $r\geq 9$, that $L=\Osh_{\PP^1\times\PP^1}(e_1,e_2)$ with $e_1,e_2\geq 1$, and
that $\frac{e_2}{e_1}\not\in (\frac{1}{s(r)},s(r))$, where $s(r) = \frac{1}{4}(2r-\sqrt{12r+1}+1)$.  Then

\vspace{-0.25cm}
\begin{equation}
\label{eqn:Outer-nef-odd}
\ep_r(L) = \left\{ 
\begin{array}{cl}
e_2 & \mbox{if $\frac{e_2}{e_1}\leq \frac{2}{r+1}$}, \\
\frac{2e_1+(r-1)e_2}{2r} & \mbox{if $\frac{e_2}{e_1} \in [\frac{2}{r+1},\frac{1}{s(r)}]$}, \rule{0cm}{0.6cm}\\
\frac{(r-1)e_1+2e_2}{2r} & \mbox{if $\frac{e_2}{e_1} \in [s(r),\frac{r+1}{2}]$}, \rule{0cm}{0.6cm}\\
e_1 & \mbox{if $\frac{r+1}{2}\leq \frac{e_2}{e_1}$}. \rule{0cm}{0.6cm} \\
\end{array}
\right.
\end{equation}

\point
If one can get better lower bound on the multiplicity of a putative curve $C$ which is $K_{X}$-positive and
satisfies $C^2<0$ (and thus computes the Seshadri constant for some outer bundle), 
e.g., $m\geq 3$, $m\geq 4$, etc, then
one can replace $r-\frac{1}{4}$ in \eqref{eqn:s-inequality} by $r-\frac{2}{9}$, $r-\frac{3}{16}$, etc, and get 
solutions for $s(r)$ which are closer to $\frac{(\sqrt{r}-1)^2}{2}$, and further improve the result above.

\point 
\label{sec:no-automorphisms-when-r-is-odd}
Let us return to the question raised in (\Euler{2}) of \S\ref{sec:remarks-on-Tr}, namely showing that, other than
switching the factors, there are no automorphisms of the problem fixing $V_r$ when $r$ is odd. 

As we see in Figure \ref{fig:General-outer-picture-odd-r}({\em a}), when $r$ is odd the only nef, outer, square-zero
symmetric classes are the fibre classes $F_1$ and $F_2$ (this also holds when $r\leq 7$, see 
\S\ref{sec:r=1-3-5}--\S\ref{sec:r=7} below).
Thus, any automorphism of the problem either fixes $F_1$ and $F_2$, or swaps them.  
But, if $F_1$ and $F_2$ are fixed, in order to preserve the intersection form on $V_r$, and preserve the nef cone,
the automorphism must also fix $E$, and thus be the identity on $V_r$.  Thus, when $r$ is odd, swapping the factors is the 
only nontrivial automorphism of the problem. 

\section{Small $r$ ($1\leq r\leq 8$)}
\label{sec:Small-r} 

In this section we list the results for $r$ between $1$ and $8$.  For $r\leq 7$ the blowup of $\PP^1\times \PP^1$
at $r$ general points is Fano, and so by Mori's theorem the nef and effective cones are polyhedral.  When $r=8$
the anticanonical bundle is nef and effective, and the both the nef and effective cones are polyhedral away
from the class of $-K_{X}$. 

We also note that for $r\leq 7$ the numbers $\alpha_r$ and $\beta_r$ are complex numbers, and so the 
vectors $v_{\alpha_r}$ and $v_{\beta_r}$ are not in the real vector space $V_r$.

\bpoint{$\mathbf{r=2}$, $\mathbf{4}$, $\mathbf{6}$} 
\label{sec:r=2-4-6}

Theorem \ref{thm:Tr-transformation} is valid for all even $r$. In contrast to the cases when $r\geq 8$, where $T_r$ has
infinite order, $T_2$, $T_4$, and $T_6$ have orders $3$, $4$, and $6$ respectively, and
the $\alpha_r$ are the complex roots of unity $\alpha_2=e^{\frac{2\pi i}{3}}$, $\alpha_4=i$, 
and $\alpha_6=e^{\frac{2\pi i}{6}}$.

Theorem \ref{thm:Outer-nef} showed, when $r\geq 10$, $r$ even, that the intersection of the nef cone in $V_{r}$ and 
the half plane $K_{X}^{\leq 0}$ is generated by the classes $v_{\alpha_r}$, $v_{\beta_r}$ and $T_r^{n}(F_2)$ 
with $n\in \ZZ$, and that the intersection
of the effective cone in $V_r$ with the half plane $K_{X}^{\leq 0}$ is generated by $v_{\alpha_r}$, $v_{\beta_r}$ and
the $T_r^{n}(E)$, with $n\in \ZZ$.
The arguments in that theorem work here, with the only change being that, since $X$ is Fano, the intersection with 
$K_{X}^{< 0}$ is all of the nef and effective cones respectively, and $v_{\alpha_r}$ and $v_{\beta_r}$ do not
appear at all.

Thus, in the cases $r\in \{2,4,6\}$, the nef cone is spanned by the classes $T_{r}^n(F_2)$, for $n=0$, $\ldots$, 
$\ord(T_r)-1$ (i.e., to $2$, $3$, and $5$ respectively).   Similarly, the effective cones in these cases are spanned 
by $T_r^{n}(E)$ for $n=0$, \ldots, $\ord(T_r)-1$. 

Figure \ref{fig:Cones-for-2-4-6} below shows the nef cones and orbits of these classes.  The classes
represented by a hollow circle are the orbits of $F_2$, and the classes represented by a solid circle are the
orbits of $E$ (with the exception of the class $E$ itself, which is not shown).   

\vspace{0.25cm}

\hspace{-1cm}
\begin{tabular}{ccccc}
\begin{tabular}{c}
\psset{yunit=3cm, xunit=5cm}
\begin{pspicture}(0,0)(1,1)
\psset{linecolor=gray}
\psset{linecolor=black}
\psset{linestyle=dashed,dash = 3pt 2pt, linecolor=lightgray}
\TanToArc{1}{0}{-0.1}{0.8}
\TanToArc{0}{1}{-0.8}{0.1}
\TanToArc{1}{1}{-0.75}{0.75}
\psset{linestyle=solid}
\psset{fillcolor=vlgray,linecolor=gray,fillstyle=solid}
\pspolygon(!1 0 0 \SliceCoords)(!1 1 1 2 \RSliceCoords)(!0 1 0 \SliceCoords)
\psset{linecolor=black,fillstyle=none}
\DrawArc
\psset{fillstyle=solid,fillcolor=white,linecolor=black,linestyle=solid,linecolor=black}
\psset{fillstyle=solid,fillcolor=black}
\pscircle(!2 0 1 2 \RSliceCoords){0.04}
\rput(!2 0 1 2 \RSliceCoords 0 0.1 \Trans){\tiny $(2,0,-1)$}
\pscircle(!0 2 1 2 \RSliceCoords){0.04}
\rput(!0 2 1 2 \RSliceCoords 0 0.1 \Trans){\tiny $(0,2,-1)$}
\psset{fillstyle=solid,fillcolor=white}
\pscircle(!1 0 \ArcCoords){0.04}
\pscircle(!0 1 \ArcCoords){0.04}
\pscircle(!1 1 1 2 \RSliceCoords){0.04}
\SLabel{1}{\tiny $1$}
\end{pspicture} 
\end{tabular}
& \rule{0.5cm}{0cm} & 
\begin{tabular}{c}
\psset{yunit=3cm, xunit=5cm}
\begin{pspicture}(0,0)(1,1)
\psset{linecolor=gray}
\psset{linecolor=black}
\psset{linestyle=dashed,dash = 3pt 2pt, linecolor=lightgray}
\TanToArc{1}{0}{-0.1}{0.55}
\TanToArc{0}{1}{-0.55}{0.1}
\TanToArc{2}{1}{-0.55}{0.3}
\TanToArc{1}{2}{-0.3}{0.55}
\psset{linestyle=solid}
\psset{fillcolor=vlgray,linecolor=gray,fillstyle=solid}
\pspolygon(!1 0 0 \SliceCoords)(!2 1 \ArcCoords)(!1 2 \ArcCoords)(!0 1 0 \SliceCoords)
\psset{linecolor=black,fillstyle=none}
\DrawArc
\psset{fillstyle=solid,fillcolor=white,linecolor=black,linestyle=solid,linecolor=black}
\psset{fillstyle=solid,fillcolor=black}
\pscircle(!4 0 1 4 \RSliceCoords){0.04}
\rput(!4 0 1 4 \RSliceCoords 0 0.1 \Trans){\tiny $(4,0,-1)$}
\pscircle(!0 4 1 4 \RSliceCoords){0.04}
\rput(!0 4 1 4 \RSliceCoords 0 0.1 \Trans){\tiny $(0,4,-1)$}
\pscircle(!4 4 3 4 \RSliceCoords){0.04}
\rput(!4 4 3 4 \RSliceCoords 0 0.1 \Trans){\tiny $(4,4,-3)$} 
\psset{fillstyle=solid,fillcolor=white}
\pscircle(!1 0 \ArcCoords){0.04}
\pscircle(!0 1 \ArcCoords){0.04}
\pscircle(!2 1 \ArcCoords){0.04}
\pscircle(!1 2 \ArcCoords){0.04}
\SLabel{0.5}{\tiny $\frac{1}{2}$}
\SLabel{2}{\tiny $2$}
\end{pspicture} 
\end{tabular}
& \rule{0.5cm}{0cm} & 
\begin{tabular}{c}
\psset{yunit=3cm, xunit=5cm}
\begin{pspicture}(0,0)(1,1)
\psset{linecolor=gray}
\psset{linecolor=black}
\psset{linestyle=dashed,dash = 3pt 2pt, linecolor=lightgray}
\TanToArc{1}{0}{-0.1}{0.5}
\TanToArc{0}{1}{-0.5}{0.1}
\TanToArc{3}{1}{-0.5}{0.2}
\TanToArc{1}{3}{-0.2}{0.5}
\TanToArc{4}{3}{-0.3}{0.2}
\TanToArc{3}{4}{-0.2}{0.3}
\psset{linestyle=solid}
\psset{fillcolor=vlgray,linecolor=gray,fillstyle=solid}
\pspolygon(!1 0 0 \SliceCoords)(!3 1 \ArcCoords)(!4 3 \ArcCoords)(!3 4 \ArcCoords)(!1 3 \ArcCoords)(!0 1 0 \SliceCoords)
\psset{linecolor=black,fillstyle=none}
\DrawArc
\psset{fillstyle=solid,fillcolor=white,linecolor=black,linestyle=solid,linecolor=black}
\psset{fillstyle=solid,fillcolor=black}
\pscircle(!6 0 1 6 \RSliceCoords){0.04}
\rput(!6 0 1 6 \RSliceCoords 0 0.1 \Trans){\tiny $(6,0,-1)$}
\pscircle(!0 6 1 6 \RSliceCoords){0.04}
\rput(!0 6 1 6 \RSliceCoords 0 0.1 \Trans){\tiny $(0,6,-1)$}
\pscircle(!6 12 5 6 \RSliceCoords){0.04}
\rput(!6 12 5 6 \RSliceCoords 0.15 0.1 \Trans){\tiny $(6,12,-5)$} 
\pscircle(!12 6 5 6 \RSliceCoords){0.04}
\rput(!12 6 5 6 \RSliceCoords -0.15 0.1 \Trans){\tiny $(12,6,-5)$} 
\pscircle(!12 12 7 6 \RSliceCoords){0.04}
\rput(!12 12 7 6 \RSliceCoords 0 0.15 \Trans){\tiny $(12,12,-7)$} 
\psset{fillstyle=solid,fillcolor=white}
\pscircle(!1 0 \ArcCoords){0.04}
\pscircle(!0 1 \ArcCoords){0.04}
\pscircle(!3 1 \ArcCoords){0.04}
\pscircle(!1 3 \ArcCoords){0.04}
\pscircle(!3 4 \ArcCoords){0.04}
\pscircle(!4 3 \ArcCoords){0.04}
\SLabel{0.3333333}{\tiny $\frac{1}{3}$}
\SLabel{0.75}{\tiny $\frac{3}{4}$}
\SLabel{1.3333333}{\tiny $\frac{4}{3}$}
\SLabel{3}{\tiny $3$}
\end{pspicture} 
\end{tabular} \\
$r=2$ & & $r=4$ & & $r=6$ \rule{0cm}{1.0cm}\\
\multicolumn{5}{c}{\Fig\label{fig:Cones-for-2-4-6}} \\
\end{tabular}

\vspace{0.25cm} 
The Seshadri constants for $L=\Osh_{Y}(e_1,e_2)$ in these cases are~:

\hspace{-1.5cm}
\begin{tabular}{ccc}
\begin{tabular}{c}
$\ep_{2}(L) = \left\{
\begin{array}{cl}
e_2 & \mbox{if $\frac{e_2}{e_1}\in (0,1]$} \\[2mm]
e_1 & \mbox{if $\frac{e_2}{e_1}\in [1,\infty)$} \\
\end{array}
\right.
$
\end{tabular};
&
\begin{tabular}{c}
$\ep_{4}(L) = \left\{
\begin{array}{cl}
e_2 & \mbox{if $\frac{e_2}{e_1}\in (0,\frac{1}{2}]$} \\[2mm]
\frac{e_1+e_2}{3}  & \mbox{if $\frac{e_2}{e_1}\in [\frac{1}{2},2]$} \\[2mm]
e_1 & \mbox{if $\frac{e_2}{e_1}\in [2,\infty)$} \\
\end{array}
\right.
$
\end{tabular};
&
\begin{tabular}{c}
$\ep_{6}(L) = \left\{
\begin{array}{cl}
e_2 & \mbox{if $\frac{e_2}{e_1}\in (0,\frac{1}{3}]$} \\[2mm]
\frac{e_1+2e_2}{5} & \mbox{if $\frac{e_2}{e_1}\in [\frac{1}{3},\frac{3}{4}]$} \\[2mm]
\frac{2e_1+2e_2}{7} & \mbox{if $\frac{e_2}{e_1}\in [\frac{3}{4},\frac{4}{3}]$} \\[2mm]
\frac{2e_1+e_2}{5} & \mbox{if $\frac{e_2}{e_1}\in [\frac{4}{3},3]$} \\[2mm]
e_1 & \mbox{if $\frac{e_2}{e_1}\in [3,\infty)$} \\
\end{array}
\right.
$
\end{tabular}.
\end{tabular}

\newpage
\bpoint{$\mathbf{r=8}$} 
\label{sec:r=8}

When $r=8$, $\alpha_{8}=\beta_{8}=1$, and $v_{\alpha_{8}}=v_{\beta_{8}}=-\frac{1}{4}K_{X}$.
The transformation $T_{8}$ is unipotent with a single Jordan block.

As in the previous even cases, the argument of Theorem \ref{thm:Outer-nef} shows that the classes $T_{8}^{n}(F_2)$, 
$n\in \ZZ$, and $v_{\alpha_8}$ ($=v_{\beta_8}$) generate the intersection of the nef cone with $K_{X}^{\leq 0}$.
The plane $K_{X}^{\perp}$ meets the square zero cone (and the nef cone, and the effective cone) only along the ray
spanned by $v_{\alpha_8}$.

\vspace{-0.5cm}
\hfill
\begin{tabular}{c}
\begin{tabular}{c}
\psset{yunit=3cm, xunit=5cm}
\begin{pspicture}(0,-0.2)(1,1)
\psset{linecolor=gray}
\psline(! 8 \AlphaVal 1 \ArcCoords exch pop 0 exch)(! 8 \AlphaVal 1 \ArcCoords exch pop 1 exch)
\rput(!1 8 \AlphaVal \ArcCoords 0.6 0 \Trans){\tiny $K_X^{\perp}$}
\psset{linecolor=black}
\psset{linestyle=dashed,dash = 3pt 2pt, linecolor=lightgray}
\TanToArc{1}{0}{-0.1}{0.45}
\TanToArc{4}{1}{-0.4}{0.2}
\TanToArc{1}{4}{-0.2}{0.4}
\TanToArc{0}{1}{-0.45}{0.1}
\psset{linestyle=solid}
\psset{fillcolor=vlgray,linecolor=gray,fillstyle=solid}
\pscustom{
\psline(!0 \EightNefPoint)(!-1 \EightNefPoint)(!-2 \EightNefPoint)(!-3 \EightNefPoint)(!-4 \EightNefPoint)(!-5 \EightNefPoint)
\DrawPartArc{0.49836}{0.590164}
\psline(!6 \EightNefPoint)(!5 \EightNefPoint)(!4 \EightNefPoint)(!3 \EightNefPoint)(!2 \EightNefPoint)(!1 \EightNefPoint)(!0 \EightNefPoint) 
}  
\psset{linecolor=black,fillstyle=none}
\DrawArc
\psset{fillstyle=solid,fillcolor=white,linecolor=black,linestyle=solid,linecolor=black}
\psset{fillstyle=solid,fillcolor=black}
\pscircle(!0 8 1 8 \RSliceCoords){0.04}
\rput(!0 8 1 8 \RSliceCoords 0.10 0 \Trans){\tiny $C_1$} 
\pscircle(!8 24 7 8 \RSliceCoords){0.04}
\pscircle(!24 48 17 8 \RSliceCoords){0.04}
\pscircle(!48 24 17 8 \RSliceCoords){0.04}
\pscircle(!24 8 7 8 \RSliceCoords){0.04}
\pscircle(!8 0 1 8 \RSliceCoords){0.04}
\rput(!8 0 1 8 \RSliceCoords -0.10 0 \Trans){\tiny $C_{-1}$} 
\psset{fillstyle=solid,fillcolor=white}
\multido{\n=0+1}{30}{%
\PutEightCircle{\n}
}
\multido{\n=-1+-1}{30}{%
\PutEightCircle{\n}
}
\pscircle(!1 1 \ArcCoords){0.04}
\rput(!1 1 \ArcCoords 0 0.1 \Trans){\tiny $v_{\alpha}=v_{\beta}$} 
\rput(!1 \EightNefPoint 0.05 -0.05 \Trans){\tiny $\xi_0$} 
\rput(!0 \EightNefPoint -0.05 -0.05 \Trans){\tiny $\xi_{-1}$} 
\rput(!2 \EightNefPoint 0.07 0.01 \Trans){\tiny $\xi_{1}$} 
\rput(!-1 \EightNefPoint -0.08 0.01 \Trans){\tiny $\xi_{-2}$} 
\SLabel{0.25}{\tiny$\frac{1}{4}$}
\SLabel{0.4444444}{\tiny$\frac{4}{9}$}
\rput(0.45,-0.12){\tiny$\cdots$}
\SLabel{0.625}{\tiny$\frac{9}{16}$}
\rput(0.57,-0.12){\tiny$\cdots$}
\SLabel{1}{\tiny $1$} 
\SLabel{1.7777777}{\tiny$\frac{16}{9}$}
\SLabel{2.25}{\tiny$\frac{9}{4}$}
\SLabel{4}{\tiny$\frac{4}{1}$}
\end{pspicture} 
\end{tabular}\\ 
\Fig\label{fig:cone-when-r=8} \\
\end{tabular}

\vspace{-3.75cm}

\parshape 1 0cm 10cm
As in \S\ref{thm:Outer-nef}, for each $n\in \ZZ$ we  set $\xi_n:=T^{n}_{8}(F_2)$ and $C_n:=T_{8}^{n}(E)$.
Because $T_{8}$ is unipotent, the coordinates of $\xi_n$ and $C_n$ are quadratic functions of $n$. 
Specifically, 

\vspace{0.25cm}
\hspace{1cm}
$\xi_n=\left(n^2,(n+1)^2,-\binom{n+1}{2}\right)$,

and

\hspace{1cm}
$C_n=\left(4n(n-1),4n(n+1),1-2n^2\right)$.

\vspace{0.25cm}
As in the proof of Theorem \ref{thm:Outer-nef} we have  $\xi_{n-1}\cdot C_n=0=\xi_{n}\cdot C_n$ for all $n\in \ZZ$.
Thus the $C_{n}$, along with the limiting class $v_{\alpha_8}$, are dual to the nef cone, and so generate the effective 
cone of $X$.  The picture in this case is shown in Figure \ref{fig:cone-when-r=8} above.

Given an ample line bundle $L=\Osh_{Y}(e_1,e_2)$ on $Y$, then as long as $\frac{e_2}{e_1}\neq 1$, 
the Seshadri constant of $L$ is computed by one of the curves
$C_n$ above, by the argument in the proof of Corollary \ref{cor:Outer-nef}.  
To find $C_{n}$ we look for the value of $n$ 
so that the slope of $L$ is between the slope of $\xi_{n}$ and the slope of $\xi_{n-1}$.  
That is, a value of $n$ so that 

\begin{equation}
\label{eqn:r=8-Cn-bounds}
 \frac{(n+1)^2}{n^2} \leq \frac{e_2}{e_1} \leq \frac{n^2}{(n-1)^2}.
\end{equation}

Here $\frac{1}{0}$ is interpreted as $\infty$ if necessary, and the conditions imply that $n\geq 1$ if 
$\frac{e_2}{e_1}>1$ while $n\leq -1$ if $\frac{e_2}{e_1}< 1$. 
For example, $C_1$ computes the Seshadri constant for those $L$ whose slopes are in $[4,\infty)$, and $C_{-1}$
computes the Seshadri constant for those $L$ whose slopes are in $(0,\frac{1}{4}]$. 

Using \eqref{eqn:r=8-Cn-bounds}, one possible formula for $n$ in these cases is 
$n = \left\lceil\frac{1}{\sqrt{\frac{e_2}{e_1}}-1} \right\rceil$.   

The reason for excluding $\frac{e_2}{e_1}=1$ is that this is the slope of $v_{\alpha_8}$ (i.e., the limit
of the $\xi_n$, up to scaling, as $n\to+\infty$), and also the slope of $v_{\beta_{8}}$ (i.e., the limit of the
$\xi_{n}$ as $n\to -\infty$), and so there are no $\xi_{n}$ and $\xi_{n-1}$ whose slopes bracket $1$.  But, 
$v_{\alpha_8}$ has slope $1$, is nef, and on the square zero cone.  So, for line bundles $L$ of slope $1$ the
Seshadri constant is the maximum possible value $\ep_{8}(L) = \eta_{8}(L)= \frac{e_1}{2}$. 

In summary, 

$$\ep_{8}(L) = \left\{
\begin{array}{cl}
\frac{e_1}{2} & \mbox{if $e_1=e_2$} \\[4mm]
\frac{n(n+1)(e_1+e_2)}{2(2n^2-1)} & \mbox{if $e_1\neq e_2$, with $n=\left\lceil\frac{1}{\sqrt{\frac{e_2}{e_1}}-1}\right\rceil$} \\ 
\end{array}\right..$$

Note that the method of finding $n$ by using $\phi_r$ as defined by \eqref{eqn:eta-r-def} does not work when $r=8$.
Since $\alpha_8=1$, the denominator of \eqref{eqn:eta-r-def} is zero.

\bpoint{$\mathbf{r=1}$, $\mathbf{3}$, $\mathbf{5}$} 
\label{sec:r=1-3-5}

Recall that in \S\ref{sec:Outer-Odd} we defined the curve classes 

\begin{equation}
\label{eqn:four-curve-classes}
\rule{1cm}{0cm}C_1:=(r,0,-1),\,\,C_2:=\left(\tfrac{r-1}{2},1,-1\right),\,\,
C_3:=\left(1,\tfrac{r-1}{2},1,-1\right),\,\,\mbox{and}\,\,C_4:=(0,r,-1),
\end{equation}

and that for odd $r\geq 9$ these classes determine the nef cone for all outer line bundles.  When $r<8$ all ample
line bundles on $\PP^1\times\PP^1$ are outer (since the plane $K_{X}^{\perp}$ does not intersect the square-zero cone).
For $r\in \{1,3,5\}$ the argument in \S\ref{sec:Outline-of-outer-odd} shows that no other symmetric curve class affects
the Seshadri constant of an outer bundle (i.e., any ample bundle in this case).  Thus, for $r\in \{1,3,5\}$ the curve
classes above determine the entire nef cone.  

One other difference in these cases is that when $r\in \{1,3\}$ some of these curve classes coincide. 
Specifically, when $r=1$ we have $C_3=C_1$ and $C_2=C_4$, and when $r=3$ we have $C_2=C_3$. 
The pictures of these curves, and the corresponding nef cones cut out by the half planes $C_i^{\leq 0}$, $i=1$,\ldots
$4$ is shown below.  With the exception of the fibre classes, the small white circles on the square-zero cone are not
nef, but are rather classes whose tangent lines contain one of the $C_i$.   As explained in 
\S\ref{sec:graphical-arguments}(a) these classes determine the planes $C_i^{\perp}$. 

\vspace{0.75cm} 

\hspace{-1.0cm}
\begin{tabular}{ccccc}
\begin{tabular}{c}
\psset{yunit=3cm, xunit=5cm}
\begin{pspicture}(0,0)(1,1)
\psset{linecolor=gray}
\psset{linecolor=black}
\psset{linestyle=dashed,dash = 3pt 2pt, linecolor=lightgray}
\TanToArc{1}{0}{-0.1}{1.1}
\TanToArc{0}{1}{-1.1}{0.1}
\TanToArc{1}{0.5}{-0.25}{1.1}
\TanToArc{0.5}{1}{-1.1}{0.25}
\psset{linestyle=solid}
\SecToArc{1}{0}{0.5}{1}{-0.1}{1.1}
\SecToArc{0}{1}{1}{0.5}{-0.1}{1.1}
\psset{fillcolor=vlgray,linecolor=gray,fillstyle=solid}
\pspolygon(!1 0 \ArcCoords)(!1 1 1 \SliceCoords)(!0 1 \ArcCoords)
\psset{linecolor=black,fillstyle=none}
\DrawArc
\psset{fillstyle=solid,fillcolor=white,linecolor=black,linestyle=solid,linecolor=black}
\psset{fillstyle=solid,fillcolor=black}
\pscircle(!1 0 1 \SliceCoords){0.04}
\rput(!1 0 1 \SliceCoords 0 0.1 \Trans){\tiny $(1,0,-1)$}
\pscircle(!0 1 1 \SliceCoords){0.04}
\rput(!0 1 1 \SliceCoords 0 0.1 \Trans){\tiny $(0,1,-1)$}
\psset{fillstyle=solid,fillcolor=white}
\pscircle(!1 0 \ArcCoords){0.04}
\pscircle(!0 1 \ArcCoords){0.04}
\pscircle(!0.5 1 \ArcCoords){0.04}
\pscircle(!1 0.5 \ArcCoords){0.04}
\psset{fillcolor=lightgray}
\pscircle(!1 1 1 \SliceCoords){0.04}
\SLabel{1}{\tiny $1$}
\end{pspicture} 
\end{tabular}
& \rule{0.5cm}{0cm} & 
\begin{tabular}{c}
\psset{yunit=3cm, xunit=5cm}
\begin{pspicture}(0,0)(1,1)
\psset{linecolor=gray}
\psset{linecolor=black}
\psset{linestyle=dashed,dash = 3pt 2pt, linecolor=lightgray}
\TanToArc{1}{0}{-0.1}{0.7}
\TanToArc{0}{1}{-0.7}{0.1}
\TanToArc{3}{2}{-0.65}{0.2}
\TanToArc{2}{3}{-0.2}{0.65}
\TanToArc{2 3 sqrt add}{1}{-0.2}{0.55}
\TanToArc{1}{2 3 sqrt add}{-0.55}{0.2}
\psset{linestyle=solid}
\SecToArc{1}{0}{3}{2}{-0.1}{1.1}
\SecToArc{0}{1}{2}{3}{-0.1}{1.1}
\SecToArc{2 3 sqrt add}{1}{1}{2 3 sqrt add}{-0.1}{1.1}
\psset{fillcolor=vlgray,linecolor=gray,fillstyle=solid}
\pspolygon(!1 0 \ArcCoords)(!2 1 1 3 \RSliceCoords)(!1 2 1 3 \RSliceCoords)(!0 1 \ArcCoords)
\psset{linecolor=black,fillstyle=none}
\DrawArc
\psset{fillstyle=solid,fillcolor=white,linecolor=black,linestyle=solid,linecolor=black}
\psset{fillstyle=solid,fillcolor=black}
\pscircle(!3 0 1 3 \RSliceCoords){0.04}
\rput(!3 0 1 3 \RSliceCoords 0 0.1 \Trans){\tiny $(3,0,-1)$}
\pscircle(!0 3 1 3 \RSliceCoords){0.04}
\rput(!0 3 1 3 \RSliceCoords 0 0.1 \Trans){\tiny $(0,3,-1)$}
\pscircle(!1 1 1 3 \RSliceCoords){0.04}
\rput(!1 1 1 3 \RSliceCoords 0 0.1 \Trans){\tiny $(1,1,-1)$}
\psset{fillstyle=solid,fillcolor=white}
\pscircle(!1 0 \ArcCoords){0.04}
\pscircle(!0 1 \ArcCoords){0.04}
\pscircle(!3 2 \ArcCoords){0.04}
\pscircle(!2 3 \ArcCoords){0.04}
\pscircle(!2 3 sqrt add 1 \ArcCoords){0.04}
\pscircle(!1 2 3 sqrt add \ArcCoords){0.04}
\psset{fillcolor=lightgray}
\pscircle(!2 1 1 3 \RSliceCoords){0.04}
\pscircle(!1 2 1 3 \RSliceCoords){0.04}
\SLabel{0.5}{\tiny $\frac{1}{2}$}
\SLabel{2}{\tiny $2$}
\end{pspicture} 
\end{tabular}
& \rule{0.5cm}{0cm} & 
\begin{tabular}{c}
\psset{yunit=3cm, xunit=5cm}
\begin{pspicture}(0,0)(1,1)
\psset{linecolor=gray}
\psset{linecolor=black}
\psset{linestyle=dashed,dash = 3pt 2pt, linecolor=lightgray}
\TanToArc{1}{0}{-0.1}{0.5}
\TanToArc{0}{1}{-0.5}{0.1}
\TanToArc{5}{2}{-0.50}{0.2}
\TanToArc{2}{5}{-0.2}{0.50}
\TanToArc{3 5 sqrt add}{4}{-0.1}{0.4}
\TanToArc{3 5 sqrt sub}{4}{-0.4}{0.1}
\TanToArc{4}{3 5 sqrt add}{-0.4}{0.1}
\TanToArc{4}{3 5 sqrt sub}{-0.1}{0.4}
\psset{linestyle=solid}
\SecToArc{1}{0}{5}{2}{-0.1}{1.1}
\SecToArc{0}{1}{2}{5}{-0.1}{1.1}
\SecToArc{3 5 sqrt add}{4}{3 5 sqrt sub}{4}{-0.1}{1.1}
\SecToArc{4}{3 5 sqrt add}{4}{3 5 sqrt sub}{-0.1}{1.1}
\psset{fillcolor=vlgray,linecolor=gray,fillstyle=solid}
\pspolygon(!1 0 \ArcCoords)(!3 1 1 5 \RSliceCoords)(!5 5 3 5 \RSliceCoords)(!1 3 1 5 \RSliceCoords)(!0 1 \ArcCoords)
\psset{linecolor=black,fillstyle=none}
\DrawArc
\psset{fillstyle=solid,fillcolor=white,linecolor=black,linestyle=solid,linecolor=black}
\psset{fillstyle=solid,fillcolor=black}
\pscircle(!5 0 1 5 \RSliceCoords){0.04}
\rput(!5 0 1 5 \RSliceCoords 0 0.1 \Trans){\tiny $(5,0,-1)$}
\pscircle(!0 5 1 5 \RSliceCoords){0.04}
\rput(!0 5 1 5 \RSliceCoords 0 0.1 \Trans){\tiny $(0,5,-1)$}
\pscircle(!2 1 1 5 \RSliceCoords){0.04}
\rput(!2 1 1 5 \RSliceCoords 0 0.1 \Trans){\tiny $(2,1,-1)$}
\pscircle(!1 2 1 5 \RSliceCoords){0.04}
\rput(!1 2 1 5 \RSliceCoords 0 0.1 \Trans){\tiny $(1,2,-1)$}
\psset{fillstyle=solid,fillcolor=white}
\pscircle(!1 0 \ArcCoords){0.04} 
\pscircle(!0 1 \ArcCoords){0.04} 
\pscircle(!5 2 \ArcCoords){0.04} 
\pscircle(!2 5 \ArcCoords){0.04} 
\pscircle(!3 5 sqrt add 4 \ArcCoords){0.04}  
\pscircle(!3 5 sqrt sub 4 \ArcCoords){0.04}  
\pscircle(!4 3 5 sqrt add \ArcCoords){0.04}  
\pscircle(!4 3 5 sqrt sub \ArcCoords){0.04}  
\psset{fillcolor=lightgray}
\pscircle(!3 1 1 5 \RSliceCoords){0.04}
\pscircle(!5 5 3 5 \RSliceCoords){0.04}
\pscircle(!1 3 1 5 \RSliceCoords){0.04}
\SLabel{0.333333333}{\tiny $\frac{1}{3}$}
\SLabel{1}{\tiny $1$}
\SLabel{3}{\tiny $3$}
\end{pspicture}
\end{tabular} \\
$r=1$ & & $r=3$ & & $r=5$ \rule{0cm}{1.0cm}\\[3mm]
\multicolumn{5}{c}{\Fig\label{fig:cones-for-1-3-5}} \\
\end{tabular}

For a line bundle $L=\Osh_{Y}(e_1,e_2)$ the corresponding Seshadri constants are~:

\hspace{-1.5cm}
\begin{tabular}{ccc}
\begin{tabular}{c}
$\ep_{1}(L) = \left\{
\begin{array}{cl}
e_2 & \mbox{if $\frac{e_2}{e_1}\in (0,1]$} \\[2mm]
e_1 & \mbox{if $\frac{e_2}{e_1}\in [1,\infty)$} \\
\end{array}
\right.
$
\end{tabular};
&
\begin{tabular}{c}
$\ep_{3}(L) = \left\{
\begin{array}{cl}
e_2 & \mbox{if $\frac{e_2}{e_1}\in (0,\frac{1}{2}]$} \\[2mm]
\frac{e_1+e_2}{3}  & \mbox{if $\frac{e_2}{e_1}\in [\frac{1}{2},2]$} \\[2mm]
e_1 & \mbox{if $\frac{e_2}{e_1}\in [2,\infty)$} \\
\end{array}
\right.
$
\end{tabular};
&
\begin{tabular}{c}
$\ep_{5}(L) = \left\{
\begin{array}{cl}
e_2 & \mbox{if $\frac{e_2}{e_1}\in (0,\frac{1}{3}]$} \\[2mm]
\frac{e_1+2e_2}{5} & \mbox{if $\frac{e_2}{e_1}\in [\frac{1}{3},1]$} \\[2mm]
\frac{2e_1+e_2}{5} & \mbox{if $\frac{e_2}{e_1}\in [1,3]$} \\[2mm]
e_1 & \mbox{if $\frac{e_2}{e_1}\in [3,\infty)$} \\
\end{array}
\right.
$
\end{tabular}.
\end{tabular}

\vspace{0.25cm}
It is interesting to note the similarities between the formulas above, and those in \S\ref{sec:r=2-4-6}.

\bpoint{$\mathbf{r=7}$}
\label{sec:r=7}

As in \S\ref{sec:r=1-3-5} all ample bundles on $\PP^1\times \PP^1$ are outer. 
The main difference in the case $r=7$ from the cases of all other odd $r$ is that the Seshadri constants
of outer bundles (i.e., all ample bundles in this case) are determined by five curve classes. 
In addition to the four curve classes in \eqref{eqn:four-curve-classes} there is an additional class,
which for reasons of consistency in the diagram we label $C_{2\frac{1}{2}}$~:

$$C_{2\frac{1}{2}}:=(28,28,-15).$$

This curve class is the union of $7$ disjoint $(-1)$-classes.  The curve class $C'=4F_1+4F_2-E-E_1$ (i.e., 
bidegree $(4,4)$, multiplicity $2$ at $p_1$, and multiplicity $1$ at $p_2$,\ldots, $p_7$) satisfies
$(C')^2=-1$ and $C'\cdot K_{X}=-1$.  Starting with the linear series $|\Osh_{Y}(2,2)|$, and imposing multiplicity
$2$ at a point $p_1$, the general member of the resulting linear series is irreducible.  Imposing the condition
that the linear series pass through six further general points, we conclude that
the class $C'$ is represented by an irreducible $(-1)$-curve.  
The symmetrization of $C'$, as in \S\ref{sec:V_r-def}, is the class $C_{2\frac{1}{2}}$. 

The argument in \S\ref{sec:Outline-of-outer-odd}, used again in \S\ref{sec:r=1-3-5}, which shows that
$C_1$,\ldots, $C_4$ determine the Seshadri constant of all outer bundles is the following.  Each time one finds
a $(-1)$-curve class, or a class which is a symmetrization of $(-1)$-curves, one draws the corresponding half-plane
$C_i^{\leq 0}$.  Any subsequent such class has to lie in that half plane. 
When $r\in \{1,3,5\}$ the half planes corresponding to $C_1$,\ldots, $C_4$ eliminate the possibility of any other
curve class with negative self-intersection.  When $r\geq 9$ these half planes do not eliminate the possibility of
any other curve class with negative self-intersection, but do eliminate the possibility of any curve class 
with negative intersection which is $K_X$-negative.  The case $r=7$ is intermediate between these behaviours.
Here the classes $C_1$,\ldots, $C_4$ do not eliminate all classes with negative self-intersection, not even in the
half plane $K_{X}^{\leq 0}$.  
However, the addition of the new curve class $C_{2\frac{1}{2}}$ to the list is sufficient to rule out all
further possibilities. 

The picture in the case $r=7$, along with the corresponding Seshadri constants for a line bundle 
$L=\Osh_{Y}(e_1,e_2)$, appear below.

\vspace{0.75cm}

\begin{tabular}{ccc}
\begin{tabular}{c}
\psset{yunit=3cm, xunit=5cm}
\begin{pspicture}(0,-0.3)(1,1)
\psset{linecolor=gray}
\psset{linecolor=black}
\psset{linestyle=dashed,dash = 3pt 2pt, linecolor=lightgray}
\TanToArc{1}{0}{-0.1}{0.5}
\TanToArc{0}{1}{-0.5}{0.1}
\TanToArc{7}{2}{-0.50}{0.2}
\TanToArc{2}{5}{-0.2}{0.50}
\TanToArc{4 7 sqrt add}{9}{-0.1}{0.3}
\TanToArc{4 7 sqrt sub}{9}{-0.3}{0.1}
\TanToArc{9}{4 7 sqrt add}{-0.3}{0.1}
\TanToArc{9}{4 7 sqrt sub}{-0.1}{0.3}
\psset{linestyle=solid}
\SecToArc{1}{0}{7}{2}{-0.1}{1.07}
\SecToArc{0}{1}{2}{7}{-0.1}{1.07}
\SecToArc{4 7 sqrt add}{9}{4 7 sqrt sub}{9}{-0.1}{1.1}
\SecToArc{9}{4 7 sqrt add}{9}{4 7 sqrt sub}{-0.1}{1.1}
\SecToArc{8}{7}{7}{8}{-0.1}{1.1}
\psset{fillcolor=vlgray,linecolor=gray,fillstyle=solid}
\pspolygon(!1 0 \ArcCoords)(!4 1 1 7 \RSliceCoords)(!17 13 8 7 \RSliceCoords)(!13 17 8 7 \RSliceCoords)(!1 4 1 7 \RSliceCoords)(!0 1 \ArcCoords)
\psset{linecolor=black,fillstyle=none}
\DrawArc
\psset{fillstyle=solid,fillcolor=white,linecolor=black,linestyle=solid,linecolor=black}
\psset{fillstyle=solid,fillcolor=black}
\pscircle(!7 0 1 7 \RSliceCoords){0.04}
\rput(!7 0 1 7 \RSliceCoords -0.1 0.15 \Trans){\tiny $(7,0,-1)$}
\pscircle(!0 7 1 7 \RSliceCoords){0.04}
\rput(!0 7 1 7 \RSliceCoords 0.1 0.15 \Trans){\tiny $(0,7,-1)$}
\pscircle(!3 1 1 7 \RSliceCoords){0.04}
\rput(!3 1 1 7 \RSliceCoords -0.1 0.1 \Trans){\tiny $(3,1,-1)$}
\pscircle(!1 3 1 7 \RSliceCoords){0.04}
\rput(!1 3 1 7 \RSliceCoords 0.1 0.1 \Trans){\tiny $(1,3,-1)$}
\pscircle(!28 28 15 7 \RSliceCoords){0.04}
\rput(!28 28 15 7 \RSliceCoords 0 0.1 \Trans){\tiny $(28,28,-15)$}
\psset{fillstyle=solid,fillcolor=white}
\pscircle(!1 0 \ArcCoords){0.04} 
\pscircle(!0 1 \ArcCoords){0.04} 
\pscircle(!7 2 \ArcCoords){0.04} 
\pscircle(!2 7 \ArcCoords){0.04} 
\pscircle(!9 4 7 sqrt add \ArcCoords){0.04}  
\pscircle(!9 4 7 sqrt sub \ArcCoords){0.04}  
\pscircle(!4 7 sqrt add 9 \ArcCoords){0.04}  
\pscircle(!4 7 sqrt sub 9 \ArcCoords){0.04}  
\pscircle(!7 8 \ArcCoords){0.04}  
\pscircle(!8 7 \ArcCoords){0.04}  
\psset{fillcolor=lightgray}
\pscircle(!4 1 1 7 \RSliceCoords){0.04}
\pscircle(!17 13 8 7 \RSliceCoords){0.04}
\pscircle(!13 17 8 7 \RSliceCoords){0.04}
\pscircle(!1 4 1 7 \RSliceCoords){0.04}
\SLabel{0.25}{\tiny $\frac{1}{4}$}
\SLabel{13 17 div}{\tiny $\frac{13}{17}$}
\SLabel{17 13 div}{\tiny $\frac{17}{13}$}
\SLabel{4}{\tiny $4$}
\end{pspicture} \\
\Fig\label{fig:cone-for-7} \\
\end{tabular}
& \rule{2cm}{0cm} &
\begin{tabular}{c}
$ \ep_{7}(L) = \left\{
\begin{array}{cl}
e_2 & \mbox{if $\frac{e_2}{e_1}\in (0,\frac{1}{4}]$} \\[2mm]
\frac{e_1+3e_2}{7} & \mbox{if $\frac{e_2}{e_1}\in [\frac{1}{4},\frac{13}{17}]$} \\[2mm]
\frac{4e_1+4e_2}{15} & \mbox{if $\frac{e_2}{e_1}\in [\frac{13}{17},\frac{17}{13}]$} \\[2mm]
\frac{3e_1+e_2}{7} & \mbox{if $\frac{e_2}{e_1}\in [\frac{17}{13},4]$} \\[2mm]
e_1 & \mbox{if $\frac{e_2}{e_1}\in [4,\infty)$} \\
\end{array}
\right..
$
\end{tabular}
\end{tabular}

\bpoint{Remark}
The blowup of $\PP^1\times\PP^1$ at $r$ general points is the same as the blowup of $\PP^2$ at $(r+1)$ general points.
The effective cone, and the generating $(-1)$-curves, for $\PP^2$ blown up at $\leq (7+1)$ points are well-known.
For example, the generators are listed in \cite[p. 135, Proposition 26.1]{M}

Applying the change of basis formula between the Picard group of $\PP^2$ blown up at $(r+1)$ points, and
the blowup of $\PP^1\times\PP^1$ at $r$ points, one obtains generators for the effective cones of $\PP^1\times\PP^1$
blown up at $r\leq 7$ points.  These $(-1)$-curves, or their symmetrizations, give generators of the symmetric
effective cones in \S\ref{sec:r=2-4-6}, \ref{sec:r=1-3-5}, and \ref{sec:r=7}.  This is another way to arrive 
at the description of the cones given above.
Note that not all $(-1)$-curves, or their symmetrizations, appear
as boundary generators of the symmetrized effective cone.  Any class which does appear has to satisfy some 
restrictive conditions on the multiplicities, see \cite[Theorem 2.6.2]{Dio}.

\newpage

\section{A brief study of the slopes of the $\xi_n$} 
\label{sec:brief-study-of-slopes}

\point For fixed $r$, $r$ even, recall that in \S\ref{sec:full-packings} 
we have defined sequences $\{q_{n,r}\}_{n\in \ZZ}$ by setting $q_{-1,r}=1$, $q_{0,r}=0$,
$q_{1,r}=1$, and then using the recursion \eqref{eqn:sequence-recursion}.    For instance, from the recursion,
$q_{2,r}=\frac{r-2}{2}\left(q_{1,r}-q_{0,r}\right)+q_{-1,r}=\frac{r}{2}$.

In the proof of Theorem \ref{thm:Outer-nef}
we have defined classes $\xi_{n}=\xi_{n,r}$ for all $n\in \ZZ$ by $\xi_{n,r}=T_r^{n}(F_2)$ (in the proof the
dependence on $r$ was omitted from the notation).  Thus in the usual coordinates on $V_r$, 
$$
\begin{array}{rcl}
\xi_{-1,r} & = & (1,0,0)=\left(q_{-1,r},q_{0,r},-\sqrt{\frac{2 q_{-1,r}q_{1,r}}{r}}\right), \\[3mm]
\xi_{0,r} & = & (0,1,0)=\left(q_{0,r},q_{1,r},-\sqrt{\frac{2 q_{0,r}q_{1,r}}{r}}\right), \,\,\mbox{and}\\[2mm]
\xi_{1,r} & = & (1,\frac{r}{2},-1)=\left(q_{1,r},q_{2,r},-\sqrt{\frac{2 q_{1,r}q_{2,r}}{r}}\right). \\
\end{array}
$$
We note that the third coordinate is determined by the first two, since by
Theorem \ref{thm:Tr-transformation} $\xi_{n,r}^2\stackrel{\mbox{\tiny\ref{thm:Tr-transformation}}}{=}F_2^2=0$ for
each $n$.
Since the recursion \eqref{eqn:sequence-recursion} is the one given by the characteristic polynomial of $T_r$, 
we conclude that 

\begin{equation}
\label{eqn:xi-formula} 
\xi_{n,r} = \left(q_{n,r},q_{n+1,r},-\sqrt{\tfrac{2q_{n,r}q_{n+1,r}}{r}}\right)\,\,\mbox{for all $n\in \ZZ$}. 
\end{equation}
In particular, the slope of $\xi_{r,n}$ is $\frac{q_{n+1,r}}{q_{n,r}}$. 

In this section we prove a few results about these slopes for use in \S\ref{sec:Symplectic-Packing} below, as
well as demonstrate the properties of the map $\phi_{r}$ defined in \eqref{eqn:phi-r-def}.

\bpoint{Lemma} 
\label{lem:slopes-of-xi-decreasing}
Let $r\geq 8$ be even.  Then the sequence of slopes $\{\frac{q_{n+1,r}}{q_{n,r}}\}_{n\geq 1}$
is strictly decreasing. 

\hfill
\begin{tabular}{c}
\begin{tabular}{c}
\psset{yunit=3cm, xunit=5cm}
\begin{pspicture}(0,-0.3)(1,0.8)
\psset{linecolor=gray}
\psset{linecolor=black}
\psset{linecolor=gray}
\psline(!1 0 0 \SliceCoords)(!0 1 0 \SliceCoords) 
\psset{fillstyle=solid,fillcolor=white,linecolor=black,linestyle=solid,linecolor=black}
\SLabel{0}{\tiny $0$}
\SLabelInf{\tiny $\infty$} 
\psset{linecolor=black,fillstyle=none}
\DrawArc
\psset{fillstyle=none,linecolor=black,linestyle=dotted}
\psline(!175 768 \ArcCoords)(!175 768 0 \SliceCoords) 
\psline(!37 192 \ArcCoords)(! 37 192 0 \SliceCoords) 
\psline(!7 48 \ArcCoords)(! 7 48 0 \SliceCoords) 
\psline(!1 12 \ArcCoords)(! 1 12 0 \SliceCoords) 
\psline(!12 1 \ArcCoords)(! 12 1 0 \SliceCoords) 
\psline(!48 7 \ArcCoords)(! 48 7 0 \SliceCoords) 
\psline(!192 37 \ArcCoords)(! 192 37 0 \SliceCoords) 
\psline(!768 175 \ArcCoords)(! 768 175 0 \SliceCoords) 
\psset{linecolor=black,fillstyle=solid,fillcolor=white,linestyle=solid}
\pscircle(!3 1 \ArcCoords){0.04} 
\pscircle(!1 3 \ArcCoords){0.04} 
\pscircle(!175 768 \ArcCoords){0.04} 
\pscircle(!37 192 \ArcCoords){0.04} 
\pscircle(!7 48 \ArcCoords){0.04} 
\pscircle(!1 12 \ArcCoords){0.04} 
\pscircle(!0 1 \ArcCoords){0.04} 
\pscircle(!1 0 \ArcCoords){0.04} 
\pscircle(!12 1 \ArcCoords){0.04} 
\pscircle(!48 7 \ArcCoords){0.04} 
\pscircle(!192 37 \ArcCoords){0.04} 
\pscircle(!768 175 \ArcCoords){0.04} 
\rput(!0 1 \ArcCoords 0.08 0 \Trans){\tiny $\xi_{0,r}$} 
\rput(!1 12 \ArcCoords 1 1 0 \SliceCoords -1 \Scale \Trans 1.12 \Scale 1 1 0 \SliceCoords \Trans 0.02 0 \Trans){\tiny $\xi_{1,r}$}
\rput(!7 48 \ArcCoords 1 1 0 \SliceCoords -1 \Scale \Trans 1.12 \Scale 1 1 0 \SliceCoords \Trans 0.02 0 \Trans){\tiny $\xi_{2,r}$}
\rput{-37}(!175 768 \ArcCoords 1 1 0 \SliceCoords -1 \Scale \Trans 1.08 \Scale 1 1 0 \SliceCoords \Trans){\tiny $\cdots$}
\rput(!1 3 \ArcCoords 0 0.07 \Trans){\tiny $v_{\alpha}$} 
\rput(!3 1 \ArcCoords 0 0.07 \Trans){\tiny $v_{\beta}$} 
\rput(!12 1 \ArcCoords 1 1 0 \SliceCoords -1 \Scale \Trans 1.14 \Scale 1 1 0 \SliceCoords \Trans -0.02 0 \Trans){\tiny $\xi_{-2,r}$}
\rput(!48 7 \ArcCoords 1 1 0 \SliceCoords -1 \Scale \Trans 1.14 \Scale 1 1 0 \SliceCoords \Trans -0.02 0 \Trans){\tiny $\xi_{-3,r}$}
\rput{37}(!786 175 \ArcCoords 1 1 0 \SliceCoords -1 \Scale \Trans 1.10 \Scale 1 1 0 \SliceCoords \Trans){\tiny $\cdots$}
\rput(!1 0 \ArcCoords -0.09 0 \Trans){\tiny $\xi_{-1,r}$} 
\rput(!37 192 0 \SliceCoords 0 -0.15 \Trans){\small $\longleftarrow$}
\rput(!37 192 0 \SliceCoords 0 -0.25 \Trans){\tiny slopes decreasing}
\rput(!192 37 0 \SliceCoords 0 -0.15 \Trans){\small $\longleftarrow$}
\rput(!192 37 0 \SliceCoords 0 -0.25 \Trans){\tiny slopes decreasing}
\end{pspicture} 
\end{tabular}\\
\Fig\label{fig:xi-slopes}
\end{tabular}

\vspace{-4cm}
\parshape 1 0cm 10cm 
\bpf
When $r=8$, from \S\ref{sec:r=8} 
we have the explicit formula $q_n=n^2$, which immediately shows that the sequence is decreasing.
When $r\geq 10$, the argument is that it ``follows from the picture''.
As shown in Figure \ref{fig:xi-slopes}, the slopes of $\xi_{1,r}$, $\xi_{2,r}$, $\xi_{3,r}$, \ldots are decreasing.   
The figure also shows that the slopes on the other side are decreasing.  That is, restricted to the set 
$\{n\in \ZZ \st n\leq -1\}$, the function $n\mapsto \frac{q_{n+1,r}}{q_{n,r}}$ is decreasing. 
This second statement also holds when $r=8$, by the explicit formula above.
\epf

The only point where the sequence of slopes $\{\frac{q_{n+1,r}}{q_{n,r}}\}_{n\in \ZZ}$ fails to be decreasing 
is the transition from $n=-1$ (where the slope is $0$) to $n=0$ (where the slope is $\infty$).  To justify that
this picture is always correct, and thus the argument of the proof is correct, we use the properties of 
$\phi_r$, exposed below.   The justification of the picture appears in \S\ref{sec:justification-of-xi-picture}.

\bpoint{Properties of $T_r^s$} It is convenient to define $T_r^{s}$ for all $s\in \RR$, and not only $s\in \ZZ$.
This is possible since $T_r$ is diagonalizable, and all eigenvalues are positive real numbers.  Specifically,
when $r\geq 10$ the vectors $v_{\alpha_r}$, $v_{\beta_r}$, and $v_1$ of \eqref{eqn:va-vb-v1-def} are a basis
of eigenvectors of $V_r$, and we define $T_r^s$ on $V_r$ by setting

\begin{equation}
\label{eqn:Trs-def}
T_r^s(v_{\alpha_r}) = \alpha_r^s\, v_{\alpha_r},\,\,
T_r^s(v_{\beta}) = \beta_r^s\, v_{\beta},\,\,\rule{0.25cm}{0cm}\mbox{and}\rule{0.25cm}{0cm}
T_r^s(v_{1}) = 1^s\, v_{1}=v_{1}.
\end{equation}
The table of intersections in Figure \ref{fig:intersection-table} and the fact that $\alpha_r\cdot\beta_r=1$ shows
that $T_r^s$ preserves the intersection form on $V_r$.   By construction $T_r^{s}$ also fixes $v_{1}$, the class
of $K_X$. 

\bpoint{Properties of $\varphi_r$} 
\label{sec:properties-of-phi-r}
Suppose that $v\in V_r$, and write $v$ as a linear combination of the basis vectors above~:

\begin{equation}
\label{eqn:v-in-coords}
v = a\, v_{\alpha_r} + b\, v_{\beta_r} + c\, v_{1}.
\end{equation}

Assuming that $b\neq 0$, and using the intersections in Figure \ref{fig:intersection-table} we compute that
$$\frac{v\cdot v_{\beta_r}}{v\cdot v_{\alpha_r}}= \frac{a}{b}.$$
For each $s\in \RR$ let $v_{s} = T_r^{s}(v)$ (so $v_0=v$).  Then by \eqref{eqn:Trs-def}
$v_s = a\alpha_r^s\, v_{\alpha_r} + b\beta_r^s\, v_{\beta_r} + c\, v_{1}$,
and so
$$\frac{v_s\cdot v_{\beta_r}}{v_s\cdot v_{\alpha_r}}= \frac{a\alpha_r^s}{b\beta_r^s} = \frac{a}{b}\cdot \alpha_r^{2s}.$$
If, in addition, $a\neq 0$ then
$$\log\left(\tfrac{v_s\cdot v_{\beta_r}}{v_s\cdot v_{\alpha_r}}\right) = 
\log\left(\tfrac{a}{b}\right) + 2s\log\left(\alpha_r\right) = 
\log\left(\tfrac{v\cdot v_{\beta_r}}{v\cdot v_{\alpha_r}}\right)  + 2s\log\left(\alpha_r\right),$$
from which we conclude that 

\begin{equation}
\label{eqn:phi_r-T-action}
\phi_r(v_s) = \phi_r(v)+s.
\end{equation}

\point
\label{sec:properties-of-phi-r-end}
Now suppose that $v$ is on the square-zero cone. Then, in the coordinates from
\eqref{eqn:v-in-coords}, $v^2=\frac{r-8}{r}(2ab-rc^2)=0$.   If in addition $ab=0$ we conclude that $c=0$, and that $v$ 
is a multiple of either $v_{\alpha_r}$ or $v_{\beta_r}$.  Conversely, if $v$ is such that $v^2=0$ and is not a multiple
of $v_{\alpha_r}$ or $v_{\beta_r}$, then $ab\neq 0$, which shows that the formula in \eqref{eqn:phi-r-def} is
well defined. 

As already noted in \S\ref{sec:digression}, for all $\lambda\in \RR^{*}$, $\phi_r(\lambda v)=\phi_r(v)$.
Thus $\phi_r$ is a well-defined function on the square-zero cone 
minus the lines spanned by $v_{\alpha_r}$ and $v_{\beta_r}$, modulo scaling by $\RR^{*}$. 

\bpoint{Applications to the square-zero cone} 
\label{sec:applications-to-square-zero-cone} 
When $r\geq 10$ the plane $K_X^{\perp}$ intersects the square-zero cone transversely (away from zero).  Thus, the quotient above consists of two open arcs, each homeomorphic to $\RR$.

\vspace{-0.3cm}
\hfill
\begin{tabular}{c}
\begin{tabular}{c}
\psset{yunit=3cm, xunit=5cm}
\begin{pspicture}(-0,-0.8)(1,0.8)
\psset{linecolor=gray}
\psset{linecolor=black}
\rput(!1 12 \AlphaVal \ArcCoords 0.4 0 \Trans){\tiny $K_X^{\perp}$}
\psset{fillstyle=solid,fillcolor=white,linecolor=black,linestyle=solid,linecolor=black}
\psset{linecolor=black,fillstyle=none}
\DrawArc
\DrawArcNeg
\psset{linecolor=white,fillstyle=solid,fillcolor=white,linestyle=solid}
\pscircle(!12 \AlphaVal 1 \ArcCoords){0.12} 
\pscircle(!1 12 \AlphaVal \ArcCoords){0.12} 
\psset{linecolor=black}
\pscircle(!1 1 \ArcCoords ){0.04} 
\pscircle(!1 1 \ArcCoords -1 mul){0.04} 
\psset{linecolor=gray,fillstyle=none}
\psline(! 12 \AlphaVal 1 \ArcCoords exch pop -0.1 exch)(! 12 \AlphaVal 1 \ArcCoords exch pop 1.1 exch)
\rput(!1 12 \AlphaVal \ArcCoords 0 -0.09 \Trans){\tiny $v_{\alpha}$} 
\rput(!12 \AlphaVal 1 \ArcCoords 0 -0.09 \Trans){\tiny $v_{\beta}$} 
\rput(!1 1 \ArcCoords  0 -0.06 \Trans){\tiny $v^{+}$} 
\rput(!1 1 \ArcCoords -1 mul 0 0.08 \Trans){\tiny $v^{-}$} 
\psset{fillstyle=none,arrows=->,linecolor=gray}
\DrawPartArcScale{0}{ 1 12 \AlphaVal div 0.02 add }{1.1}
\psset{arrows=-}
\DrawArcNegScale{1.1}
\psset{arrows=<-}
\DrawPartArcScale{1 1 12 \AlphaVal div 0.02 add sub }{1}{1.1}
\psset{arrows=<->}
\DrawPartArcScale{ 1 12 \AlphaVal div 0.08 add}{1 1 12 \AlphaVal div 0.08 add sub }{1.1}
\rput(!1 12 \AlphaVal\ArcCoords 0.1 0.1 \Trans){\tiny\color{gray} $+\infty$} 
\rput(!1 12 \AlphaVal\ArcCoords 0.25 -0.08 \Trans){\tiny\color{gray} $+\infty$} 
\rput(!12 \AlphaVal 1 \ArcCoords -0.1 0.1 \Trans){\tiny\color{gray} $-\infty$} 
\rput(!12 \AlphaVal 1 \ArcCoords -0.25 -0.08 \Trans){\tiny\color{gray} $-\infty$} 
\end{pspicture} 
\end{tabular}\\
\Fig\label{fig:Tr-and-phir}
\end{tabular}

\nopagebreak[1]
\vspace{-5.2cm}
\parshape 1 0cm 10cm 
If $v^2=0$ then $T_r^{s}(v)^2=0$, since $T_r^{s}$ preserves the intersection form. 
Similarly, since $T_r^{s}$ preserves $K_X$, if $v$ is on the arc in $K_{X}^{>0}$ (respectively
in $K_{X}^{<0}$) then so is $T_r^{s}(v)$. 

\parshape 1 0cm 10cm 
If $v$ is on one of the arcs then as observed above, in the coordinates from 
\eqref{eqn:v-in-coords}, $ab\neq 0$, and thus (modulo scaling) as $s\to\infty$, $T_r^{s}(v)\to v_{\alpha_r}$
and as $s\to -\infty$, $T_r^{s}(v)\to v_{\beta_r}$. 

\parshape 4 0cm 10cm 0cm 10cm 0cm 10cm 0cm \textwidth 
If $v$ is as above, and so $ab\neq 0$, the formula for the action of $T_r^s$ shows that $T_r^{s}(v)$ is a scalar 
multiple of $v$ if and only if $s=0$.  Thus, the action of the group $(\RR,+)$ on each of the upper and lower 
arcs, where $s\in \RR$ acts on $v$ via $s\cdot v = T_r^{s}(v)$, is simply transitive. 

Setting $v^{\pm} = (1,1,\pm \sqrt{\frac{2}{r}})$, $v^{+}$ is on the upper arc, and $v^{-}$ is on the lower arc. 
The maps $s\to T_r^{s}(v^{+})$ and $s\to T_r^{s}(v^{-})$ are therefore continuous bijections of $\RR$ with
the upper and lower arcs respectively.  Since $\phi_r(v^{\pm})=0$, 
\eqref{eqn:phi_r-T-action} shows that $\phi_r$ provides a continuous inverse to each of the previous bijections.  
This allows us to put coordinates on the upper and lower arcs.
The situation is illustrated in Figure \ref{fig:Tr-and-phir} above.

\bpoint{Justification of the picture used in the proof of Lemma \ref{lem:slopes-of-xi-decreasing}}
\label{sec:justification-of-xi-picture}
Since $\xi_{n,r}=T_{r}^{n}(\xi_{0,r}) = T_{r}^{n}(F_2)$ for all $n\in \ZZ$, each $\xi_{n+1,r}$ is farther along
the lower arc in the positive direction than $\xi_{n,r}$, where the notion of positive is provided by 
the action of $(\RR,+)$ above.   Given the locations of $\xi_{-1}$ and $\xi_{0}$ on the square zero cone, 
this shows that the picture in Figure \ref{fig:xi-slopes} is correct, and hence the deduction from it 
used in the proof of Lemma \ref{lem:slopes-of-xi-decreasing} is also correct.

We now return to studying the slopes of the $\xi_{n}$. 

\bpoint{Definition} 
\label{def:Jm}
For each positive integer $m$, define $J_m\subseteq \RR$ by 
$$J_m:=\left(m-\tfrac{1}{2},m\right) \cup \left(m,m+\tfrac{1}{2}\right)\cup \{m+2\} 
= \left(\left(m-\tfrac{1}{2},m+\tfrac{1}{2}\right)\setminus \{m\}\right) \cup \{m+2\}.$$
We note that if $m\neq m'$, then $J_{m}\cap J_{m'}=\emptyset$. 

\bpoint{Proposition}
\label{prop:distinct-slopes-by-r}
Let $r$ be even, $r\geq 10$. Then the sequence $\{\frac{q_{n+1,r}}{q_{n,r}}\}_{n\geq 1}$ is contained in
$J_{\frac{r-4}{2}}$.

\bpf
Starting with $q_{-1,r}=1$, $q_{0,r}=0$, $q_{1,r}=1$, and using \eqref{eqn:sequence-recursion} we compute that 
$$q_{2,r} = \tfrac{r}{2},\,\, q_{3,r}=\tfrac{(r-2)^2}{4},\,\,\mbox{and}\,\,
q_{4,r} = \tfrac{r(r-4)^2}{8},$$
which then gives 
$$
\tfrac{q_{2,r}}{q_{1,r}}=\tfrac{r}{2},\,\,
\tfrac{q_{3,r}}{q_{2,r}}=\tfrac{r-4}{2}+\tfrac{2}{r},\,\,\mbox{and}\,\,
\tfrac{r-4}{2}-\tfrac{q_{4,r}}{q_{3,r}}=\tfrac{2r-8}{(r-2)^2}.
$$
Since $\frac{2r-8}{(r-2)^2}>0$, we conclude that $\frac{q_{4,r}}{q_{3,r}}< \frac{r-4}{2}$. 
As $n\to \infty$, $\frac{q_{n+1,r}}{q_{n,r}}\to \alpha_r$.  
By Lemma \ref{lem:slopes-of-xi-decreasing} the sequence $\{\frac{q_{n+1,r}}{q_{n,r}}\}_{n\geq 1}$ is decreasing.
Thus the previous calculations show that 
$\{\frac{q_{n+1,r}}{q_{n,r}}\}_{n\geq 1}$ is contained in the set 

\begin{equation}
\label{eqn:intermediate-set} 
\left(\alpha_r,\tfrac{r-4}{2}\right) \cup \left\{\tfrac{r-4}{2}+\tfrac{2}{r}\right\} \cup 
\left\{\tfrac{r}{2}\right\}.
\end{equation}

Lemma \ref{lem:alpha-estimate} gives the estimate $\frac{r-5}{2}\leq \alpha_r$, and so the set in 
\eqref{eqn:intermediate-set} is contained in the set 

\begin{equation}
\label{eqn:intermediate-set-2}
\left(\tfrac{r-5}{2},\tfrac{r-4}{2}\right) \cup \left\{\tfrac{r-4}{2}+\tfrac{2}{r}\right\} \cup 
\left\{\tfrac{r}{2}\right\}.
\end{equation}

In turn, for $m=\frac{r-4}{2}$, the set in \eqref{eqn:intermediate-set-2} is contained in $J_m$, 
proving the proposition.
\epf

\bpoint{Theorem (Uniqueness of $r$)} 
\label{thm:uniqueness-of-r}
Let $e_1$ and $e_2$ be positive integers.  Then there is at most one even $r$, $r\geq 2$, such that 
$\frac{e_2}{e_1}=\frac{q_{n+1,r}}{q_{n,r}}$ for some $n\in \ZZ$. If $r\geq 8$ then the value of $n$ is also unique. 

\bpf
By symmetry we may restrict to the case $\frac{e_2}{e_1}\geq 1$.  If 
$\frac{e_2}{e_1}=\frac{q_{n+1,r}}{q_{n,r}}$ and $r\geq 8$, this implies that $n\geq 1$.  
For $r=2$, $4$, and $6$ the sequence is periodic; see \S\ref{sec:r=2-4-6}.

By Proposition \ref{prop:distinct-slopes-by-r} for each even $r\geq 10$ the slopes 
$\{\frac{q_{n+1,r}}{q_{n,r}}\}_{n\geq 1}$ are contained in $J_{\frac{r-4}{2}}$. Since $J_{m}\cap J_{m'}=\emptyset$
if $m\neq m'$, the slopes for different $r$, $r\geq 10$ do not coincide. 

The possible slopes ($\geq 1$) when $r=2$, $4$, $6$, and $8$ are~:
$1$ ($r=2$); $2$ ($r=4$); $\frac{4}{3}$ and $6$ ($r=6$); and $\left\{\frac{(n+1)^2}{n^2}\right\}_{n\geq 1}$ ($r=8$).
See \S\ref{sec:r=2-4-6}--\S\ref{sec:r=8}.

These slopes are distinct, and none are contained in the sets $J_m$ for $m\geq 3$. 
This proves uniqueness of $r$.

If $r\geq 8$ the uniqueness of $n$ follows from the fact that, by Lemma \ref{lem:slopes-of-xi-decreasing},
the sequence of slopes $\{\frac{q_{n+1,r}}{q_{n,r}}\}_{n\geq 1}$ is strictly decreasing.
\epf

\bpoint{Remarks} (\Euler{1}) Not all positive rational numbers are the slopes of some $\xi_{n,r}$, i.e., there
are many slopes $\frac{e_2}{e_1}$ which are not of the form $\frac{q_{n+1,r}}{q_{n,r}}$ for some $n$ and $r$.
For instance, when $m\geq 3$ the only rational number in $(m,m+\frac{1}{2})$ which is a slope of some $\xi_{n,r}$
is $\frac{(m+1)^2}{(m+2)}$ (i.e., $\frac{r-4}{2}+\frac{2}{r}$, where $r=2m+4$), the slope of $\xi_{2,2m+4}$.

(\Euler{2}) The description above of the slopes, particularly the partition of the slopes into the sets $J_m$, can be 
used to give an algorithm to decide whether a particular slope $\frac{e_2}{e_1}$ is the slope of some $\xi_{n,r}$.

\section{Other arguments related to the symplectic packing problem} 
\label{sec:Symplectic-Packing}

In this section we use the results from \S\ref{sec:Outer-Even}--\S\ref{sec:brief-study-of-slopes}
to prove Theorems \ref{thm:symplectic-packing-results} and \ref{thm:full-packings}.
We also give the formulae for the packing constants when $r\leq 8$. 
We start by recalling the dictionary between the differential-geometric language and the algebro-geometric one.

\point 
If $M$ is the real manifold underlying $\PP^1\times \PP^1$, and $\omega_{M}$ a symplectic form on $M$, then
the cohomology class of $\omega_{M}$ is an element of $H^2(M,\RR)\cong \RR^2$, with a natural basis coming from
the K\"{u}nneth theorem.  In this basis the cohomology class of $\omega_M$ corresponds to a pair $(e_1,e_2)$.
By \cite[Theorem 1.1]{LM}, up to diffeomorphism of $M$, 
every such form $\omega_{M}$ comes from a K\"{a}hler form on $\PP^1\times\PP^1$, and for such a form $e_1$, $e_2>0$.
It is these numbers in the formulas below which give the packing constant associated to $\omega_{M}$.

%

\point 
Recall that by the theorem of Biran \cite[Theorem 6.A]{B1} for $Y=\PP^1\times\PP^1$ and $L$ a real ample class of
type $(e_1,e_2)$ (i.e., $e_1$, $e_2\in \RR$, $e_1$, $e_2>0$) one has

\vspace{-0.15cm}
\begin{equation}
\label{eqn:Biran-formula}
\nu_{r}(L) = \left(\frac{\eptilde_{r}(L)}{\eta_{r}(L)}\right)^2 = \frac{r(\eptilde_{r}(L))^2}{2e_1e_2}.
\end{equation}
Here $\eptilde_r(L)$ is like the Seshadri constant
but restricting the test curves to be $(-1)$-curves or, equivalently, their symmetrizations (see
\eqref{eqn-ep-tilde-r-def} for the definition of $\eptilde_{r}(L)$).  

In the calculations in \S\ref{sec:Outer-Even}--\S\ref{sec:Small-r} we have found all the $(-1)$-curves, or 
symmetrizations of $(-1)$-curves, which can affect a Seshadri constant, and thus can compute $\eptilde_{r}(L)$ 
for all $r$.  Reversing these formulae we find the conditions on $r$ for when $L$ of type $(e_1,e_2)$ admits
a full-packing. 

\newpage
\bpoint{Proof of Theorem \ref{thm:symplectic-packing-results}}

{\bf The case of odd $r$.} 
In \S\ref{sec:Outer-Odd} we have shown that the curve classes $C_1$, \ldots, $C_4$ are the only $(-1)$-curves, 
or symmetrizations of such, which affect Seshadri constants when $r\geq 9$ is odd.  

The curves $C_1$,\ldots, $C_4$
affect the Seshadri constants of bundles of slopes $(0,\frac{2}{r+1}]$, $[\frac{2}{r+1},\frac{2}{(\sqrt{r}-1)^2}]$,
$[\frac{(\sqrt{r}-1)^2}{2},\frac{r+1}{2}]$, and $[\frac{r+1}{2},\infty)$ respectively.

In \S\ref{sec:Outer-Odd} we were unable to show that $C_2$ and $C_3$ computed the Seshadri constant over their
entire respective intervals, and instead restricted ourselves to smaller
intervals where we could justify this.   However, for $\eptilde_{r}(L)$ there are no other competing curve classes 
to consider.
Thus, $\eptilde_r(L)$ is computed by the curves $C_1$,\ldots, $C_4$ on the respective intervals listed above, 
and is equal to $\nu_r(L)$ on the intermediate interval $[\frac{2}{(\sqrt{r}-1)^2},\frac{(\sqrt{r}-1)^2}{2}]$.   
This gives \eqref{eqn:Seshadri-tilde-odd-n}.

{\bf The case of even $r$.} 
In the proof of Theorem \ref{thm:Outer-nef} we showed that the only $(-1)$-curves (or symmetrizations) which 
affected the Seshadri constants of line bundles were the curve classes $C_{n}=C_{n,r}=T_r^{n}(E)$, and that these
curve classes computed the Seshadri constants for outer bundles (i.e., bundles whose slope is outside of
$[\beta_r,\alpha_r]$) and did not affect any inner bundles (bundles whose slope is in $[\beta_r,\alpha_r]$).
Thus $\eptilde_r(L)=\ep_{r}(L)$ for outer bundles, and for inner bundles $\eptilde_r(L)=\eta_{r}(L)$.
This gives \eqref{eqn:Seshadri-tilde-even-n}. \epf

\bpoint{Packing constants for $r\leq 8$} 
\label{sec:packing-constants-small-r}
In \S\ref{sec:Small-r} we computed the Seshadri constants when $r\leq 8$ for all ample $L$, and all Seshadri constants 
were computed by $(-1)$-curves or their symmetrizations.  Thus, for $r\leq 8$ we have $\eptilde_{r}(L)=\ep_{r}(L)$
for all ample $L$.   By \eqref{eqn:Biran-formula}, for $L$ a real ample bundle of type $(e_1,e_2)$, we therefore have
$\nu_{r}(L) = \frac{r(\ep_r(L))^2}{2e_1e_2}$ for all $r\leq 8$.

For convenience we list the formulae for those packing constants here.

\hspace{-2.0cm}
\begin{tabular}{ccc}
\begin{tabular}{c}
$\nu_{1}(L) = \left\{
\begin{array}{cl}
\frac{e_2}{2e_1} & \mbox{if $\frac{e_2}{e_1}\in (0,1]$} \\[2mm]
\frac{e_1}{2e_2} & \mbox{if $\frac{e_2}{e_1}\in [1,\infty)$} \\
\end{array}
\right.
$
\end{tabular};
&
\begin{tabular}{c}
$\nu_{3}(L) = \left\{
\begin{array}{cl}
\frac{3e_2}{2e_1} & \mbox{if $\frac{e_2}{e_1}\in (0,\frac{1}{2}]$} \\[2mm]
\frac{(e_1+e_2)^2}{6e_1e_2}  & \mbox{if $\frac{e_2}{e_1}\in [\frac{1}{2},2]$} \\[2mm]
\frac{3e_1}{2e_2} & \mbox{if $\frac{e_2}{e_1}\in [2,\infty)$} \\
\end{array}
\right.
$
\end{tabular};
&
\begin{tabular}{c}
$\nu_{5}(L) = \left\{
\begin{array}{cl}
\frac{5e_2}{2e_1} & \mbox{if $\frac{e_2}{e_1}\in (0,\frac{1}{3}]$} \\[2mm]
\frac{(e_1+2e_2)^2}{10e_1e_2} & \mbox{if $\frac{e_2}{e_1}\in [\frac{1}{3},1]$} \\[2mm]
\frac{(2e_1+e_2)^2}{10e_1e_2} & \mbox{if $\frac{e_2}{e_1}\in [1,3]$} \\[2mm]
\frac{5e_1}{2e_2} & \mbox{if $\frac{e_2}{e_1}\in [3,\infty)$} \\
\end{array}
\right.
$
\end{tabular};
\end{tabular}

\vspace{0.5cm}
\hspace{-2.0cm}
\begin{tabular}{ccc}
\begin{tabular}{c}
$\nu_{2}(L) = \left\{
\begin{array}{cl}
\frac{e_2}{e_1} & \mbox{if $\frac{e_2}{e_1}\in (0,1]$} \\[2mm]
\frac{e_1}{e_2} & \mbox{if $\frac{e_2}{e_1}\in [1,\infty)$} \\
\end{array}
\right.
$
\end{tabular};
&
\begin{tabular}{c}
$\nu_{4}(L) = \left\{
\begin{array}{cl}
\frac{2e_2}{e_1} & \mbox{if $\frac{e_2}{e_1}\in (0,\frac{1}{2}]$} \\[2mm]
\frac{2(e_1+e_2)^2}{9e_1e_2}  & \mbox{if $\frac{e_2}{e_1}\in [\frac{1}{2},2]$} \\[2mm]
\frac{2e_1}{e_2} & \mbox{if $\frac{e_2}{e_1}\in [2,\infty)$} \\
\end{array}
\right.
$
\end{tabular};
&
\begin{tabular}{c}
$\nu_{6}(L) = \left\{
\begin{array}{cl}
\frac{3e_2}{e_1} & \mbox{if $\frac{e_2}{e_1}\in (0,\frac{1}{3}]$} \\[2mm]
\frac{3(e_1+2e_2)^2}{25e_1e_2} & \mbox{if $\frac{e_2}{e_1}\in [\frac{1}{3},\frac{3}{4}]$} \\[2mm]
\frac{12(e_1+e_2)^2}{49e_1e_2} & \mbox{if $\frac{e_2}{e_1}\in [\frac{3}{4},\frac{4}{3}]$} \\[2mm]
\frac{3(2e_1+e_2)^2}{25e_1e_2} & \mbox{if $\frac{e_2}{e_1}\in [\frac{4}{3},3]$} \\[2mm]
\frac{3e_1}{e_2} & \mbox{if $\frac{e_2}{e_1}\in [3,\infty)$} \\
\end{array}
\right.
$
\end{tabular};
\end{tabular}

\hspace{-2.0cm}
\begin{tabular}{ccc}
\begin{tabular}{c}
$\nu_{7}(L) = \left\{
\begin{array}{cl}
\frac{7e_2}{2e_1} & \mbox{if $\frac{e_2}{e_1}\in (0,\frac{1}{4}]$} \\[2mm]
\frac{(e_1+3e_2)^2}{14e_1e_2} & \mbox{if $\frac{e_2}{e_1}\in [\frac{1}{4},\frac{13}{17}]$} \\[2mm]
\frac{56(e_1+e_2)^2}{225e_1e_2} & \mbox{if $\frac{e_2}{e_1}\in [\frac{13}{17},\frac{17}{13}]$} \\[2mm]
\frac{(3e_1+e_2)^2}{14e_1e_2} & \mbox{if $\frac{e_2}{e_1}\in [\frac{17}{13},4]$} \\[2mm]
\frac{7e_1}{2e_2} & \mbox{if $\frac{e_2}{e_1}\in [4,\infty)$} \\
\end{array}
\right.;
$
\end{tabular}
& and & 
\begin{tabular}{c}
$\nu_{8}(L) = \left\{
\begin{array}{cl}
1 & \mbox{if $e_1=e_2$} \\[4mm]
\frac{n^2(n+1)^2(e_1+e_2)^2}{(2n^2-1)^2e_1e_2} & \mbox{if $e_1\neq e_2$, with $n=\left\lceil\frac{1}{\sqrt{\frac{e_2}{e_1}}-1}\right\rceil$} \\ 
\end{array}\right..$
\end{tabular}
\end{tabular}

\vspace{0.5cm}

\bpoint{Conditions for a full packing} 
By \eqref{eqn:Biran-formula} one has a full packing if and only if $\eptilde_{r}(L)=\eta_{r}(L)$.  In terms of our
graphical description of the Seshadri constants (c.f.  \S\ref{sec:graphical-arguments}({\em c}) or
Figure \ref{fig:Closeup-of-Cn}) this means that the bundle $L_{\gamma}$, with $\gamma=\eptilde_{r}(L)$ has
reached the square-zero cone without crossing any plane of the form $C^{\perp}$, where $C$ is a $(-1)$-curve
or symmetrization of a $(-1)$-curve. 
If we are in a region where the Seshadri constant is computed by such curves, this means that $L_{\gamma}$ is
a nef class on the square-zero cone. 

{\bf Proof of Theorem \ref{thm:full-packings} when $r$ is odd.}
In \S\ref{sec:Small-r} we
have seen that for $r$ odd, $r\leq 7$ (where all Seshadri constants are determined by $(-1)$-curves or
their symmetrizations) there are no nef square-zero classes. In other words, in the pictures in \S\ref{sec:Small-r}, the
nef cone never reaches the square-zero cone, although this is a bit difficult to see in the picture for $r=7$.
Thus, to have a full packing when $r$ is odd, one needs at least $r\geq 9$. 

When $r\geq 9$ and odd we have seen in \S\ref{sec:Outer-Odd}, that the curve classes $C_1$,\ldots, $C_4$ only affect 
Seshadri constants for bundles with slopes outside $[\frac{2}{(\sqrt{r}-1)^2},\frac{(\sqrt{r}-1)^2}{2}]$, 
and that outside that interval the nef cone never reaches the square-zero cone. 

Thus, when $r$ is odd a full packing occurs for a real line bundle of type $(e_1,e_2)$ if and only if $r\geq 9$ and 
$\frac{2}{(\sqrt{r}-1)^2}\leq \frac{e_2}{e_1}\leq \frac{(\sqrt{r}-1)^2}{2}$. This is equivalent to the condition 
$$r\geq \max\left(\left(\sqrt{\tfrac{2e_2}{e_1}}+1\right)^2, 
\left(\sqrt{\tfrac{2e_1}{e_2}}+1\right)^2,9\right)$$
of Theorem \ref{thm:full-packings}.

\vspace{0.25cm}
{\bf Proof of Theorem \ref{thm:full-packings} when $r$ is even.}
In the proof of Theorem \ref{thm:Outer-nef}, when $r$ is even, $r\geq 10$, 
we have seen that the curve classes $C_{n,r}:=T_r^{n}(E)$, each the union of $r$ disjoint $(-1)$-curves, 
determine the Seshadri constants for all outer bundles, i.e., bundles whose slope is outside $[\beta_r,\alpha_r]$,
and that no $(-1)$-curve or symmetrization affects the Seshadri constant of inner bundles. 

We have also seen that when $r\geq 10$ the classes $\xi_{n,r}$ are the only nef classes on the square-zero cone
with slope outside of $[\beta_r,\alpha_r]$.  In \S\ref{sec:Small-r}, when $r$ is even, $r\leq 8$,
we have similarly seen that the classes $\xi_{n,r}$ (and $v_{\alpha_8}$) are the only nef classes on the square-zero
cone. 

Thus, when $r$ is even a full packing occurs for a real line bundle of type $(e_1,e_2)$ if and only if 

\RomanList
\begin{enumerate}
\item $\beta_r\leq \frac{e_2}{e_1}\leq \alpha_r$, or 
\item $r$ is a value for which $\frac{e_2}{e_1}$ is equal to the slope of a $\xi_{n,r}$ for some $n$.
\rule{0cm}{0.5cm}
\end{enumerate}

These are equivalent to the conditions 

\RomanList
\begin{enumerate}
\item $r\geq \dfrac{2(e_1+e_2)^2}{e_1e_2}$, or 
\item $r$ is a value for which $\frac{e_2}{e_1}$ is equal to $\frac{q_{n+1,r}}{q_{n,r}}$ for some $n$,
\rule{0cm}{0.5cm}
\end{enumerate}
\AlphaList

appearing in Theorem \ref{thm:full-packings}. 

In addition, we recall that by Theorem \ref{thm:uniqueness-of-r}, for a fixed $(e_1,e_2)$ there is at most
one value of $r$ where condition ({\em ii}) occurs. \epf

This completes the arguments related to the symplectic packing problem.

\section{Petrakiev reflections}
\label{sec:Reflections}

\point 
A common technique for investigating the cone of effective classes on a surface
blown up at $r$ general points is to specialize the points to lie on some fixed curve $G$ in such a way
that the proper transform $G_0$ of $G$ has negative self intersection.    

If $C$ is a class which is effective
when the points are in general position, then 
under the condition that $G_0\cdot G_0<0$ 
the class $C - \left(\frac{C\cdot G_0}{G_0\cdot G_0}\right) G_0$ is
effective when the points are in special position, and this can be used to deduce restrictions on $C$. 

A beautiful argument of I.\ Petrakiev allows one, under mild conditions on $G$ and $r$,  
to double the coefficient of $G_0$ subtracted in the formula above.  
The formula then becomes that for the reflection of $C$ in $G_0$.  Dually, one may
reflect classes which are nef on the specialization to get classes which are nef when the points are in general
position.

In this section we record this result, and in the next use it to produce inner nef classes on the square-zero cone 
for even and odd $r$. 

\bpoint{Theorem (Reflection Theorem)}
\label{thm:Petrakiev-reflection}
Let $Y$ be a smooth surface, and $G\subset Y$ a smooth irreducible curve.  We use $X$ to denote the blowup of $Y$ at
$r$ general points, and $X_0$ to denote the blowup of $Y$ at $r$ general points of $G$.  
We identify the N\'eron-Severi groups of $X$ and $X_0$, along with their intersection forms, via the isomorphisms
$$ \NS(X)_{\QQ} \cong \NS(Y) \bigoplus \oplus_{i=1}^{r} \QQ E_i \cong \NS(X_0)_{\QQ}, $$
and use $V$ to denote the common inner product space. 
We additionally assume the following numerical conditions :
that $r \geq |G\cdot G|$; that $(G\cdot G)-r < 0$ (thus $r=|G\cdot G|$ is
only allowed when $G\cdot G <0$); and, if $G$ has genus $0$, that $(G\cdot G)-r$ is even.  
We denote by $G_0$ the proper transform of $G$ in $X_0$ (therefore $G_0\cdot G_0=(G\cdot G)-r <0$), 
and by $\phi_{G_0}\colon V\longrightarrow V$ the isometry
``reflection in $G_0$'' given by the formula

\begin{equation}
\phi_{G_0}(\xi) := \xi - 2\left(\frac{\xi\cdot G_0}{G_0\cdot G_0}\right)G_0.
\end{equation}
Under the numerical conditions above, 

\AlphaList
\begin{enumerate}
\item If $C$ is an effective class on $X$, then $\phi_{G_0}(C)$ is an effective class on $X_0$.
\item If $\xi_0$ is a nef class on $X_0$, then $\xi:=\phi_{G_0}(\xi_0)$ is a nef class on $X$,
and $\xi\cdot \xi = \xi_0\cdot \xi_0$.
\end{enumerate}

\bpf We first show that ({\em a}) implies ({\em b}). 
We recall that

\RomanList
\begin{enumerate}
\item $\phi_{G_0}$ is an isometry, i.e., $\phi_{G_0}(\xi_1)\cdot\phi_{G_0}(\xi_2)=\xi_1\cdot \xi_2$ for
all $\xi_1$, $\xi_2\in V$, and
\item $\phi_{G_0}$ is self adjoint, i.e., $\xi_1\cdot\phi_{G_0}(\xi_2)=\phi_{G_0}(\xi_1)\cdot \xi_2$ for
all $\xi_1$, $\xi_2\in V$.
\end{enumerate}

(The two statements are equivalent, since $\phi_{G_0}^2 = \Id_{V}$.) 

If $C$ is an effective class on $X$, then by ({\em a}), $\phi_{G_0}(C)$ is an effective class on $X_0$, 
and thus, since $\xi_0$ is a nef class on $X_0$, $\phi_{G_0}(C)\cdot \xi_0\geq 0$.
Therefore
$$C\cdot \xi = C\cdot \phi_{G_0}(\xi_0) \xeq{\mbox{\small ({\em ii})}} \phi_{G_0}(C)\cdot \xi_0\geq 0.$$
Since $C\cdot \xi\geq 0$ for all effective classes $C$ on $X$ we conclude that $\xi$ is a nef class on $X$. 
Furthermore, $\xi\cdot \xi = \xi_0\cdot \xi_0$ since $\phi_{G_0}$ is an isometry.

Part ({\em a}) is proved in \cite{P}, although the result is not explicitly stated in that form.  We
will outline the argument, indicating where in \cite{P} one can find the proofs of the claims. 

To study the degeneration one starts with the product $\cY:=Y\times \Delta$, with $\Delta \subset \CC$ the unit disc.
Before blowing up, the normal bundle of $G\times \{0\}$ in $\cY$ is 
$N_{G/Y} \oplus \Osh_{G}$, where $N_{G/Y}$ is the normal bundle of $G$ in $Y$. 

To blow up, one chooses sections $p_i(t)$, $i=1$,\ldots $r$, in $\cY$, which are general points of $Y$ for general 
$t\in \Delta$, and are general points of $G$ when $t=0$.  Letting $\cX$ be the blowup of $\cY$ along the sections, 
one checks that the normal bundle of $G_0\times \{0\}$ in $\cX$ is obtained by performing $r$ elementary transformations
on $N_{G/Y}\oplus \Osh_{G}$.    Under the numerical conditions above, 
($r\geq |G\cdot G|$; $(G\cdot G)-r<0$ and even if $G$ has genus zero),
a generic choice of degeneration will ensure that 
these elementary transformations result in a semistable bundle. 
This result is a combination of \cite[Lemma 3.4, Corollary 3.5, and Lemma 3.6]{P}. 

Choose such a generic degeneration and let $\cN$ be the resulting normal bundle of $G_0$ in $\cX$.   
Since $\deg_{G_0}(\cN)=\deg(N_{G/Y})-r=G_0\cdot G_0$, $\cN$ is a semistable bundle
of slope $(G_0\cdot G_0)/2$. 

Next, if $C$ is an effective class on the blowup of $r$ general points, we may choose a family of effective curves 
$\cC_t$, $t\in \Delta$,  such that each $\cC_t\subset \cX_t$ has class $C$.   
Let $\mult_{G_0}(\cC_0)$ denote the multiplicity of $G_0$ in the limiting curve $\cC_0$. 
Under the condition that $\cN$ (and thus also its dual $\cN^{*}$) is semistable, the statement of 
\cite[Corollary 2.2]{P} is that 
$$\slope(\cN^{*})\mult_{G_0}(\cC_0) \geq - (C\cdot G_0).$$
Since $\slope(\cN^{*})=-(G_0\cdot G_0)/2$, we conclude that 
$\mult_{G_0}(\cC_0)\geq 2\left(\frac{C\cdot G_0}{G_0\cdot G_0}\right)$, and thus 
that $C-2\left(\frac{C\cdot G_0}{G_0\cdot G_0}\right) G_0$ is an effective class on $X_0$. \epf

\bpoint{Remarks} \RemNum  As is clear from the proof, the extra factor of $2$ comes from the semistability
of $\cN^{*}$ combined with the multiplicity estimate of \cite[Corollary 2.2]{P}. 

\RemNum Since the multiplicity of $G_0$ in $\cC_0$ is an integer, the estimate 
$\mult_{G_0}(\cC_0)\geq 2\left(\frac{C\cdot G_0}{G_0\cdot G_0}\right)$ implies the stronger result that 

\AlphaList
\begin{enumerate}
\item[($\mbox{\em a}'$)] If $C$ is an effective class on $X$, then $
C-\left\lceil{ 2\left(\tfrac{C\cdot G_0}{G_0\cdot G_0}\right)}\right\rceil G_0$
 is an effective class on $X_0$.
\end{enumerate}

Here $\lceil \cdot \rceil$ denotes the round-up. 
We have chosen to record the result without the round-up for several reasons. First, because of the elegance of 
expressing the result as a reflection, which leads easily to the statement in part ({\em b}) of the theorem. 
Second, if $C$ is sufficiently divisible (e.g., if $C\cdot G_0$ is a multiple of $G_0\cdot G_0$) then
the round up makes no difference.  

However, if one knows more precise information about the class $C$ (enough to determine that 
$(C\cdot G_0)/(G_0\cdot G_0)$ is not an integer), then the version above may be useful.

\RemNum If $\xi_0$ is an integral class, $\phi_{G_0}(\xi_0)$ may only be a rational class, and then it is natural
to scale to make $\xi$ integral.  Thus, the fact that $\xi\cdot\xi = \xi_0\cdot\xi_0$ essentially only ensures
that, after scaling, $(m\xi)\cdot (m\xi)$ has the same sign (here meaning $>0$ or $=0$) as $\xi_0\cdot\xi_0$. 
Since we are mostly interested in reflecting square-zero classes, the scaling makes no difference. 

\RemNum Since the multiplicities of $G_0$ are symmetric (they are all $1$), reflection in $G_0$ preserves
the subspace of symmetric classes.   Therefore if $\xi_0$ is a symmetric nef class on $X_0$, $\phi_{G_0}(\xi_0)$
is a symmetric nef class on $X$. 

\RemNum More generally, one could also obtain symmetric nef classes by finding nef classes on $X$ (for instance
by reflection), and restricting them to the subspace $V_r\subset V$ of symmetric classes.     Because $V_r$ is self-dual
under the intersection product, restriction of the linear form defined by a class $\xi$ amounts to orthogonal projection
of $\xi$ onto $V_r$.

Specifically, if $\xi\in V$ is a nef class on $X$, then $\xi$ decomposes as $\xi=\xi_r + \xi^{\perp}$, 
with $\xi_r\in V_r$, and $\xi^{\perp}$ in the orthogonal subspace $V_{r}^{\perp}$ to $V_r$.
($V_{r}^{\perp}$ consists of those classes of the form $\sum m_i E_i$, with $\sum m_i=0$).  
Since the decomposition is orthogonal, we have $\xi^2 = (\xi_r)^2 + (\xi^{\perp})^2$.  
Since the intersection form is negative definite on $V^{\perp}$, 
if $\xi\not\in V_r$  (i.e., if $\xi^{\perp}\neq 0$) then $\xi_r^2 > \xi^2$. 
 
Thus, if one is searching for nef classes in $V_r$ which are square-zero (classes imposing the strongest conditions, 
since they are clearly on the boundary of the nef cone), one cannot do it by restricting nef classes to $V_r$ unless 
they were already in $V_r$ to begin with.

\section{Inner square-zero nef classes via reflections} 
\label{sec:Inner-square-zero-classes-reflections}

\bpoint{Theorem}   
\label{thm:inner-nef-classes}
Let $X$ be the blowup of $\PP^1\times \PP^1$ at $r$ general points, $r\geq 9$, 
and let $e$ be an integer such that
\AlphaList
\begin{enumerate}
\item $2\leq e\leq \frac{r-4}{2}$ if $r$ is even, or
\item $2\leq e < \frac{r}{4}$ if $r$ is odd. \rule{0cm}{0.4cm}
\end{enumerate}

Then the class $\xi(e,r):=(2e^2,r,-2e)$ is an inner nef class on $X$ and
satisfies $(\xi(e,r))^2=0$ (i.e., $\xi(e,r)$ is on the square-zero cone). 
The same statements hold for the class $(r, 2e^2, -2e)$ obtained by switching the bidegrees. 

\bpf
Let $(e_1,e_2)=(e,1)$ if $r$ is even, and $(e_1,e_2)=(e,2)$ if $r$ is odd. 
The class $\xi(e,r)$ is obtained, up to multiple of $r-2e_1e_2$, by reflecting the fibre class $\xi_0:=F_2=(0,1,0)$ 
in the curve $G_0$ obtained by specializing the $r$ general points to lie on a smooth curve $G$ of bidegree 
$(e_1,e_2)$. 

We now check that such a curve $G$ satisfies the numerical conditions of Theorem \ref{thm:Petrakiev-reflection}. 
We recall that a smooth curve of bidegree $(e_1,e_2)$ on $\PP^1\times\PP^1$ has genus $(e_1-1)(e_2-1)$. 

\begin{itemize}
\item If $r$ is even, then $(G\cdot G)-r = 2e-r$ is clearly even, a necessary condition since a curve of bidegree $(e,1)$
has genus $0$.  It is also clear that $(G\cdot G)-r<0$ whenever $e< \frac{r}{2}$, a weaker condition than the
condition $e\leq \frac{r-4}{2}$ being imposed above in ({\em a}). 

\smallskip
\item If $r$ is odd, then since $e\geq 2$, the curve $G$ of bidegree $(e,2)$ has genus $\geq 1$, and thus there
is no restriction on the parity of $(G\cdot G)-r$.  The only condition needed to apply the theorem is that
$(G\cdot G)-r = 4e-r <0$.  This clearly holds whenever $e< \frac{r}{4}$, which is one of the conditions being imposed
in ({\em b}). 
\end{itemize}

Thus we may apply the theorem.
Since the class $F_2$ is clearly nef on the specialization $X_0$, 
part ({\em b}) of Theorem \ref{thm:Petrakiev-reflection}, along with the fact that $r-2e_1e_2>0$  guarantees that 
$\xi(e,r) = (r-2e_1e_2)\phi_{G_0}(\xi_0)$ is nef on $X$. 

It is easy to check directly that $\xi(e,r)$ is square-zero.  On the other hand, this also follows from 
part ({\em b}) of Theorem \ref{thm:Petrakiev-reflection}.
Since $\xi_0=F_2$ is square-zero, so is its reflection $\phi_{G_0}(\xi_0)$, and therefore so is $\xi(e,r)$, since
$\xi(e,r)$ is a multiple of $\phi_{G_0}(\xi_0)$. 

We have 
$\xi(e,r)\cdot K_X = -4e^2+2er-2r = -4\left((e-1)^2-\left(\frac{r-4}{2}\right)(e-1)+1\right).$
Thus, to have $\xi(e,r)\cdot K_X> 0$ (i.e., in order that $\xi(e,r)$ be an inner class), we need
$$(e-1)^2-\left(\tfrac{r-4}{2}\right)(e-1)+1 < 0.$$
With the substitution $t=e-1$, the left hand side of the inequality above becomes the minimal polynomial for $\alpha_r$
and $\beta_r$ from \S\ref{sec:alpha-r-def}.  Thus $\xi(e,r)\cdot K_X>0$ when $e\in (\beta_r+1,\alpha_r+1)$.

Since $\beta_r\in (0,1)$, and $e$ is an integer, the condition $\beta_{r}+1< e$ is equivalent to $2\leq e$.  
If $r\geq 9$ then $e<\frac{r}{4} < \frac{r-3}{2} \leq \alpha_r+1$, where the last inequality comes from 
Lemma \ref{lem:alpha-estimate}.  
Thus when $r$ is odd the conditions in ({\em b}) imply that $e\in (\beta_r+1,\alpha_r+1)$.

If $r$ is even then Lemma \ref{lem:alpha-estimate} shows 
that $\lfloor \alpha_r +1 \rfloor = \frac{r-4}{2}$, and thus since $e$ is an integer (and $\alpha_r$ not an integer
when $r\geq 10$), the condition $e< \alpha_r+1$ is equivalent to the condition $e\leq\frac{r-4}{2}$.
Thus when $r$ is even the conditions in ({\em a}) are equivalent to the condition that $e\in (\beta_r+1,\alpha_r+1)$.

This finishes the proof that the class $\xi(e,r)$ has the properties claimed.  By symmetry the class
$(r, 2e^2, -2e)$ also has these properties.  
\epf

\bpoint{Corollary} 
\label{cor:inner-nef-classes}
Suppose that $r$ is even, $r\geq 10$, and that $e$ is an integer satisfying $2\leq e \leq \frac{r-4}{2}$.  Then
for all $n\in \ZZ$ the classes $T_r^{n}(\xi(e,r))$ are square-zero inner nef classes.

\bpf
Clear since $T_r$ is an automorphism of the problem, preserving nef classes, self intersections, and the canonical
bundle. 
\epf

\bpoint{Remarks} 

\RemNum
In Corollary \ref{cor:inner-nef-classes} the orbits of the classes $\xi(e,r)$ under $T_r$ are distinct, 
but if one also allows scaling by positive integers, then there are fewer equivalence classes.

For instance, when $r=10$, $\xi(2,10)$ and $\xi(3,10)$ (the only values of $e$ possible in this case)
and the classes obtained by switching the bidegrees, are all, up to scaling, in the same orbit.
Similarly, when $r=12$ the classes $\xi(2,12)$, $\xi(3,12)$, and $\xi(4,12)$ and the classes obtained by switching the 
bidegrees are, up to scaling,  in the same orbit.   When $r=14$ there are two equivalence classes, those for which the 
minimal integral ray generator on the ray spanned by $\xi(e,14)$ has intersection $6$ with $K_X$ ($e=2$, $5$), 
and those for which the intersection is $10$ ($e=3$, $4$). 

\RemNum The conclusion of Theorem \ref{thm:inner-nef-classes} also holds when $r=8$.  Then we must have $e=2$,
and $\xi(2,8)$ is the class $(8,8,-4)$, a multiple of $(2,2,-1)$, which is also a multiple of $v_{\alpha_8}=v_{\beta_8}$.
This class is inner, but not strictly inner 
($\alpha_r$ and $\beta_r$ are rational numbers only when $r=8$ or $r=9$, so in the even case
$r=8$ is only case a rational class could have slope $\alpha_r$ or $\beta_r$).   
See also \S\ref{sec:r=8}.  

\RemNum Here is an explanation of the curve classes used in Theorem \ref{thm:inner-nef-classes}. 
Let $G$ have bidegree $(e_1,e_2)$.  In order to apply 
Theorem \ref{thm:Petrakiev-reflection} we need to have $(G\cdot G)-r=2e_1e_2-r<0$.  This puts restrictions on the 
sizes of $e_1$ and $e_2$. 

The one class which is clearly nef on the specialization $X_0$ is the fibre class $F_2$
(or symmetrically, $F_1$).  With $G_0$ the curve obtained as the proper transform of $G$ when the points are specialized
onto $G$, we have
$$\phi_{G_0}(F_2) = \frac{1}{(r-2e_1e_2)}\left(2e_1^2,\, r,\, -2e_1\right).$$
Surprisingly (after scaling), the result only depends on $e_1$.  
Thus, in light of the condition $2e_1e_2< r$, in order to obtain as wide a range of reflections as possible, 
we should make $e_2$ as small as possible. 

When $r$ is even, we may take $e_2=1$, since  $(G\cdot G)-r$ will always be even.  When $r$ is odd, we must 
take $e_2\geq 2$ (and $e_1\geq 2$) in order to ensure that $G$ has genus $\geq 1$.  
From these choices the rest of the conditions in Theorem \ref{thm:inner-nef-classes}({\em a}) follow directly 
from the further requirement that the reflection of $F_2$ be an inner class.

\RemNum 
When $r$ is even, if one takes $e=\frac{r-2}{2}$, then by the arguments in the proof of 
Theorem \ref{thm:inner-nef-classes}, the class $\xi(e,r)$ is a nef {\em outer} class on the square-zero cone.   
These have all been determined in Theorem \ref{thm:Outer-nef}, and all are positive multiples of 
the classes $T_r^{n}(F_2)$, $n\in \ZZ$.  In this case one can check that 
$\xi\left(\frac{r-2}{2},r\right) = 2\cdot T_r^{-2}(F_2)$.

\bpoint{On the possibility of reflecting other classes} 
As remark (\Euler{3}) above explains, Theorem \ref{thm:inner-nef-classes} is based on the fact that the class
$F_2$ is guaranteed to be nef on $X_0$.  When $r$ is even, and the points in general position we know that
the classes $T_r^{n}(F_2)$, $n\in \ZZ$, are nef classes on $X$, and on the square-zero cone.  If we knew that some or 
all of these classes remain nef after we specialize the points (i.e., on $X_0$), then we could reflect those classes,
and obtain other square-zero inner nef classes.     
We conclude this section with some particularly interesting examples of such potential calculations. 

\point In the case $r=10$, let us consider specializing the points onto a curve $G$ of bidegree $(3,1)$.  In this
case, we know that the classes $T_{10}^{n}(F_2)$ are not nef on $X_0$ when $n\leq -2$, since we can compute that those 
particular classes intersect negatively with $G_0$. 

We (the authors) do not know if the classes $T_{10}^{n}(F_2)$ remain nef on $X_0$ when $n\geq 1$.   Suppose for a moment
that they do.  Then, by continuity, the limiting class $v_{\alpha_{10}}$ would also be nef
on $X_0$, and therefore its reflection $\phi_{G_0}(v_{\alpha_{10}})$ would be a square-zero inner nef class 
on $X$.  However, up to multiple, $\phi_{G_0}(v_{\alpha_{10}})$ is the class $(1,1,-\frac{1}{\sqrt{5}})$, and thus
we would have shown the existence of an irrational Seshadri constant, something suspected, but not yet known to
exist. 

That is, the argument above is 

\begin{itemize}
\item If the classes $T_{10}^{r}(F_2)$ are nef on $X_0$ for all $n\geq 1$ (or all sufficiently large $n$), 
then $v_{\alpha_{10}}$ is nef on $X_0$.
\item If $v_{\alpha_{10}}$ is nef on $X_0$ then the class $(1,1,-\frac{1}{\sqrt{5}})$ is nef when the points are 
in general position, and therefore irrational Seshadri constants exist.
\end{itemize}

We do not know if this observation is a step forward in producing an irrational Seshadri constant, or just another
way of hiding the crucial issue.

The surface $X_0$ can also be realized as the blowup of the Hirzebruch surface $\FF_{4}$ at $10$ general points.  
(The curve of self intersection $(-4)$ is $G_0$, and the class $F_1$ is the fibre class of the morphism 
$\FF_4\longrightarrow \PP^1$.)   For the record we include the change of basis matrix 
on the symmetrized N\'{e}ron-Severi lattice for these two realizations of the surface. 
$$M = 
\begin{bNiceMatrix}[first-row,first-col]
& \mbox{\color{gray}\scriptsize $F_1$} & \mbox{\color{gray}\scriptsize $F_2$} & \mbox{\color{gray}\scriptsize $E$} \\
\mbox{\color{gray}\scriptsize $B$} & 0 & \phantom{-}1 & \phantom{-}0 \\
\mbox{\color{gray}\scriptsize $F'$} & 1 & \phantom{-}7 & \phantom{-}10 \\
\mbox{\color{gray}\scriptsize $E'$} & 0 & -1 & -1 \\
\end{bNiceMatrix}
$$
Here $B$ is the class of the curve of self intersection $-4$, $F'$ is the fibre class, and $E'$ the sum of
the exceptional divisors. 

Thus, the question above is whether, in the basis coming from $\FF_4$, the classes $M\cdot T_{10}^{n}(F_2)$, $n\geq 1$, 
are nef on the blowup of $\FF_4$ at $10$ general points.  For instance, when $n=1$ this is asking if the class
$5B+26F'-4E'$ is nef. 

\point 
\label{sec:more-general-question}
More generally, for even $r\geq 10$, if one knew that that the class $v_{\alpha_{r}}$ was nef on the surface
$X_0$ which resulted from specializing the points to lie on a curve of bidegree $(\frac{r-4}{2},1)$ (the last
possible case in ({\em a}) of Theorem \ref{thm:inner-nef-classes}), then one would know that the reflection 
$\phi_{G_0}(v_{\alpha_{r}})$ was nef when the points were in general position.  This reflection is, up to multiple,
the class $(\frac{r-8}{2}, 1, -\sqrt{\frac{r-8}{r}})$, and would again provide an example of an irrational 
Seshadri constant. 

This surface $X_0$ can again be realized as the blowup of the Hirzebruch surface $\FF_4$ at $r$ general points (the
curve $G_0$ again has self-intersection $-4$.)  Thus, it seems very interesting to investigate nef classes
on the blowup of $\FF_4$ at $r$ general points, $r\geq 10$, $r$ even. 
The classes in question all lie in the $K_X$-negative part of the effective cone. However, a curve which
shows that such a class is not nef must be $K_X$-null or $K_{X}$-positive, which is where the difficulty of the 
question lies. 

\point Here is a potentially useful restatement of the above question.  As the proof of Theorem \ref{thm:Outer-nef}
(implicitly) shows, each class $\xi_{n,r}=T_{r}^{n}(F_2)$ is a convex combination of $C_{n,r}$ and $C_{n+1,r}$, where 
$C_{n,r}=T_{r}^{n}(E)$.  Specifically, since $E+T_r(E) = (0,0,1)+(0,r,-1) = (0,r,0) = r\xi_{0,r}$, it follows that

\begin{equation}
\label{eqn:xi-n-in-terms of C-n} 
\xi_{n,r} = \tfrac{1}{r}\left(C_{n,r}+C_{n+1,r}\right) \rule{0.25cm}{0cm}\mbox{for all $n\in \ZZ$}.
\end{equation}

The class $E$ is the disjoint union of the $r$ exceptional divisors, and so each class $C_{n,r}$ is the disjoint
union of $r$ $(-1)$-curves.   If these components of $C_{n,r}$ and $C_{n+1,r}$ remain irreducible when 
specializing the points to general points of a curve of bidegree $\left(\frac{r-4}{2},1\right)$, then it would
follow that $\xi_{n,r}$ is a nef class on the specialization. Thus, the more general question from 
\S\ref{sec:more-general-question} above can be rephrased as~:

{\bf Q :} For some even $r$, $r\geq 10$, is it true that for sufficiently large $n$, the $(-1)$-curves which 
are the components of $C_{n,r}$ remain irreducible when the $r$ points are specialized to general points of
a curve of bidegree $(\frac{r-4}{2},1)$?

One can also rephrase this even more explicitly.  Since the specialization is isomorphic to the blowup of $\FF_4$
at $r$ general points, we can write the $(-1)$-classes in the natural basis from this point of view.  
The question then becomes ``are these $(-1)$-classes irreducible for all sufficiently large $n$?''.

As discussed above, a positive answer, for any fixed even $r\geq 10$, would establish the existence of an irrational 
Seshadri constant.

\newpage
\section{Inner square-zero nef classes via pullbacks} 
\label{sec:Inner-square-zero-classes-pullbacks}

In this section we use pullback maps to produce inner square-zero nef classes. When $r$ has a factor $r_0$
which is even and $\geq 8$, this allows us to produce such classes different from the 
classes in \S\ref{sec:Inner-square-zero-classes-reflections}. 

\point 
\label{sec:pullback-and-pushforward-formulae}
Given positive integers $a$, $b$, let $\phi_{a,b}\colon \PP^1\times\PP^1\longrightarrow \PP^1\times \PP^1$ be a map 
which is of degree $a$ on the first $\PP^1$ factor, and degree $b$ on the second.  The
map $\phi_{a,b}$ has degree $ab$.  For instance, $([x_0:x_1],[y_0:y_1])\mapsto ([x_0^a:x_1^a],[y_0^b:y_1^b])$ 
is such a map.

Fix a positive integer $r_0$, and let $X_{r_0}$ be the blowup of $Y=\PP^1\times\PP^1$ at $r_0$ general points, 
$p_1$, \ldots, $p_{r_0}$.  Since the points are general we may assume that they do not lie in the branch locus 
of $\phi_{ab}$, and so each of the points $p_i$ pulls back to $ab$ distinct points. 

Let $X_{r}$ be the blowup of $\PP^1\times\PP^1$ at the resulting $r:=abr_0$ points, and 
$$\psi_{a,b}\colon X_{r}\longrightarrow X_{r_0}$$
the induced morphism. 

As in the rest of the paper, we use $V_{r}$ and $V_{r_0}$ respectively for the subspaces of $H^2(X_{r},\RR)$
and $H^2(X_{r_0},\RR)$ generated by the fibre classes and the sum of the exceptional divisors.  The morphism
$\psi_{a,b}$ induces pullback and pushforward morphisms between these spaces.  Specifically,

\begin{equation}
\label{eqn:pullback-and-pushforward-formulae}
\begin{array}{lcccr}
& V_{r} & \xleftarrow{\rule{0.75cm}{0cm}} & V_{r_0} & \colon \psi_{a,b}^{*} \\
& (ad_1, bd_2, -m) & \rotatedown{\xmapsto{\rule{0.5cm}{0cm}}} & (d_1, d_2, -m) & \\
\\
\hline
\\
\psi_{a,b*}\colon & V_{r} & \xrightarrow{\rule{0.75cm}{0cm}} & V_{r_0} & \\
& (d_1',d_2',-m') & \xmapsto{\rule{0.5cm}{0cm}} & (bd_1', ad_2', -abm') \\
\end{array}
\end{equation}
Here, as before, the coordinates on $V_{r}$ and $V_{r_0}$ are with respect to the basis 
consisting of the fibre classes and the sum of the exceptional divisors. 

The pullback and pushforward morphisms are adjoint with respect to the inner products on the two spaces.  
For $v=(d_1,d_2,-m)\in V_{r_0}$ and $w=(d_1',d_2',-m')\in V_{r}$,

\begin{equation}
\label{eqn:adjointness-of-pushforward-and-pullback}
\left\langle \phi_{a,b}^{*}(v),\,w\right\rangle_{r} = ad_1d_2' + bd_2d_1' -abr_0mm' 
= \left\langle v,\,\phi_{a,b*}(w)\right\rangle_{r_0}.
\end{equation}
In the equation above we have used $\langle\smallbullet,\smallbullet\rangle$ instead of $\cdot$ for the inner product, 
to allow us to write a subscript indicating on which space the inner product is being evaluated. 

For classes $v_1$, $v_2\in V_{r_0}$, and $w_1$, $w_2\in V_{r}$, we also have 

\vspace{-0.25cm} 
\begin{equation}
\label{eqn:other-intersection-formulae}
\rule{1.5cm}{0cm} \left\langle \psi_{a,b}^{*}(v_1),\, \psi_{a,b}^{*}(v_2)\right\rangle_{r} 
= ab\left\langle v_1,\, v_2\right\rangle_{r_0}
\rule{0.25cm}{0cm}\mbox{and}\rule{0.25cm}{0cm} 
\left\langle \psi_{a,b*}(w_1),\, \psi_{a,b*}(w_2)\right\rangle_{r_0} 
= ab\left\langle w_1,\, w_2\right\rangle_{r}.
\end{equation}

The $r$ points at which we are blowing up are not general.  However, effective classes remain effective
under specialization, and dually, classes which are nef when the points are specialized are nef when the 
points are in general position.  Thus, if $\xi_0$ is a nef class in $V_{r_0}$, $\psi_{a,b}^{*}(\xi_{0})$ is a nef class
on $\PP^1\times\PP^1$ blown up at $r$ general points.  Here we are identifying the intersection spaces
for the blowup at general points, and the blowup at special points as in Theorem \ref{thm:Petrakiev-reflection}.

\newpage
\bpoint{Proposition} 
\label{prop:pullback-of-inner-nef}
Let $\xi_{0}=(d_1,d_2,-m)\in V_{r_0}$ be a point in the octant where $d_1$, $d_2$, and 
$m$ are $\geq 0$.  Fix positive integers $a$ and $b$, and set $\xi=\psi_{a,b}^{*}(\xi_0)\in V_{r}$.

Then
\begin{enumerate}
\item If $\xi_{0}^2=0$ then $\xi^2=0$. 
\item If $\xi_{0}$ is nef, then $\xi$ is a nef class on $\PP^1\times\PP^1$ blown up at $r$ general points.
\item If $\xi_{0}$ is $K_{X_{r_0}}$ positive, then $\xi$ is $K_{X_{r}}$ positive.
\end{enumerate}

That is, if $\xi_{0}$ is square zero, nef, or $K_X$-positive, the same is true of the pullback $\xi$.
In particular, inner square-zero nef classes on $X_{r_0}$ pull back to inner square-zero nef classes on 
$\PP^1\times\PP^1$ blown up at $r$ general points. 

\bpf
({\em a}) 
By \eqref{eqn:other-intersection-formulae} we have $\left\langle \xi,\,\xi\right\rangle_{r} = ab
\left\langle\xi_{0},\,\xi_{0}\right\rangle_{r_0} = ab\,\xi_{0}^2=0$, so $\xi$ is a square-zero class.
({\em b}) The class $\xi$ is nef by the argument above~: when the $r$ points are specialized, 
the class $\xi$ is nef on the specialization, and hence nef when the points are in general position. 

({\em c}) 
To see that $\left\langle K_{X_{r}},\xi\right\rangle_{r}>0$, we use 
\eqref{eqn:adjointness-of-pushforward-and-pullback}, and show that 
$\left\langle \psi_{a,b*}(K_{X_{r}}),\,\xi_{0}\right\rangle_{r_0}>0$. 
Since $K_{X_{r}}=(-2,-2,1)$, and assuming by symmetry that $a\leq b$, we have 
$$
\begin{array}{rcl}
\psi_{a,b*}(K_{X_{r}}) & \xeq{\mbox{\tiny \eqref{eqn:pullback-and-pushforward-formulae}}} & 
(-2b,-2a,ab) \\[2mm]
& = & b\,(-2,-2,1) + (b-a)\,(0,2,0) +  (a-1)b\,(0,0,1) \\[2mm]
& = & b\, K_{X_{r_0}} + 2(b-a)\, F_2 + (a-1)b\, E,
\end{array}
$$
all the coefficients above are $\geq 0$ and $b>0$.  By assumption 
$\left\langle K_{X_{r_0}},\xi_{0}\right\rangle_{r_0}>0$. Furthermore $\langle F_2, \xi_{0}\rangle = d_1$ and
$\langle E, \xi_{0}\rangle = m$, both of which are $\geq 0$ by assumption on the octant. 
Thus $\left\langle K_{X_{r}},\xi\right\rangle_{r}>0$. \epf

\bpoint{Remarks} 
\label{sec:remarks-on-pullback-maps}
(\Euler{1}) 
As clear from the proof of ({\em c}), even a class $\xi_0$ such that $\langle K_{X_{r_0}},\xi_0\rangle_{r_0}<0$
can pull back to a class $\xi$ with $\langle K_{X_{r}},\xi\rangle_{r}>0$, as long as the intersections of $\xi_0$
with $2(b-a)F_2$ and $(a-1)bE$ are sufficiently positive to make up for the negativity of the first intersection. 

In particular as long as $(a,b)\neq (1,1)$ (i.e., as long as $\psi_{a,b}$ is not the identity map), 
the $K_{X_{r_0}}$-null classes $v_{\alpha_{r_0}}$ and $v_{\beta_{r_0}}$pull back to classes which are 
$K_{X_{r}}$ positive.

(\Euler{2}) We also note that if $\xi_0$ is nef, $\xi_0$ must be in the octant with $d_1$, $d_2$, and $m\geq 0$
(\S\ref{sec:airplane-hanger}).

\point 
Let $a$, $b$, $r_0$, and $r=abr_0$ be as above.   Starting with a class $\xi(e_0,r_0)=(2e_0^2,r_0,-2e_0)$ produced by
Theorem \ref{thm:inner-nef-classes}, we compute that 
$$\psi_{a,b}^{*}(\xi(e_0,r_0)) = (2ae_0^2, br_0, -2e_0) = \tfrac{1}{a}(2(ae_0)^2, abr_0,-2ae_0)
= \tfrac{1}{a} \xi(e,r)
$$
where $e=ae_0$.  From the fact that $e_0$ satisfies the inequalities necessary to apply 
Theorem \ref{thm:inner-nef-classes} (i.e., that $2\leq e_0\leq\frac{r_0-4}{2}$ or
$2\leq e_0 < \frac{r_0}{4}$ depending on the parity of $r_0$), we see that $e$ also satisfies the 
corresponding inequalities.

That is, pulling back the classes as produced by Theorem \ref{thm:inner-nef-classes} only gives classes
of the same type, up to scalar multiple.

However, if $r_0$ is even, $r_0\geq 8$, then as recorded in Corollary \ref{cor:inner-nef-classes}, we may apply
powers of $T_{r_0}$ to each $\xi(e_0,r_0)$ to obtain infinitely many other inner nef square-zero classes. 
We may then pull back these classes to $X_{r}$ and apply powers of $T_r$.  In general the pullback of
points in the $T_{r_0}$-orbit of $\xi(e_0,r_0)$ are not in the $T_r$ orbit of the pullback of any $\xi(e_0',r_0)$,
and this allows us to produce infinitely many new inner nef square-zero classes.

To illustrate the idea, we look at one of the smallest possible examples.

\bpoint{Example} Let $r_0=10$, and $(a,b)=(2,1)$, so that $r=abr_0=20$. 
Taking $e_0=2$, $3$ in Theorem \ref{thm:inner-nef-classes}({\em a}) the classes $\frac{1}{2}\xi(2,10)=(4,5,-2)$
and $\frac{1}{2}\xi(3,10)=(9,5,-3)$ are square-zero inner nef classes, as are the classes $(4,5,-2)$ and
$(9,5,-3)$ obtained by switching the first two factors.  Here we have divided by $2$ to remove common factors
among the coordinates, i.e., to replace each class by the integral generator of the ray it spans.

Let us just focus on one of these, the class $\xi_0:=\frac{1}{2}\xi(2,10)=(4,5,-2)$.  

The forward orbit $T_{10}^{n}(\xi_0)$, $n\geq 1$, of $\xi_0$ converges, up to scaling, to $v_{\alpha_{10}}$.
As noted in \S\ref{sec:remarks-on-pullback-maps}, $v_{\alpha_{10}}$ will pull back to a $K_{X}$-positive class
(and also a nef, square zero class, by Proposition \ref{prop:pullback-of-inner-nef}).   The sequence 
$\xi_{n,10}=T_{10}^{n}(F_2)$, $n\geq 1$ also converges, modulo scaling, to $v_{\alpha_{10}}$, although
from the ``other side''.  It follows that for $n$ sufficiently large, the $\xi_{n,10}$ also pull back to 
$K_{X}$-positive classes.

Thus, by pulling back, we obtain infinitely many inner nef square-zero classes converging (modulo scaling)
to $\psi_{2,1}^{*}(v_{\alpha_{10}})$, and converging from both sides. 
The example is illustrated in Figure \ref{fig:pullback-illustration}.

\begin{tabular}{ccc}
\begin{tabular}{c}
\psset{yunit=3cm, xunit=5cm}
\begin{pspicture}(0,-0.2)(1,1.1)
\psset{linecolor=gray}
\psset{linecolor=black}
\psset{linecolor=gray}
\psline(!1 0 0 \SliceCoords)(!0 1 0 \SliceCoords) 
\psline(! 20 \AlphaVal 1 \ArcCoords exch pop 0 exch)(! 20 \AlphaVal 1 \ArcCoords exch pop 1 exch)
\SLabel{0}{\tiny $0$}
\SLabel{1 20 \AlphaVal div}{\tiny $\beta_{20}$} 
\SLabel{20 \AlphaVal }{\tiny $\alpha_{20}$} 
\SLabelInf{\tiny $\infty$} 
\psset{linecolor=black,fillstyle=none}
\DrawArc
\psset{fillstyle=solid,fillcolor=white,linecolor=black,linestyle=solid,linecolor=black}
\pscircle(!2 5 \ArcCoords){0.036} 
\pscircle(!10 9 \ArcCoords){0.034} 
\pscircle(!18 20 \ArcCoords){0.032} 
\pscircle(!40 49 \ArcCoords){0.030} 
\pscircle(!2 5 \ArcCoords){0.036} 
\pscircle(!10 16 \ArcCoords){0.034} 
\pscircle(!32 45 \ArcCoords){0.032} 
\pscircle(!2 10 \AlphaVal \ArcCoords){0.04} 
\pscircle(!20 \AlphaVal 1 \ArcCoords){0.04}
\rput(!20 \AlphaVal 1 \ArcCoords -0.05 0.1 \Trans){$v_{\beta_{20}}$}
\pscircle(!1 20 \AlphaVal \ArcCoords){0.04}
\rput(!1 20 \AlphaVal \ArcCoords 0.05 0.1 \Trans){$v_{\alpha_{20}}$}
\rput(!1 20 \AlphaVal \ArcCoords 0.20 0 \Trans){\tiny $K_X^{\perp}$}
\rput(!2 10 \AlphaVal \ArcCoords 0 0.3 \Trans){\tiny $\psi^{*}_{2,1}(v_{\alpha_{10}})$} 
\psset{fillstyle=none,arrows=->,linecolor=gray}
\DrawPartArcScale{0.45}{0.52}{1.1}
\psset{arrows=<-}
\DrawPartArcScale{0.54}{0.61}{1.1}
\end{pspicture} \\
\end{tabular} 
&
\begin{tabular}{c}
\begin{pspicture}(-2.0,-0.2)(2.0,0.2)
\psline[arrows=->](1.0,0)(-1.0,0)
\rput(0,0.3){\small $\psi_{2,1}^{*}$} 
\end{pspicture}
\end{tabular}
&
\begin{tabular}{c}
\psset{yunit=3cm, xunit=5cm}
\begin{pspicture}(0,-0.2)(1,0.8)
\psset{linecolor=gray}
\psset{linecolor=black}
\psset{linecolor=gray}
\psline(!1 0 0 \SliceCoords)(!0 1 0 \SliceCoords) 
\psline(! 10 \AlphaVal 1 \ArcCoords exch pop 0 exch)(! 10 \AlphaVal 1 \ArcCoords exch pop 1 exch)
\SLabel{0}{\tiny $0$}
\SLabel{1 10 \AlphaVal div}{\tiny $\beta_{10}$} 
\SLabel{10 \AlphaVal }{\tiny $\alpha_{10}$} 
\SLabelInf{\tiny $\infty$} 
\rput(!1 10 \AlphaVal \ArcCoords 0.35 0 \Trans){\tiny $K_X^{\perp}$}
\psset{linecolor=black,fillstyle=none}
\DrawArc
\psset{fillstyle=solid,fillcolor=white,linecolor=black,linestyle=solid,linecolor=black}
\pscircle(!4 5 \ArcCoords){0.036} 
\pscircle(!5 9 \ArcCoords){0.034} 
\pscircle(!9 20 \ArcCoords){0.032} 
\pscircle(!20 49 \ArcCoords){0.030} 
\pscircle(!1 5 \ArcCoords){0.036} 
\pscircle(!5 16 \ArcCoords){0.034} 
\pscircle(!16 45 \ArcCoords){0.032} 
\pscircle(!10 \AlphaVal 1 \ArcCoords){0.04}
\rput(!10 \AlphaVal 1 \ArcCoords 0.12 sub){$v_{\beta_{10}}$}
\pscircle(!1 10 \AlphaVal \ArcCoords){0.04}
\rput(!1 10 \AlphaVal \ArcCoords 0.12 sub){$v_{\alpha_{10}}$}
\psset{fillstyle=none,arrows=->,linecolor=gray}
\DrawPartArcScale{0.55}{0.63}{1.1}
\psset{arrows=<-}
\DrawPartArcScale{0.67}{0.75}{1.1}
\end{pspicture} \\
\end{tabular} \\
\small $V_{20}$ & & \small $V_{10}$ \\
\multicolumn{3}{c}{
\small \Fig \label{fig:pullback-illustration} Illustration of pullback map}
\end{tabular}

\vspace{0.25cm} 

Now we can apply $T_{20}^{m}$ to the pullbacks.  Applying $T_{20}^m$
to $\psi_{2,1}^{*}(v_{\alpha_{10}})$ we obtain a sequence of inner, nef, square-zero classes
converging (modulo scaling) to $v_{\beta_{20}}$ or $v_{\alpha_{20}}$ as $m\to-\infty$ or
$m\to \infty$ respectively.

But, each of the elements in this sequence itself has a sequence of inner, nef, square-zero classes 
converging to it, from both sides.  Specifically, fixing $m$, the sequences
$$T^{m}_{20}(\psi_{2,1}^{*}(T_{10}^{n}(\xi_0)))
\rule{0.25cm}{0cm}\mbox{and}\rule{0.25cm}{0cm} 
T^{m}_{20}(\psi_{2,1}^{*}(\xi_{n,10}))$$ 
converge (modulo scaling) to $T^m_{20}(\psi_{2,1}^{*}(v_{\alpha_10}))$ as $n\to\infty$.  For all $n$ the classes
of the first type are inner, nef, and square-zero classes.  Classes of the second type are nef and
square-zero, and are inner (i.e., $K_{X}$-positive) for sufficiently large $n$.  In this specific example, 
$n\geq 1$ is large enough. 

We can also apply this construction to the three other classes (e.g., $(9,5,-3)$) listed above.

\point 
By repeated pullbacks we can thus arrive at an $r$ where we can find an infinite sequence 
of inner, nef, square zero classes, each member of which has an infinite sequence of such classes converging to it,
and each member of those previous sequences has an infinite sequence of such classes converging to it, \ldots,
and so on, up to a finite number of such steps.

As discussed in \S\ref{sec:inner-bundle-summary}, a consequence of the SGSH conjecture is that some portion
of the nef cone should be round.
If this implication of the SHGH conjecture is not true, then the above examples suggest that the actual 
description of the nef cone is likely to be quite complicated. 

\bpoint{Use in lower bounds for the Seshadri constant} 
\label{sec:comparison-of-lower-bounds}
The paper \cite{DTG} establishes lower bounds on the Seshadri constants for line bundles whose Seshadri constants
are not affected by $(-1)$-curves.  As noted in  \S\ref{sec:other-work-II}, for a line bundle $L$ of type
$(e_1,e_2)$, these bounds are 

\begin{equation}
\label{eqn:bounds-formulae}
\begin{array}{rcll}
\ep_r(L) & \geq & \eta_{r}(L)\left(1-\tfrac{1}{5r}\right)^{\frac{1}{2}} & \mbox{for $r$ odd, 
$\frac{e_2}{e_1}\in [\frac{2}{(\sqrt{r}-1)^2},\frac{(\sqrt{r}-1)^2}{2}]$} \\
\ep_r(L) & \geq & \eta_{r}(L)\left(1-\tfrac{2}{9r}\right)^{\frac{1}{2}} & \mbox{for $r$ even, 
$\frac{e_2}{e_1}\in [\beta_r,\alpha_r]$}.  \\
\end{array}
\end{equation}

In \S\ref{sec:Inner-square-zero-classes-reflections} we have constructed examples of inner
bundles where $\ep_r(L)=\eta_r(L)$ (all of the classes $\xi(e,r)$, and, when $r$ is even, their orbits under $T_r$).
Applying pullbacks we can construct even more such classes when $r$ has an even factor $\geq 8$.   Thus, the
methods of these sections produce examples of bundles whose Seshadri constants are larger than the lower
bounds above.

The convex hull of such classes (for fixed $r$) then also provides a lower bound on the Seshadri constant. This 
lower bound is exact for bundles of the type we have constructed (those where $\ep_r(L)=\nu_r(L)$), 
and improves on the bounds in \eqref{eqn:bounds-formulae}, at least in the neighbourhood of such bundles.

The authors have been unable to find a useful way to describe and organize all the bundles produced by these
procedures, and thus are unable to give a short formula for a better lower bound.   Thus, the bounds
in \eqref{eqn:bounds-formulae} seem, at the moment, to be the most generally useful.  They are also quite strong.
For instance, when $r=20$, the lower bound in \eqref{eqn:bounds-formulae} is that $\ep_r(L)$ is differs
from $\eta_{r}(L)$ by a factor of no worse than $\sqrt{\frac{89}{90}}\approx 0.994428\ldots$. 

We end this section by giving an application of the formulas in \S\ref{sec:pullback-and-pushforward-formulae}
to establish a reasonably strong family of bounds on the symmetric effective cone, valid (at least) whenever $8\mid r$. 

\bpoint{Theorem} 
Suppose that $r_0$ and $\xi_0\in V_{r_0}$ are such that for all effective classes $C\in V_{r_0}$ one has

\begin{equation}
\label{eqn:effective-inequality-r0}
-r_0C^2 \leq (C\cdot \xi_0)^2.
\end{equation}
Then

\begin{enumerate}
\item For any positive $a$, $b$, setting $r=abr_0$ and $\xi=\psi^{*}_{a,b}(\xi_0)$, for
all effective classes $C\in V_{r}$ we similarly have 
\PauseEnumerate

\begin{equation}
\label{eqn:effective-inequality-r}
-rC^2 \leq (C\cdot \xi)^2.
\end{equation}

\ResumeEnumerate
\item If $r_0$ is even, then for each $n\in \ZZ$ \eqref{eqn:effective-inequality-r} holds with $r=r_0$,  
$\xi=T_{r_0}^{n}(\xi_0)$, and for all effective classes $C\in V_{r_0}$.

\smallskip
\item Condition \eqref{eqn:effective-inequality-r0} of the theorem holds when $r_0=8$ and $\xi_0=-K_{X_8}$. 
\end{enumerate}

Here the condition on effectivity means 
``for $\PP^1\times\PP^1$ blown up at $r_0$ (respectively $r$) general points''.

\bpf
({\em a})
Let $C$ be an effective class in $V_r$.
Since $C$ is a class which is effective on $\PP^1\times\PP^1$ blown up at $r$ general points,
then it is also an effective class when blowing up at special points.  Thus, $\psi_{a,b*}(C)$ is also an effective
class in $V_{r_0}$. 

We then have 
$$-r \langle C,C\rangle_r \xeq{\mbox{\tiny \eqref{eqn:other-intersection-formulae}}} 
-r_0 \langle \psi_{a,b*}(C),\psi_{a,b*}(C)\rangle_{r_0}
\stackrel{\mbox{\tiny \eqref{eqn:effective-inequality-r0}}}{\leq}
\langle \psi_{a,b*}(C),(-K_{X_8})\rangle_{r_0}^2
\xeq{\mbox{\tiny \eqref{eqn:adjointness-of-pushforward-and-pullback} }} \langle C,\xi\rangle_{r}^{2},
$$
or $-rC^2\leq (C\cdot\xi)^2$, which was the inequality to be proved in this case. 

({\em b}) This follows from \eqref{eqn:effective-inequality-r} and the fact that, by 
Theorem \ref{thm:Tr-transformation}, $T_{r_0}$ preserves the intersection form and the property of being effective.

({\em c}) 
Let $X_8$ denote $\PP^1\times \PP^1$ blown up at $8$ general points. 
For the curve classes $C_{n}=(4n(n-1),\, 4n(n+1),\, 1-2n^2)$ of \S\ref{sec:r=8}, we compute that
$C_n^2=-8$ and $C_n\cdot(-K_{X_8})=8$.  Thus these curves, which, by the arguments in 
\S\ref{sec:r=8} form the boundary of the effective cone, satisfy $-8(C_n^2)=\left(C_n\cdot(-K_{X_8})\right)^2$.
By convexity, it follows that for all effective classes $C\in V_{8}$ we have 
$ -8C^2 \leq \left(C\cdot (-K_{X_8})\right)^2.$
\epf

\bpoint{Remarks} 
(\Euler{1}) The bound \eqref{eqn:effective-inequality-r} is homogeneous in $C$, but not in $\xi$. 

(\Euler{2}) 
We recall that $-K_{X_8}$ is a square-zero nef class.  Thus by
Proposition \ref{prop:pullback-of-inner-nef} and Theorem \ref{thm:Tr-transformation} 
each class $\xi=T_{r}^{n}(\psi^{*}_{a,b}(-K_{X_8}))$ is also a square-zero nef class, and is an inner
class as long as $ab\neq 1$ (i.e., $r\neq 8$).  

(\Euler{3}) The form of the theorem is set up to be able to iterate the process.  For instance, if $8\mid r$ and
if $\frac{r}{8}$ has many factors, one can choose different combinations of ({\em a}) and ({\em b}) to step from
$r_0=8$ to $r$.

(\Euler{4}) Here is an illustration of how this bound works.   In the diagram below, in $V_{8}$, one can see
the class $\xi_{0}=-K_{X}$, the corresponding line $\xi_{0}^{\perp}$ (tangent to the square-zero cone), and the
dashed curve, also tangent to the square-zero cone at $\xi_0$, which contains the boundary effective classes when 
$r_0=8$.

\begin{tabular}{ccc}
\begin{tabular}{c}
\psset{yunit=3cm, xunit=5cm}
\begin{pspicture}(0,-0.2)(1,1.1)
\psset{linecolor=gray}
\psset{linecolor=black}
\psset{linecolor=gray}
\psline(!1 0 0 \SliceCoords)(!0 1 0 \SliceCoords) 
\psline(! 24 \AlphaVal 1 \ArcCoords exch pop 0 exch)(! 24 \AlphaVal 1 \ArcCoords exch pop 1 exch)
\TanToArc{6}{2}{-0.5}{0.5}
\psset{linestyle=dashed,dash = 3pt 2pt, linecolor=lightgray}
\parametricplot{0.2}{1.1}{t 1 t sub 8 t dup mul mul 9 sub 2 t mul 3 sub 24 mul div 24 \RSliceCoords}
\psset{linecolor=black,fillstyle=none,linestyle=solid}
\DrawArc
\psset{fillstyle=solid,fillcolor=white,linecolor=black,linestyle=solid,linecolor=black}
\pscircle(!6 2 1 24 \RSliceCoords){0.04}
\rput(!6 2 1 24 \RSliceCoords -0.04 0.05 \Trans){\tiny $\xi$}
\rput(!1 24 \AlphaVal \ArcCoords 0.20 0 \Trans){\tiny $K_{X_{24}}^{\perp}$}
\end{pspicture} \\
\end{tabular} 
&
\begin{tabular}{c}
\begin{pspicture}(-2.0,-0.2)(2.0,0.2)
\psline[arrows=->](1.0,0)(-1.0,0)
\rput(0,0.3){\small $\psi_{3,1}^{*}$} 
\end{pspicture}
\end{tabular}
&
\begin{tabular}{c}
\psset{yunit=3cm, xunit=5cm}
\begin{pspicture}(0,-0.2)(1,1)
\psset{linecolor=gray}
\psline(!1 0 0 \SliceCoords)(!0 1 0 \SliceCoords) 
\psline(! 8 \AlphaVal 1 \ArcCoords exch pop 0 exch)(! 8 \AlphaVal 1 \ArcCoords exch pop 1 exch)
\rput(!1 8 \AlphaVal \ArcCoords 0.6 0 \Trans){\tiny $K_{X_8}^{\perp}$}
\psset{linecolor=black}
\psset{linestyle=dashed,dash = 3pt 2pt, linecolor=lightgray}
\parametricplot{-0.1}{1.1}{t 1 t sub 2 copy mul 2 div 0.125 add 8 \RSliceCoords}
\psset{linestyle=solid}
\psset{fillcolor=vlgray,linecolor=gray,fillstyle=solid}
\psset{linecolor=black,fillstyle=none}
\DrawArc
\psset{fillstyle=solid,fillcolor=white,linecolor=black,linestyle=solid,linecolor=black}
\psset{fillstyle=solid,fillcolor=black}
\pscircle(!0 8 1 8 \RSliceCoords){0.04}
\rput(!0 8 1 8 \RSliceCoords 0.10 0 \Trans){\tiny $C_1$} 
\pscircle(!8 24 7 8 \RSliceCoords){0.04}
\pscircle(!24 48 17 8 \RSliceCoords){0.04}
\pscircle(!48 80 31 8 \RSliceCoords){0.04}
\pscircle(!80 120 49 8 \RSliceCoords){0.04}
\pscircle(!120 168 71 8 \RSliceCoords){0.04}
\pscircle(!168 120 71 8 \RSliceCoords){0.04}
\pscircle(!120 80 49 8 \RSliceCoords){0.04}
\pscircle(!80 48 31 8 \RSliceCoords){0.04}
\pscircle(!48 24 17 8 \RSliceCoords){0.04}
\pscircle(!24 8 7 8 \RSliceCoords){0.04}
\pscircle(!8 0 1 8 \RSliceCoords){0.04}
\rput(!8 0 1 8 \RSliceCoords -0.10 0 \Trans){\tiny $C_{-1}$} 
\psset{fillcolor=white}
\pscircle(!1 1 \ArcCoords){0.04}
\rput(!1 1 \ArcCoords 0.1 add){\tiny $-K_{X_8}$} 
\end{pspicture} 
\end{tabular}  \\
$V_{24}$ & & $V_{8}$ \\
\multicolumn{3}{c}{\Fig\label{fig:illustration-of-theorem}  Illustration of the bound} \\
\end{tabular}

Pulling this back via $\psi_{3,1}$, we obtain a class $\xi$ in $V_{24}$, which is an inner square-zero nef class.
The dashed curve pulls back to a curve (the curve $24C^2=(C\cdot \xi)^2$) which bounds all effective classes.
That curve is tangent to the square-zero cone at $\xi$, and in a neighbourhood of $\xi$ stays very close to
the cone. 

The curve does a worse job of bounding the effective classes farther away from $\xi$.  But, by applying powers of
$T_{24}$ we can shift $\xi$, and the bounding curve, and obtain a family of bounds, which together tightly restrict
the possible $K_X$-positive effective curve classes.

\end{document}